\documentclass[11pt]{amsart}
\usepackage{amsxtra}

\theoremstyle{plain}

\def\CC {{\mathbb C}}
\def\RR {{\mathbb R}}
\def\NN {{\mathbb N}}
\def\ZZ {{\mathbb Z}}
\def\PP {{\mathbb P}}

\def\be {\begin{eqnarray}}
\def\ben {\begin{eqnarray*}}
\def\ee {\end{eqnarray}}
\def\een {\end{eqnarray*}}

\def\AAA{\kern-0.3em}
\def\AA{\kern-0.18em}
\def\AC{\kern-0.14em}
\def\AB{\kern-0.22em}

\newcommand \nc {\newcommand}

\newtheorem{theorem}{Theorem}[section]
\newtheorem{lemma}[theorem]{Lemma}
\newtheorem{proposition}[theorem]{Proposition}
\newtheorem{corollary}[theorem]{Corollary}
\newtheorem{definition}[theorem]{Definition}
\newtheorem{example}[theorem]{Example}
\newtheorem{remark}[theorem]{Remark}
\nc \bth[1] { \begin{theorem}\label{t#1} } \nc \ble[1] {
\begin{lemma}\label{l#1} } \nc \bpr[1] {
\begin{proposition}\label{p#1} } \nc \bco[1] {
\begin{corollary}\label{c#1} } \nc \bde[1] {
\begin{definition}\label{d#1}\rm } \nc \bex[1] {
\begin{example}\label{e#1}\rm } \nc \bre[1] {
\begin{remark}\label{r#1}\rm } \nc \bcon[1] {
\medskip\noindent{\it{Conjecture #1}} } \nc \bqu[1]  {
\medskip\noindent{\it{Question #1}} }
\nc {\ethe} { \end{theorem} }
 \nc {\ele} { \end{lemma} } \nc {\epr}
{ \end{proposition} } \nc {\eco} { \end{corollary} } \nc {\ede} {
\end{definition} } \nc {\eex} { \end{example} } \nc {\ere} {
\end{remark} } \nc {\econ} {\smallskip} \nc {\equ} {\smallskip}
 \nc \thref[1]{Theorem \ref{t#1}}
\nc \leref[1]{Lemma \ref{l#1}} \nc \prref[1]{Proposition
\ref{p#1}} \nc \coref[1]{Corollary \ref{c#1}} \nc
\deref[1]{Definition \ref{d#1}} \nc \exref[1]{Example \ref{e#1}}
\nc \reref[1]{Remark \ref{r#1}}
\def \B {{\mathcal B}}

\def \L {{\mathcal L}}

\def \diag { {\mathrm{diag}} }
\def \res { {\mathrm{Res}} }

 \def\AA  {\kern-0.1em}
 \def\BB  {\kern+0.1em}
 \def\BBB {\kern+0.15em}
 \def\K   {\kern+0.05em}
 \def\MK  {\kern-0.07em}
 \def\MKK {\kern-0.04em}

\begin{document}

\vspace{0.5cm}

\title[ Stokes matrices via monodromy matrices ]
{  Stokes matrices of a reducible   equation with two irregular singularities of Poincar\'{e} rank 1
	via monodromy
	matrices of a reducible Huen type equation}

\author[Tsvetana  Stoyanova]{Tsvetana  Stoyanova}

\date{30.07.2020}

 \maketitle

\begin{center}
{Department of Mathematics and Informatics,
Sofia University,\\ 5 J. Bourchier Blvd., Sofia 1164, Bulgaria, 
cveti@fmi.uni-sofia.bg}
\end{center}

\vspace{0.5cm}

{\bf Abstract.} We consider a second order reducible equation having non-resonant irregular singularities at $x=0$
and $x=\infty$. Both of them are of Poincar\'{e} rank 1. We introduce a small complex parameter $\varepsilon$ that splits
together $x=0$ and $x=\infty$ into four different Fuchsian singularities $x_L=-\sqrt{\varepsilon}, x_R=\sqrt{\varepsilon}$, 
and $x_{LL}=-1/\sqrt{\varepsilon}, x_{RR}=1/\sqrt{\varepsilon}$, respectively. The perturbed equation is a second
order reducible Fuchsian equation with 4 different singularities, i.e. a Heun type equation.
Then we prove that when the perturbed equation has exactly two resonant singularities of different type,
all the Stokes matrices of the initial equation are realized as a limit of the nilpotent parts of the monodromy
matrices of the perturbed equation when $\varepsilon \rightarrow 0$ in the real positive direction.
To establish this result we combine a direct computation with a theoretical approach.



\headsep 10mm \oddsidemargin 0in \evensidemargin 0in

\section{Introduction}

      With this paper we continue our research on the Stokes phenomenon for reducible linear scalar
      equations using perturbative approach.   
	 In this paper we consider the   
		second order reducible equation
		 \be\label{npe}
		  L(y)=L_2(L_1(y))=0
		\ee
  where 
        \be\label{op-ini}
		   L_j=\partial - \left(\frac{\alpha_j}{x} + \frac{\beta_j}{x^2} + \gamma_j\right)\,,\quad j=1, 2\,,\quad
			\partial=\frac{d}{d x}
		\ee
		and $\alpha_j, \beta_j, \gamma_j\in\CC$ are arbitrary bounded parameters.
		In general the equation \eqref{npe} - \eqref{op-ini} has two irregular singular points over $\CC\PP^1$
		of Poincar\'{e} rank 1:\, $x=0$ and $x=\infty$.

		Introducing  a  small complex  parameter $\varepsilon$ we perturb equation \eqref{npe} -\eqref{op-ini} to the equation
	   \be\label{pe-ope}
	     L(y, \varepsilon)=L_{2, \varepsilon}(L_{1, \varepsilon}(y))=0
	    \ee
	    where 
	     	\be\label{pe}
	     	L_{j, \varepsilon}  &=& 
	     	\partial - \frac{\alpha_j}{2}\,\left(
	     	\frac{1}{x-\sqrt{\varepsilon}} + \frac{1}{x+\sqrt{\varepsilon}}\right)-
	     	\frac{\beta_j}{2 \sqrt{\varepsilon}}\,\left(
	     	\frac{1}{x-\sqrt{\varepsilon}} - \frac{1}{x+\sqrt{\varepsilon}}\right)-\\[0.85ex]
	     	&-&
	     	\frac{\gamma_j}{2 \sqrt{\varepsilon}}\,\left(
	     	-\frac{1}{x-\frac{1}{\sqrt{\varepsilon}}} + \frac{1}{x+\frac{1}{\sqrt{\varepsilon}}}\right)\nonumber
	     	\ee
	    such that
		 $$\,
		   L_{j, \varepsilon} \longrightarrow L_j\quad
			\textrm{as}\quad \varepsilon \rightarrow 0\,.
		\,$$
		In general the equation \eqref{pe-ope} - \eqref{pe} is a second order Fuchsian equation having 
		five regular singularities over $\CC\PP^1$  -- four of them are finite and the fifth is at $x=\infty$.
		In this article we restrict our attention to the equation \eqref{pe-ope} - \eqref{pe} which has exactly four
		singular points.
		In particular we study such an equation \eqref{pe} for
		which  the coefficients of the  equation \eqref{npe} - \eqref{op-ini}  do not depend
		on the parameter of perturbation $\varepsilon$. It turns out that there is an unique family of parameters,
		 which satisfies this condition -- the family
		with $\alpha_1=0, \alpha_2=-2$ and arbitrary $\beta_j, \gamma_j$  (see \thref{HE}).  The first order differential 
		operators $L_j$ and $L_{j, \varepsilon}$ which determine this particular family become
		 \be\label{in-par}
		    L_1 = \partial - \left(\frac{\beta_1}{x^2} + \gamma_1\right),\quad
		    L_2 = \partial - \left(-\frac{2}{x} + \frac{\beta_2}{x^2} + \gamma_2\right)
		 \ee
		and
         	 \be\label{pe-par}
         	 L_{1, \varepsilon} &=& \partial - \frac{\beta_1}{2 \sqrt{\varepsilon}}
         	 \left(\frac{1}{x-\sqrt{\varepsilon}} - \frac{1}{x + \sqrt{\varepsilon}}\right)
         	 -\frac{\gamma_1}{2 \sqrt{\varepsilon}}
         	 \left(-\frac{1}{x -\frac{1}{\sqrt{\varepsilon}}} + \frac{1}{x + \frac{1}{\sqrt{\varepsilon}}}\right)\,,\\[0.2ex]
         		 L_{2, \varepsilon} &=& \partial 
             		 - \left(\frac{\beta_2}{2 \sqrt{\varepsilon}} - 1\right)
         		 \frac{1}{x-\sqrt{\varepsilon}}
         			 + \left(\frac{\beta_2}{2 \sqrt{\varepsilon}} + 1\right)
         			 \frac{1}{x+\sqrt{\varepsilon}}-\nonumber\\[0.2ex]
         			 &-& 
         		 \frac{\gamma_2}{2 \sqrt{\varepsilon}}
         		 \left(-\frac{1}{x -\frac{1}{\sqrt{\varepsilon}}} + \frac{1}{x + \frac{1}{\sqrt{\varepsilon}}}\right)\nonumber
         	 \ee		
          Through this paper we call the equation \eqref{npe} with \eqref{in-par} the initial equation.
          Since every four distinct points over the Riemann sphere can be fixed by a M\"{o}bius transformation at
          $0, 1, \infty$ and $a\neq 0, 1, \infty$ through this paper we call the  equation \eqref{pe-ope} with \eqref{pe-par} 
           the perturbed equation or the Heun type equation. In this paper we assume that the point $x=0$ is a 
           non-resonant irregular point for the initial equation. This assumption is equivalent to the condition
            \be\label{nonr}
               \beta_1 \neq \beta_2\,.
            \ee
           
           		To each singular point of the initial equation, with respect to a given fundamental
           		matrix, we associate the so called Stokes matrices $St^{\theta}_j, j=0, \infty$ corresponding to the singular 
           		direction $\theta$. Due to the reducibility the initial equation  admits actual
           		fundamental matrices $\Phi_0(x, 0)$ at $x=0$ and $\Phi_{\infty}(x, 0)$ at $x=\infty$ with respect to which
           		the Stokes matrices $St^{\theta}_j$ have the upper-triangular form
           		$$\,
           		St^{\theta}_j=\left(\begin{array}{cc}
           		1   &\mu^{\theta}_j\\
           		0   &1
           		\end{array}
           		\right)\,.     
           		\,$$ 
           		Moreover the reducibility  implies that the initial equation can have only one singular direction
           		$\theta=\arg (\beta_1-\beta_2)$ at the origin (see \thref{A-S}) and only one singular direction 
           		$\theta=\arg (\gamma_2-\gamma_1)$ at $x=\infty$ (see \thref{A-S-inft}).
           Denote by $x_R=\sqrt{\varepsilon}, x_L=-\sqrt{\varepsilon}, x_{RR}=1/\sqrt{\varepsilon},
           x_{LL}=-1/\sqrt{\varepsilon}$ the singular points of the perturbed equation. To each of them we associate the so
           called monodromy matrices $M_j(\varepsilon)$ and $M_{jj}(\varepsilon), j=L, R$.
           Through this paper we call a double resonace these values of the parameters $\beta_j, \gamma_j, \varepsilon$
           for which the solution of the perturbed  equation can contain logarithmic terms near exactly two
           singular points of different type. Under this we mean that one of the point is of the type $x_j, j=R, L$
           and the other is of the type $x_{jj}, j=R, L$. It turns out that the perturbed equation has a double resonance if and only
           if $\varepsilon\in\RR^{+}$ (see \prref{two}). Due to the reducibility diring a double resonance the Heun type equation admits two different fundamental matrices $\Phi_0(x, \varepsilon)$ and $\Phi_{\infty}(x, \varepsilon)$
           which respect to which the monodromy matrices $M_j(\varepsilon)$ and $M_{jj}(\varepsilon)$ have
           the simple form (see \thref{M})
            \be\label{dec}
             M_j(\varepsilon)=e^{\pi\,i (\Lambda + \frac{1}{x_j}\,B)}\,e^{2 \pi\,i\,T_j},\quad
                M_{jj}(\varepsilon)=e^{-\pi\,i\,x_{jj}\,G}\,e^{2 \pi\,i\,T_{jj}},\quad
                j=R, L\,.
            \ee
            The constant diagonal matrices $\Lambda, B$ and $G$ come from the initial equation (see \prref{p2}). The matrices
            $\Lambda/2 + B/2 x_j$ and $-x_{jj} G/2$ are their perturbed analog. The matrices $T_j$ and $T_{jj}$ are nilponent
             $$\,
               T_j=\left(\begin{array}{cc}
                  0  &d_j\\
                  0  &0
                          \end{array}
                     \right),\qquad
                    T_{jj}=\left(\begin{array}{cc}
                    0  &d_{jj}\\
                    0  &0
                    \end{array}
                    \right),\quad j=L, R         
             \,$$
             and they measure the existence of logarithmic terms around the singular points.
          Then the main result of this paper demonstrates in an explicit way that during a double resonance both Stokes matrices $St^{\theta}_0$ and $St^{\theta}_{\infty}$ of the initial equation are realized as a limit of the matrices
          $e^{2 \pi\,i\,T_j}$ and $e^{2 \pi\,i\,T_{jj}}$, i.e.
           $$\,
            e^{2 \pi\,i\,T_j} \longrightarrow St^{\theta}_0,\quad
            e^{2 \pi\,i\,T_{jj}} \longrightarrow St^{\theta}_{\infty},\quad j=R, L
           \,$$
        when $\sqrt{\varepsilon} \rightarrow 0$.

           Glutsyuk make the start for studying the Stokes phenomenon of linear systems with an  irregular singularity at the
           origin by a perturbative approach with the papers \cite{AG, AG1, AG2}. More precisely, he considers a generic perturbation 
           depending on a parameter $\varepsilon$ that splits the irregular singularity at the origin of the initial system into
           Fuchsian singularities. Then Glutsyuk  prove that the Stokes operators of the initial
           system with an irregular singularity at the origin of Poincar\'{e} rank $k \geq 1$ are limits of the 
           transition operators of the perturbed  system. Recently, Lambert, Rousseau, Hurtubise and Klime\u{s}
           \cite{CL-CR, CL-CR1, HLR, Kl-2} using a different approach extend the set of the values of the parameter of
           perturbation $\varepsilon$ over a whole neighborhood of $\varepsilon=0$. They  study only non-resonant initial system
           with an irregular singularity at the origin of Poincar\'{e} rank $k\in\NN$. In the case when $k=1$ the perturbation 
           splits the origin into two Fuchsina singularity $x_L=\sqrt{\varepsilon}$ and $x_R=-\sqrt{\varepsilon}$ \cite{CL-CR}.
           In \cite{CL-CR} Lambert and Rousseau prove
           that the monodromy operator acting on the fundamental matrix solutions of the 
           perturbed systems decomposes into the Stokes operator multiplied by the classical monodromy operator acting
           on the branch of $(x-x_L)^{\Lambda/2 + B/2 x_L}\,(x-x_R)^{\Lambda/2 +B/2 x_R}$. In addition, in \cite{CL-CR1}
           they prove that the so called unfolded Stokes matrices $St_j(\varepsilon), j=L, R$ of the perturbed system 
           depend analytically on the 
           parameter of perturbation $\varepsilon$ and converge when $\varepsilon \rightarrow 0$ to the Stokes matrices
           $St_j, j=L, R$ of the initial system.
            Later in \cite{Kl-2} Klime\u{s}
           specifies this result expressing the acting of the monodromy operators on the solutions of the perturbed systems
           by the monodromy matrices $M_j(\varepsilon)$, the unfolded Stokes matrices $St_j(\varepsilon)$ and the
           matrices $e^{\pi\,i (\Lambda + B/x_j)}, j=L, R$. In \cite{Re} Remy also is interested in Stokes phenomenon for
           a linear differential system with an irregular singularity at the origin $x^{r+1}\,Y'(x)=A(x)\,Y(x)$ and
           an arbitrary single level $r \geq 1$ from a perturbative point of view. But he does not split the origin into
           Fuchsian singularities. He considers a regular holomorphic perturbation of the coefficients of the initial system
           $x^{r+1}\,Y'(x)=A^{\varepsilon}(x)\,Y(x),\,A^{\bf 1}(x)=A(x)$ which preserves the single level $r$ of the initial system.
           He proves that the Stokes-Ramis matrices of the initial system are limits of convenient products of the
           Stokes-Ramis matrices of the perturbed system. 
           
            Parallel to the study of linear systems Ramis \cite{R},
           Zhang \cite{Z}, Duval \cite{Du} and latter Lambert and Rousseau \cite{CL-CR2} dealt with the Stokes phenomenon
           in the confluence of the classical hypergeometric equation and generalized hypergeometric family.
           In \cite{Du} Duval study the family $D_{p+1, p}(\underline{\alpha}; \underline{\beta})$ of generalized hypergeometric
           equations for $p \geq 2$. This equation has two Stokes matrices $St_0$ and $St_{\pi}$ and can be obtained
           from the Fuchsian equation $D_{p+1, p+1}(\underline{\alpha}; \underline{\beta})$ by a confluence procedure. 
           More precisely, Duval  apply  such a confluence procedure by the change $z=t/b$ and then making 
           $b \rightarrow \infty$. The change $z=t/b$ takes the singular points $z=0, 1, \infty$ of the equation 
           $D_{p+1, p+1}(\underline{\alpha}; \underline{\beta})$ into  the points $0, b, \infty$.
           When $b \rightarrow \infty$ in a non-real direction Duval prove, by a direct calculation.
            that the Stokes matrices can be obtained as limit
           of the connection matrices linking well chosen fundamental set of solutions around $b$ and $\infty$ of the
           Fuchsian equation. When $b \rightarrow \infty$ in a real direction she express, by a direct calculation, 
           the Stokes matrices as limits of the monodromy matrices around $b$ and $\infty$ with respect to a "mixed"
           basis of solutions. 
           In \cite{St1} we consider a particular family of third
           order linear reducible scalar equation having a non-resonant irregular singularity at the origin of 
           Poincar\'{e} rank 1. We prove that the Stokes matrices of the initial equation are limits of the 
           nilpotent parts of the monodromy matrices of a Fuchsian equation with respect to a "mixed" basis of solutions  
           that contains logarithmic terms. This Fuchsian equation was obtained from the initial equation by introducing a 
           small real parameter
           $\varepsilon$ that splits the origin into two regular singularities $x_L=-\sqrt{\varepsilon}$ and
           $x_R=\sqrt{\varepsilon}$. 
           
          Very recently in \cite{M} Malek is interested in  the effect of the unfolding in families of singular PDEs from 
          an asymptotic point of view. He unfolds 
          singularly perturbed differential operator of irregular type $\varepsilon\,t^2\,\partial_t$ into a family
          of singular operators 
          $\mathcal{D}_{\varepsilon, \alpha}(\partial_t)=(\varepsilon\,t^2 - \varepsilon^{\alpha})\,\partial_t$
          of Fuchsian type. Then he studies how this unfolding changes the asymptotic properties of holomorphic
          solutions of the unfolded equation in comparison to the ones of the initial equations.

         In this paper, as in our previous two papers \cite{St1, St2} we deal with reducible scalar equations.
         This time the initial equation has two non-resonant irregular singular points of Poincar\'{e} rank 1, taken at
         $x=0$ and $x=\infty$. Introducing a small complex parameter of perturbation $\varepsilon$ we split 
         the non-resonant singularity at the origin
         into two finite Fuchsian singularities $x_R=\sqrt{\varepsilon}$ and $x_L=-\sqrt{\varepsilon}$. At the same
         time we also split $x=\infty$ into again two finite Fuchsian singularities $x_{RR}=1/\sqrt{\varepsilon}$
         and $x_{LL}=-1/\sqrt{\varepsilon}$. Since this situation is more complicated we restrict our attention
         to second order equations. Due to the reducibility both initial and perturbed equations admit 
         upper-triangular fundamental matrices, whose off-diagonal element can be expressed as an integral.
         As in \cite{St1} we fully exploit this phenomenon in purpose of explicitly computing the Stokes matrices
         of the initial equation and the monodromy matrices of the perturbed equation. In fact, in this paper we use 
         two different fundamental matrices $\Phi_0(x, 0)$ and $\Phi_{\infty}(x, 0)$ of the initial equation
         depending of the path of integration. These solutions correspond to the point $x=0$ and $x=\infty$. Here
         we are not interested in the connection problem. The off-diagonal element of the actual $\Phi_0(x, 0)$ and 
         $\Phi_{\infty}(x, 0)$ is obtained from the integral representation by two steps. In the first step we expres the
         integral as formal power series in $x$ (resp. $x^{-1}$). It turns out that these series are in general divergent.
         In the second step, utilizing the summability theory, we lift these formal series to actual solutions.  
         With respect to the actual $\Phi_0(x, 0)$ and $\Phi_{\infty}(x, 0)$
         we explicitly compute the corresonding Stokes matrices $St^{\theta}_0$ and $St^{\theta}_{\infty}$, where
         $\theta$  is a singular direction. To each fundamental matrix of the initial equation corresonds a
         fundamental matrix $\Phi_j(x, \varepsilon), j=0, \infty$ of the perturbed equation, such that
         $\lim_{\sqrt{\varepsilon} \rightarrow 0} \Phi_j(x, \varepsilon)=\Phi_j(x, 0), j=0, \infty$. Moreover, this
         choice of fundamental matrices leads to the so called "mixed" basis of solutions for which one of the
         solutions is always a eigenvector of the local monodromy around the singular points (see \reref{mixed}). With respect
         to exactly these fundamental matrices we explicitly compute the monodromy matrices $M_j(\varepsilon)$
         and $M_{jj}(\varepsilon), j=R, L$ of the perturbed equation only during a double resonance. Then
         we explicitly show that during a double resonance the both Stokes matrices of the initial equation are
         realized as a limit of the nilpotent parts of the suitable monodromy matrices of the perturbed equation
        when $\sqrt{\varepsilon} \rightarrow 0$ (see \thref{main}). 
        Since we compute the Stokes and monodromy matrices explicitly by hand one can say that our approach is closer
        to the works of Duval \cite{Du}, Ramis \cite{R} and Zhang \cite{Z}. Moreover, our result almost repeats the result
        of Duval when $b$ (resp. $\sqrt{\varepsilon}$) tends to $\infty$ in the real direction. On the other hand combining
        the decomposition results of Klime\u{s} (see \prref{un}) and our formulas \eqref{dec} we find that the unfolded
        Stokes matrices $St_j(\varepsilon)$ are equal to the matrices $e^{2 \pi\,i\,T_j}, j=R, L$ ( see \coref{T}).
        Moreover, we extend this result to a neigborhood of $x=\infty$ and get the same equivalence between the 
        unfolded Stokes matrices $St_{jj}(\varepsilon)$ and the matrices $e^{2 \pi\,i\,T_{jj}}, j=R, L$ (see \coref{T-inf}).
       Then using that the unfolded Stokes matrices converge analytically to the Stokes matrices we establish the main 
       result of this paper. So, we can say that our approach is new cleverly weaves the direct computation
       (Duva, Ramis, Zhang) into the theoretical findings (Glutsyuk, Hurtubise, Klime\u{s}, Lambert, Rousseau).
       Moreover our approach allows us to split more that one irregular singularity into several Fuchsian singularities.

      In addition, we find two particular families of parameters $\beta_j, \gamma_j$ for which the initial equation
      has trivial Stokes matrices at both singular points. The first family is determined by the condition
       $$\,
        \gamma_1=\gamma_2 
       \,$$
       and $\beta_j$'s are arbitrary such that satisfy the non-resonant condition \eqref{nonr}.
      This family is expected. For these values of the parameters the elements of the matrices $\Phi_0(x, 0)$ and
      $\Phi_{\infty}(x, 0)$ are defined by well known  functions 
      (see \prref{p2}(3) and \prref{inf}(1). It turns out that the matrices $e^{2 \pi\,i\,T_j}$ and $e^{2 \pi\,i\,T_{jj}},\,j=R, L$
      are also trivial when $\gamma_1=\gamma_2$ (see from \thref{A.1-M}(1) to \thref{A.4-M-inf}(1)).
      The second family of parameters for which the both Stokes matrices are trivial is a very interesting. 
      We are surprised to find that for parameters $\beta_j, \gamma_j$ such that satisfy the condition \eqref{nonr} and
       \be\label{sur}
         \sum_{n=0}^{\infty}
         \frac{(-1)^{n+1}\,(\gamma_2-\gamma_1)^n\,(\beta_2-\beta_1)^n}{n!\,(n+1)!}=0
       \ee
      the initial equation also has trivial Stokes matrices at both singular points (see \thref{A-S} and \thref{A-S-inft}).
      The cause for this phenomenon is the convergence of the series
       $$\,
       \hat{\psi}(x)=\sum_{k=1}^{\infty}
       \frac{(-1)^k\,k!\,S_{k-1}}{(\beta_2-\beta_1)^k}\,x^k\quad
       \textrm{and}\quad
         \hat{\varphi}(x)=\sum_{k=1}^{\infty}
         \frac{\,k!\,S_{k-1}}{(\gamma_2-\gamma_1)^k}\,x^{-k}
       \,$$
       under condition \eqref{sur}. Here $S_{k-1}, k\in\NN$ are the partial sums of the number series \eqref{sur}.
       The series $\hat{\psi}(x)$ is an element of the formal matrix solution at $x=0$, and the series $\hat{\varphi}(x)$
       includes in the formal matrix solution at $x=\infty$ (see \prref{p2}(2) and \prref{inf}(2.b)). Unfortunately, till now we 
       can not point which analytic functions
       in $\CC$ and $\CC\PP^1-\{0\}$, respectively, these series represent at $x=0$ and $x=\infty$, respectively.
       Moreover, till now we can not say if the corresponding matrices $e^{2 \pi\,i\,T_j}$ and $e^{2 \pi\,i\,T_{jj}}, j=R, L$
       are trivial or not under condition \eqref{sur}. However their limits when $\sqrt{\varepsilon} \rightarrow 0$ tend
        to the identity matrices provided that the condition \eqref{sur} holds.

      On of the motivation for studying this pair of equations arises from its relation with
       Heun class of equations. We will only comment the application of the initial equation regarding it as a particular 
       double Heun equation. More precisely we connect some particular families of the initial equation with two important
       non-linear ordinary differential equations. For classes of solutions, properties and applications of the Heun class of
       equations see the book of Slavyanov and Lay \cite{S-L}.  Following \cite{S-S} the double confluent Heun equation (DCHE) is the 
       second order linear differential equation
         \be\label{dche}
          x^2\,y''(x) + (-x^2 + c x  - t)\,y'(x) + (-a x + t\,h)\,y(x)=0\,,
        \ee
       where $c, a, t, h$ are complex parameters. The DCHE admits two singular points over $\CC\PP^1$: $x=0$ and $x=\infty$.
       Both of them are irregular points. The DCHE has many applications in the mathematical physics: optics, hydrodynamics, gravity.   Comparing the coefficients $b_1(x, 0)$ and $b_2(x, 0)$ from section 2.5 we see that
       the initial equation is a DCHE \eqref{dche} if and only if the parameters $c, \beta_j, \gamma_j$ satisfy the following conditions
         \be\label{p-dche}
            & &
         \gamma_1+\gamma_2=1,\quad \beta_1+\beta_2=t,\quad c=2,\quad \beta_1 \beta_2=0,\\
            & &
         \gamma_1 \gamma_2=0,\quad \gamma_1=\frac{a}{2},\quad
         \beta_1 \gamma_2 + \beta_2 \gamma_1=-t\,h\,.\nonumber
        \ee
        In the same paper \cite{S-S}
       Salatich and Slavyanov relate by an antiquantization procedure the DCHE \eqref{dche} with the third Painlev\'{e} equation
       $P_{III}$
        $$\,
          \ddot{q} - \frac{\dot{q}^2}{q} + \frac{\dot{q}}{t} - \gamma\,q^3 
          -\frac{\alpha q^2 + \beta}{t} - \frac{\delta}{q}=0\,.
        \,$$
         In section 2.2 repeating their approach we obtain two particular families of $P_{III}$ equation from the particular
         families of the initial equation that satisfy the conditions \eqref{p-dche} provided that the non-resonant condition
         \eqref{nonr} remains valid.    
       
           Recently in several works \cite{B-G, B-G-1, B-T} Buchstaber, Glutsyuk and Tertychnyi study the following family of DCHE
          \be\label{dche-1}
            z^2\,E'' + (n\,z + \mu (1-z^2))\,E' + (\lambda - \mu\,n\,z)\,E=0,\quad
            n, \lambda, \mu\in\CC,\quad \mu \neq 0\,.
          \ee        
        When the parameters are real such that $\lambda + \mu^2 > 0$ the equation \eqref{dche-1} appears as a linearization
        of the family of nonlinear equations on two-torus that model the Josephson effect in superconductively.
        The initial equation coincides with the equation \eqref{dche-1} if and only if
         $$\,
           n=2,\,\beta_1+\beta_2=-\mu,\,\gamma_1+\gamma_2=\mu,\,\beta_1 \beta_2=\gamma_1 \gamma_2=0,\,\gamma_1=\mu,\,
           \beta_1 \gamma_2 + \beta_2 \gamma_1 =\lambda\,.
         \,$$
       In particular, following \cite{B-T} the initial equation
         $$\,
           x^2\,y''(x) + (-2 x + \mu - \mu x^2)\,y'(x) - 2 \mu\,x\,y(x)=0,\quad
           \mu \neq 0
         \,$$
       appears as a linearization of the equation
        \be\label{Jo}
         \dot{\varphi}(t) + \sin \varphi(t)=-3 \omega + cos \omega t\,,
        \ee
        where $\omega > 0$ is a real parameter, related to $\mu$ by $\mu=1/2 \omega$. The equation {Jo} is a particular family of the 
        non-linear equation
         $$\,
          \dot{\varphi}(t) + \sin \varphi(t)=B + A \,\cos \omega t,\quad
          A, \omega > 0,\quad B \geq 0\,.
         \,$$
         Here $\omega > 0$ is a fixed constant, $A, B$ are parameters. The last equation arises in several models in physics
         (Josephson junction in superconductively), mechanics, geometry.
        Among many results, in the pointed papers the authors study the problems of the existence of a holomorphic solution
        on $\CC$ of \eqref{dche-1} and the existence of the eigenfunctions of the monodromy operators with a given eigenvalue. 
        
        This article is organized as follows. In the next section we recall the required facts and definitions from the theory
        of differential equations with irregular singularities, as well as from the theory of Fuchsian equations. We also introduce
        the perturbed equation as a Heun type equation and define global fundamental matrices of both equations, which we use
        to build local fundamental matrix solutions. 
        In section 3 we compute the Stokes matrices of the initial equation under the restriction \eqref{nonr}. In section 4
        we compute the monodromy matrices of the perturbed equation provided that condition \eqref{nonr} is valid.
        In the last section 5 we establish the main results of this paper.
       
%%%%%%%%%%%%%%%%%%%%%%%%%%%%%%%%%%%%%%%%%%%%%%%%%
%theory
%%%%%%%%%%%%%%%%%%%%%%%%%%%%%%%%%%%%%%%%%%%%%%%%

   \section{Preliminaries}

    \subsection {The  equation \eqref{pe-ope} - \eqref{pe} as a Heun type equation}

      In this section we introduce the perturbed equation as a Heun type equation.
      
     In general the  equation \eqref{pe-ope} - \eqref{pe} is a second order Fuchsian equation
		with 5 regular singularities over $\CC\PP^1$ taken at $x_1=\infty,\,
		x_2=\sqrt{\varepsilon},\,x_3=-\sqrt{\varepsilon},\,x_4=1/\sqrt{\varepsilon},\,
		x_5=-1/\sqrt{\varepsilon}$. The next theorem describes all the cases when the 
		this equation has only four regular singular points.
		
		\bth{HE}
		 The  equation \eqref{pe-ope} - \eqref{pe} has only four regular singular points over $\CC\PP^1$
		if and only if the parameters $\alpha_j,\,\beta_j,\,\gamma_j$ and $\varepsilon$
		satisfy one of the following conditions:
		  \begin{itemize}
			  \item[(I.)]\,
			     	$\,
			     	\alpha_1=0,\,\alpha_2=-2,
			     	\,$ and
			     	$\beta_j$ and $\gamma_j$ are arbitrary.
				\item[(II.)]\,
					$\,
					\alpha_1=\alpha_2=-1,\, \beta_1-\beta_2 + \frac{\gamma_2-\gamma_1}{\varepsilon}=0\,.
					\,$
				\item[(III.)]\,
				 $\gamma_1=-\gamma_2=-2 \sqrt{\varepsilon}\,,\quad
			\frac{\varepsilon(\beta_1 - \beta_2) + \sqrt{\varepsilon}(\alpha_1-\alpha_2)}{\varepsilon (1-\varepsilon^2)}=1$\,.
		     	\item[(IV.)]\,
				 $\gamma_2=-\gamma_1=-2 \sqrt{\varepsilon}\,,\quad
			\frac{\varepsilon(\beta_1 - \beta_2) - \sqrt{\varepsilon}(\alpha_1-\alpha_2)}{\varepsilon (1-\varepsilon^2)}=-1$\,.
		      \item[(V)]\,
					 $\beta_1 = - \alpha_1\,\sqrt{\varepsilon},\,
					\beta_2 = - \alpha_2\,\sqrt{\varepsilon}$, and
					 $\gamma_j$ are arbitrary.
					\item[(VI.)]\,
					  $\beta_1=-(\alpha_1-2)\sqrt{\varepsilon},\,
				\beta_2=-(\alpha_2+2)\sqrt{\varepsilon}$, and
				$\frac{\alpha_2-\alpha_1}{2} + \frac{\sqrt{\varepsilon}}{1-\varepsilon^2}
				(\gamma_2-\gamma_1)=-1$.
				   \item[(VII.)]\,
					  $\beta_1=\alpha_1\sqrt{\varepsilon},\,
			 \beta_2=\alpha_2\sqrt{\varepsilon}$ and $\gamma_j$ are arbitrary.
		        \item[(VIII.)]\,
						  $\beta_1=(\alpha_1-2)\sqrt{\varepsilon},\,
					\beta_2=(\alpha_2+2)\sqrt{\varepsilon}$, and
					$\frac{\alpha_2-\alpha_1}{2} - \frac{\sqrt{\varepsilon}(\gamma_2-\gamma_1)}{1-\varepsilon^2}=-1$.
		\end{itemize}
		\ethe

    \proof
		 
		The transformation $x=1/t$ changes the  equation  \eqref{pe-ope} - \eqref{pe} into the equation
		 \be\label{ipe}
		  \ddot{y}(t) &+&
			\left[\frac{p_1}{t} + \frac{p_2}{t-\sqrt{\varepsilon}} + \frac{p_3}{t+\sqrt{\varepsilon}}
			+\frac{p_4}{t-\frac{1}{\sqrt{\varepsilon}}} + \frac{p_5}{t+\frac{1}{\sqrt{\varepsilon}}}\right]\,
			\dot{y}(t) +\nonumber\\[0.3ex]
			            &+&
			\Big[\frac{q_{10}}{t^2} + \frac{q_{11}}{t} + \frac{q_{20}}{(t-\sqrt{\varepsilon})^2}
			+ \frac{q_{21}}{t-\sqrt{\varepsilon}} + \frac{q_{30}}{(t+\sqrt{\varepsilon})^2}
			+ \frac{q_{31}}{t+\sqrt{\varepsilon}} + \\[0.3ex]
			           &+&
				\frac{q_{40}}{(t-\frac{1}{\sqrt{\varepsilon}})^2} + \frac{q_{41}}{t-\frac{1}{\sqrt{\varepsilon}}} +
				\frac{q_{50}}{(t+\frac{1}{\sqrt{\varepsilon}})^2}	+ \frac{q_{51}}{t+\frac{1}{\sqrt{\varepsilon}}}\Big]\,
								y(t)=0\,,\nonumber		
		\ee
 where the coefficients $p_j,\,q_{j0},\,q_{j1}$ depend on the parameters $\alpha_i,\,\beta_i, \gamma_i$ and
$\varepsilon$. The equation \eqref{ipe} has 5 regular singularities over $\CC\PP^1$ taken at
$t_1=0, \,t_2=\sqrt{\varepsilon},\,t_3=-\sqrt{\varepsilon},\,t_4=1/\sqrt{\varepsilon},\,t_5=-1/\sqrt{\varepsilon}$. 
 Each point $t=t_j,\, 1 \leq j \leq 5$ is an ordinary point for the equation \eqref{ipe} if and only if
the coefficients $p_j,\,q_{j0}$ and $q_{j1}$ satisfy the simultaneous conditions
   \ben
	   \left|\begin{array}{l}
		  p_j=0\,,\\
			q_{j0}=0\,,\\
			q_{j1}=0\,.
			     \end{array}
			\right.		
	\een
	
In particular, the point $t_1=0$ is an ordinary point for the equation \eqref{ipe} if and only if
  \ben
	  \left|\begin{array}{l}
		p_1=\alpha_1 + \alpha_2 +2=0\,,\\
	   q_{10}=\alpha_1\,(\alpha_2+1)=0\,,\\
			q_{11}=2 \beta_1 + \alpha_1\,\beta_2 + \alpha_2\,\beta_1 -
			\frac{2 \gamma_1}{\varepsilon} - \frac{\alpha_1\,\gamma_2}{\varepsilon} -
			\frac{\alpha_2\,\gamma_1}{\varepsilon}=0\,.
		      \end{array}
		\right.			
	\een
  This system of conditions has two families of solutions. The first one is given by
	$\alpha_1=0,\,\alpha_2=-2$ and the parameters $\beta_j,\,\gamma_j$ and $\varepsilon$
	are arbitrary. This family defines the case $(I.)$.
	The second family of solutions
	is defined by $\alpha_1=\alpha_2=-1$ and
	 $$\,
	  \beta_1-\beta_2 + \frac{\gamma_2-\gamma_1}{\varepsilon}=0\,,
	\,$$
	which we call the case $(II.)$.
	In the cases $(I.)$ and $(II.)$  the equation  \eqref{pe-ope} - \eqref{pe} 
	admits only finite regular singular points taken at $x_2, x_3, x_4$ and $x_5$. 
	
   Next, the point $t_2=\sqrt{\varepsilon}$ is an ordinary point for the equation \eqref{ipe} if and only if
		\ben
		 \left|\begin{array}{l}
		   p_2=\frac{\gamma_1+\gamma_2}{2 \sqrt{\varepsilon}}=0\,,\\[0.35ex]
			 q_{20}=-\frac{\gamma_1}{4} (2 \sqrt{\varepsilon} - \gamma_2)=0\,,\\[0.35ex]
			 q_{21}=\frac{\gamma_1}{2 \varepsilon\,\sqrt{\varepsilon}} (2 \sqrt{\varepsilon} - \gamma_2)
			+\frac{\gamma_1\,\gamma_2}{4 \varepsilon}+
			\frac{\beta_1\,\gamma_2\,\sqrt{\varepsilon}+\beta_2\,\gamma_1\,\sqrt{\varepsilon}
			+\alpha_1\,\gamma_2 + \alpha_2\,\gamma_1}{2 \varepsilon (1-\varepsilon^2)}=0\,.
			    \end{array}
			\right.		
		\een
		This conditions imply that either
		 $$\,
		  \gamma_1=\gamma_2=0,\, \quad \textrm{and}\quad \alpha_j, \beta_j,\,\varepsilon \quad 
			\textrm{are arbitrary}\,,
		\,$$
	 or
	  $$\,
		  \gamma_1=-\gamma_2=-2 \sqrt{\varepsilon}\,,\quad
			\frac{\varepsilon(\beta_1 -\beta_2) + \sqrt{\varepsilon}(\alpha_1-\alpha_2)}{\varepsilon (1-\varepsilon^2)}=1\,.
		\,$$
		In the firs case when $\gamma_1=\gamma_2=0$ the equation  \eqref{pe-ope} - \eqref{pe} admits 3 regular singularities taken
		at $x_1, x_2$ and $x_3$. Because of this reason we do not include it in the list of the statement
		of the theorem. Note also that in this case the point $x=\infty$ is a regular point for the  equation \eqref{npe} - \eqref{op-ini}.
    	We call the second  case $(III)$. In this  case the  equation  \eqref{pe-ope} - \eqref{pe} 
		admits regular singular points taken at $x_1, x_2, x_3$ and $x_5$.

		Similarly, the point $t_3=-\sqrt{\varepsilon}$ is an ordinary point for the equation \eqref{ipe} if and only if
		\ben
		  \left|\begin{array}{l}
			 p_3=-\frac{\gamma_1+\gamma_2}{2 \sqrt{\varepsilon}}=0\,,\\[0.35ex]
			 q_{30}=\frac{\gamma_1}{4}(2 \sqrt{\varepsilon}+\gamma_2)=0\,,\\[0.35ex]
			 q_{31}=\frac{\gamma_1}{2 \varepsilon \sqrt{\varepsilon}}(2 \sqrt{\varepsilon}+\gamma_2)-
			\frac{\gamma_1 \gamma_2}{4 \varepsilon} -
			\frac{\beta_1\,\gamma_2\,\sqrt{\varepsilon}+\beta_2\,\gamma_1\,\sqrt{\varepsilon}
			-\alpha_1\,\gamma_2-\alpha_2\,\gamma_1}{2 \varepsilon (1-\varepsilon^2)}=0\,.
			     \end{array}
			\right.		
		\een
		This system of conditions  implies either
		$$\,
		   \gamma_1=\gamma_2=0,\,\quad \textrm{and}\quad
			\alpha_j,\,\beta_j,\,\varepsilon \quad \textrm{are arbitrary}\,
		\,$$
		or
		$$\,
		  \gamma_2=-\gamma_1=-2 \sqrt{\varepsilon}\,,\quad
			\frac{\varepsilon(\beta_1 - \beta_2) -\sqrt{\varepsilon} (\alpha_1 -\alpha_2)}{\varepsilon (1-\varepsilon^2)}=-1\,.
		\,$$
		Just before  the first case when $\gamma_1=\gamma_2=0$ and $\alpha_j,\,\beta_j$ are
		arbitrary non-zero parameters is out of the statement of the theorem. We call the second case 
		the case $(IV.)$. 
		
		The point $t_4=1/\sqrt{\varepsilon}$ is an ordinary point for the equation \eqref{ipe} if and only if
		\ben
		  \left|\begin{array}{l}
			 p_4=-\left(\frac{\alpha_1+\alpha_2}{2}+\frac{\beta_1+\beta_2}{2 \sqrt{\varepsilon}}\right)=0\,,\\[0.35ex]
			  q_{40}=\left(\frac{\alpha_1}{2} + \frac{\beta_1}{2 \sqrt{\varepsilon}}\right)
			 \left(1 +\frac{\alpha_2}{2} + \frac{\beta_2}{2 \sqrt{\varepsilon}}\right)=0\,,\\[0.45ex]
				q_{41}=-2 \sqrt{\varepsilon}
				 \left(\frac{\alpha_1}{2} + \frac{\beta_1}{2 \sqrt{\varepsilon}}\right)
				  \left(1 +\frac{\alpha_2}{2} + \frac{\beta_2}{2 \sqrt{\varepsilon}}\right)
				+ \frac{\beta_1 \beta_2 - \varepsilon \alpha_1 \alpha_1}{4 \sqrt{\varepsilon}}
				 - \frac{\varepsilon (\alpha_1 \gamma_2 + \alpha_2 \gamma_1)}{2 (1-\varepsilon^2)}
				  - \frac{\varepsilon (\beta_1 \gamma_2 + \beta_2 \gamma_1)}{2 \sqrt{\varepsilon }(1-\varepsilon^2)}=0\,.
	   \end{array}
			\right.			
		\een
     This system admits two families of solutions. The first one is given by
				$$\,
				  \beta_1 = - \alpha_1\,\sqrt{\varepsilon},\,
					\beta_2 = - \alpha_2\,\sqrt{\varepsilon}, \quad \textrm{and}\quad
					 \gamma_j,\,\varepsilon \quad \textrm{are arbitrary}\,.
				\,$$
		The second one is defined by
		  $$\,
			  \beta_1=-(\alpha_1-2)\sqrt{\varepsilon},\,
				\beta_2=-(\alpha_2+2)\sqrt{\varepsilon},\,\,\,\,
				\frac{\alpha_2-\alpha_1}{2} + \frac{\sqrt{\varepsilon}}{1-\varepsilon^2}
				(\gamma_2-\gamma_1)=-1\,.
			\,$$
			We call these two cases the case $(V.)$ and the case $(VI.)$, respectively. In both cases the  equation
		 \eqref{pe-ope} - \eqref{pe} 	admits regular singular points at $x_1, x_3, x_4$ and $x_5$.
			
		The last singular point $t_5=-1/\sqrt{\varepsilon}$ is an ordinary point  if and only if
		 \ben
		  \left|\begin{array}{l}
		    	 p_5=-\left(\frac{\alpha_1+\alpha_2}{2}-\frac{\beta_1+\beta_2}{2 \sqrt{\varepsilon}}\right)=0\,,\\[0.35ex]
		    	 q_{50}=\left(\frac{\alpha_1}{2} - \frac{\beta_1}{2 \sqrt{\varepsilon}}\right)
		    	 	  \left(1 +\frac{\alpha_2}{2} - \frac{\beta_2}{2 \sqrt{\varepsilon}}\right)=0\,,\\[0.45ex]
		    	 q_{51}=-2 \sqrt{\varepsilon}
		    	 	  \left(\frac{\alpha_1}{2} - \frac{\beta_1}{2 \sqrt{\varepsilon}}\right)
		    	 	  \left(1 +\frac{\alpha_2}{2} - \frac{\beta_2}{2 \sqrt{\varepsilon}}\right)
		    	 	  - \frac{\beta_1 \beta_2 - \varepsilon \alpha_1 \alpha_1}{4 \sqrt{\varepsilon}}
		    	 	  - \frac{\varepsilon (\alpha_1 \gamma_2 + \alpha_2 \gamma_1)}{2 (1-\varepsilon^2)}
		    	 	  + \frac{\varepsilon (\beta_1 \gamma_2 + \beta_2 \gamma_1)}{2 \sqrt{\varepsilon }(1-\varepsilon^2)}=0\,.
			\end{array}
			\right.			
		\een
		Again we find two families of solutions of this system. The first one is given by
		 $$\,
		   \beta_1=\alpha_1\sqrt{\varepsilon},\,
			 \beta_2=\alpha_2\sqrt{\varepsilon},\,
		\,$$
		and $\gamma_j,\,\varepsilon$ are arbitrary. 
    The second family is defined by
				$$\,
				  \beta_1=(\alpha_1-2)\sqrt{\varepsilon},\,
					\beta_2=(\alpha_2+2)\sqrt{\varepsilon},\,\,\,\,
					\frac{\alpha_2-\alpha_1}{2} - \frac{\sqrt{\varepsilon}(\gamma_2-\gamma_1)}{1-\varepsilon^2}=-1\,.
				\,$$
				These are the last two cases  -- the case $(VII.)$ and the case $(VIII.)$. In both cases the 
				equation  \eqref{pe-ope} - \eqref{pe} admits 4 regular singularities taken at $x_1, x_2, x_4$ and $x_5$.
				
				This ends the proof.
		\qed

       In this paper we study the equation \eqref{pe-ope} - \eqref{pe} under the condition $(I.)$ when
       the coefficients of the equation \eqref{npe} - \eqref{op-ini} do not depend on the parameter of perturbation 
       $\varepsilon$. 
      It is well known that every 4 distinct points over $\CC\PP^1$ can be fixed by 
      a M\"{o}bius transformation at $0, 1, \infty$ and $a\neq 0, 1, \infty$.
      Recall that the Heun equation is a second order Fuchsian equation with four singular points over $\CC\PP^1$
      taken at $0, 1, \infty$ and $a \neq 0, 1, \infty$, \cite{S-L}. Since a M\"{o}bius transformation can take the perturbed
      equation into the Heun equation through this paper we call the perturbed equation Heun type equation.      
      We prefer to work with the perturbed equation rather than work with Heun equation because of the symmetries.
      The next proposition gives all the more reason for working with the perturbed equation instead of the Heun equation.
      
       \bpr{heun}      
           The M\"{o}bius transformation
            $$\,
              x= \sqrt{\varepsilon}\,\frac{t+ \frac{1-\varepsilon}{1+\varepsilon}}{t - \frac{1-\varepsilon}{1+\varepsilon}}
            \,$$
            takes the perturbed equation into the Heun equation
            \ben
             \ddot{y}
              &+&
              \left[\left(\frac{\beta_1+\beta_2}{2 \sqrt{\varepsilon}} +1\right) \frac{1}{t}
               + \frac{\gamma_1+\gamma_2}{2 \sqrt{\varepsilon}} \frac{1}{t-1}
                - \frac{\gamma_1+\gamma_2}{2 \sqrt{\varepsilon}}
                \frac{1}{t- \left(\frac{1-\varepsilon}{1+\varepsilon}\right)^2}\right] \dot{y}+\\[0.25ex]
                &+&
             \Big[\frac{\beta_1 \beta_2}{4 \varepsilon} \frac{1}{t^2} +
                \frac{\gamma_1}{2 \sqrt{\varepsilon}} \left(\frac{\gamma_2}{2 \sqrt{\varepsilon}}-1\right)\frac{1}{(t-1)^2}\Big]   
             + \frac{\gamma_1}{2 \sqrt{\varepsilon}} \left(\frac{\gamma_2}{2 \sqrt{\varepsilon}}+1\right)
             \frac{1}{\left(t-\left(\frac{1-\varepsilon}{1+\varepsilon}\right)^2\right)^2} +\\[0.25ex]
               &+&
               \left(\frac{\beta_1 \gamma_2 + \beta_2 \gamma_1}{4 \varepsilon}
               + \frac{\gamma_1}{2 \sqrt{\varepsilon}}\right) \frac{4 \varepsilon}{(1-\varepsilon)^2} \frac{1}{t}
               +  \left(\frac{\beta_1 \gamma_2 + \beta_2 \gamma_1}{4 \varepsilon}+
                \frac{\gamma_1}{2 \sqrt{\varepsilon}} - \frac{(1+\varepsilon)^2}{8 \varepsilon^2} \gamma_1 \gamma_2\right) 
                 \frac{1}{t-1}+\\[0.25ex]
                &+&
                \left(-\frac{\beta_1 \gamma_2 + \beta_2 \gamma_1}{4 \varepsilon}\left(\frac{1+\varepsilon}{1-\varepsilon}\right)^2
                -\frac{\gamma_1}{2 \sqrt{\varepsilon}}\left(\frac{1+\varepsilon}{1-\varepsilon}\right)^2
                +\frac{(1+\varepsilon)^2}{8 \varepsilon^2} \gamma_1 \gamma_2\right)
                \frac{1}{t - \left(\frac{1-\varepsilon}{1+\varepsilon}\right)^2} \Big] y=0\,,   
            \een
            which does not have limit when $\varepsilon \rightarrow 0$.
       \epr
       
     \proof
     The proof is straightforward.\qed

   We leave to the reader as an easy exercise to check that there is no Heun equation, obtained from the perturbed equation
   by a M\"{o}bius transformation, that has a limit when $\varepsilon \rightarrow 0$.

%%%%%%%%%%%%%%%%%%%%%%%%%%%%%%%%%%%%%%%%%%%%%%%%%%
%painleve
%%%%%%%%%%%%%%%%%%%%%%%%%%%%%%%%%%%%%%%%%
           \subsection{The initial equation and the third Painlev\'{e} equation}
           
           In \cite{S-S} Salatich and Slavyanov applying an antiquantization procedure to the DCHE \eqref{dche}
           obtain a particular case of the third Painlev\'{e} equation $P_{III}$
            \be\label{p3}
              \ddot{q} - \frac{\dot{q}^2}{q} + \frac{\dot{q}}{t} - \gamma\,q^3 
              - \frac{\alpha\,q^2 + \beta}{t} - \frac{\delta}{q}=0\,,
            \ee
           where $\alpha, \beta, \gamma, \delta\in\CC$ are parameters and $\cdot=\frac{d}{d t}$. In this section
           following their approach we connect by an antiquantization particular families of the initial equation
           with particular families of the $P_{III}$ equation \eqref{p3}.
           
           Under the non-resonant condition \eqref{nonr} there are four particular families of the initial equation
           that satisfy the conditions \eqref{p-dche}. They are defined as follows:
            \begin{enumerate}
            	\item\,
            	$\gamma_2=\beta_1=0,\, \gamma_1=1,\,c=2, \,a=2,\, h=-1$ and $\beta_2=t \neq 0$ is an arbitrary. Then 
            	the corresponding DCHE is
            	 \be\label{1}
            	   x^2\,y''(x) + (-x^2 + 2 x - t)\,y'(x) + (-2 x + t)\,y(x)=0\,.
            	 \ee
            	 
            	 \item\,
            	 $\gamma_2=\beta_2=0,\,\gamma_1=1, \,c=2,\,a=2, \,h=0$ and $\beta_1=t\neq 0$ is an arbitrary. Then the 
            	 corresponding DCHE is
            	 	 \be\label{2}
            	 	 x^2\,y''(x) + (-x^2 + 2 x - t)\,y'(x) -2 x \,y(x)=0\,.
            	 	 \ee
            	 	 
            	 	\item\,
            	$\gamma_1=\beta_1=0,\,\gamma_2=1,\,c=2,\,a=h=0$ and $\beta_2=t\neq 0$ is an arbitrary. Then the
            	corresponding DCHE is
            	 	 \be\label{3}
            	 	 x^2\,y''(x) + (-x^2 + 2 x - t)\,y'(x) =0\,.
            	 	 \ee 	 
            	 	 
            	\item\,
            	$\gamma_1=\beta_2=0,\,\gamma_2=1,\,c=2,\,a=0,\,h=-1$ and $\beta_1=t \neq 0$ is an arbitrary. Then the
            	corresponding DCHE is
            	 	 \be\label{4}
            	 	 x^2\,y''(x) + (-x^2 + 2 x - t)\,y'(x) + t\,y(x)=0\,.
            	 	 \ee 	 
            \end{enumerate}
             Note that for the initial equation $c, a, h$ from \eqref{dche} are fixed and $t$ is a varying parameter.
           Denote
            $$\,
              \sigma(x)=x^2,\quad \tau(x)=-x^2 + 2 x - t,\quad \omega(x)=-a x\,.
            \,$$                     
            Following \cite{S-S} we asociate with the DCHE the Hamiltonian
             $$\,
               H(q, p)=\frac{1}{t}\left[\sigma(q)\,p^2 + \tau(q)\,p + \omega(q)\right]\,.
             \,$$
             The Hamiltonian system is written as
              \ben
                \dot{q} 
                   &=&
                  \frac{\partial H}{\partial p}=\frac{2 \sigma(q)\,p + \tau(q)}{t}=
                  \frac{2 q^2\,p - q^2 + 2 q -t}{t},\\[0.2ex]
                \dot{p}
                 &=&
                - \frac{\partial H}{\partial q}=-\frac{\sigma'(q)\,p^2 + \tau'(q)\,p + \omega'(q)}{t}=
                -\frac{2 q\,p^2 -2 q\,p + 2 p -a}{t}\,.    
              \een
            Eliminating from this system $p$ we obtain the following non-linear second order equation for $q$
            $$\,
             \ddot{q} - \frac{\dot{q}^2}{q} + \frac{\dot{q}}{t} - \frac{q^3}{t^2} -
             \frac{2 (a-1)\,q^2}{t^2} - \frac{1}{t} + \frac{1}{q}=0\,,
            \,$$
            which is called the third Painlev\'{e} $P_{III'}$ equation \cite{O}.
            The transformations
             $$\,
              t=z^2,\qquad q=z\,u
             \,$$ 
            take the last equation into the particular $P_{III}$
            $$\,
             u'' - \frac{u'^2}{u} + \frac{u'}{z} - 4 u^3 - \frac{4 (2 (a-1)\,u^2 + 1)}{z} + \frac{4}{u}=0
            \,$$
           with $\gamma=4,\,\alpha=8 (a-1),\,\beta=4$ and $\delta=-4$. As a result, the first two families of the
           initial equation are related to $P_{III}$ equation with
            $$\,
             \gamma=-4,\quad \alpha=8,\quad \beta=4,\quad \delta=-4\,.
            \,$$ 
           The last two families are related to the $P_{III}$ equation with 
                 $$\,
                 \gamma=-4,\quad \alpha=-8,\quad \beta=4,\quad \delta=-4\,.
                 \,$$ 
            These two families of the $P_{III}$ equation fall into the generic case of the  third Painlev\'{e} equation. It is defined 
            by the condition $\gamma\,\delta \neq 0$ and it is well studied in \cite{Mu, O}. Note that when $a=0$ 
            the $P_{III'}$ equation has a Riccaty type of solution. Indeed, when $a=0$ the function $p\equiv 0$ satisfies the
            second equation of the above Hamiltonian system. Then the function $q$ satisfies the Riccaty equation 
             $$\,
              t\,\dot{q}=-q^2 + 2 q -t\,.
             \,$$
            Similarly, when $a=2$ the $P_{III'}$ equation again has a Riccaty type of solution. This time the function
            $p=1$ is a particular solution of the Hamiltonian system. Then the function $q$ satisfies the Riccati equation
                  $$\,
                  t\,\dot{q}=q^2 + 2 q -t\,.
                  \,$$

%%%%%%%%%%%%%%%%%%%%%%%%%%%%%%%%%%%%%%%%%%%%
%itterated integrals
%%%%%%%%%%%%%%%%%%%%%%%%%%%%%%%%%%%%%%%%%%%%%
   \subsection{Global solutions}
   In \cite{St1} we have proved that reducible differential equations by means of the first-order
   differential equations admit a global fundamental matrix, whose off-diagonal elements are
   represented by iterated integrals. More precisely,  		

    \bth{global}
     Both initial equation  and perturbed  equations admit a global fundamental matrix $\Phi(x, \cdot)$ of the form
       \be\label{gfm}
          \Phi(x, \cdot)=\left(\begin{array}{cc}
          	 \Phi_1(x, \cdot)   &\Phi_{12}(x, \cdot)\\
          	  0                 &\Phi_2(x, \cdot)
          	                 \end{array}
          	             \right)\,.       
       \ee
       The diagonal elements $\Phi_j(x, \cdot), j=1, 2$ are the solutions od the equations
       $L_{j, \cdot}\,u=0$ with $L_{j, \varepsilon}$ and $L_{j, 0}=L_j$ given by. The 
       off-diagonal element is defined as
        \be\label{off}
          \Phi_{12}(x, \cdot)=\Phi_1(x, \cdot)
           \int_{\Gamma(x, \cdot )}
             \frac{\Phi_2(z, \cdot)}{\Phi_1(z, \cdot)}\, d z\,.
        \ee
        The paths of integration $\Gamma(x, \varepsilon)$ and $\Gamma(x, 0)$ are taken from 
        the same base point $x$ in such a way that $\Gamma(x, \varepsilon) \rightarrow \Gamma(x, 0)$
        as $\varepsilon \rightarrow 0$, and the matrices $\Phi(x, \cdot)$ are fundamental matrix
        solution of the corresponding equations.
    \ethe		
		
    As an immediate corollary we define a fundamental set of solutions of both equations.
    
    \bpr{fss}
      Bot initial and perturbed equations possess a fundamental set of solutions of the form
       $$\,
          \Phi_1(x, \cdot),\quad
                \Phi_1(x, \cdot)
                \int_{\Gamma(x, \cdot )}
                \frac{\Phi_2(z, \cdot)}{\Phi_1(z, \cdot)}\, d z\,.       
       \,$$
    \epr

%%%%%%%%%%%%%%%%%%%%%%%
%irregular
%%%%%%%%%%%%%%%%%%%%%%%%%%%%%

   \subsection{Irregular singularities}

  The initial equation is a second order linear differential equation having two irregular singularities.
  They are taken at the origin and at the infinity point and both singular points are of Poincar\'{e} rank 1.
  In this paper, as in our previous works \cite{St1}, we utilize summability theory to build actual 
  fundamental matrices at the origin and at the infinity point. Then with respect to these fundamental matrices
  we compute the corresponding Stokes matrices. In this paragraph we review some definitions, facts and notation 
  from the applications of summability theory
  to the ordinary differential equations which is required to compute the Stokes matrices of the initial equation.
  We consider in parallel the summation at the origin and at the infinity point following works of Balser \cite{WB} and
  Ramis \cite{R2} for the summation at the origin, and  Sauzin \cite{Sa} for the summation at $x=\infty$.

    We denote by $\CC[[x]]$ the field of formal power series at the origin
		$$\,
		  \CC[[x]]=\left\{ \sum_{n=0}^{\infty} f_n\,x^n\,|\,f_n\in\CC, \,n\in\NN\right\}\,,
		\,$$
		and by $\CC[[x^{-1}]]$ the field of formal power series at $\infty$
		$$\,
		\CC[[x^{-1}]]=\left\{ \sum_{n \geq 0} \varphi_n\,x^{-n}\,|\,\varphi_n\in\CC, \,n\in\NN\right\}\,.
		\,$$
    Equipping $\CC[[x]]$ (resp. $\CC[[x^{-1}]]$) with the natural derivation $\partial=\frac{d}{d x}$ such that
     $$\,
       \partial(\psi\,\varphi)=\partial (\psi)\,\varphi + \psi\,\partial(\varphi),\qquad
       \psi, \varphi\in\CC[[x]] \quad (\textrm{resp.}\,\,\, \psi, \varphi\in\CC[[x^{-1}]])
     \,$$
  we make $\CC[[x]]$ (resp. $\CC[[x^{-1}]]$) a differential algebra.
   All singular directions and sectors are defined on the Riemann surface of the natural logarithms.		

    \bde{sec}
    1.\,An open sector $S$ with vertex 0 is a set of the form
        $$\,
          S=S(\theta, \alpha, \rho)=\left\{\,x=r\,e^{i \delta}\,|\,
          0 < r < \rho,\,\theta-\alpha/2 < \delta < \theta + \alpha/2\right\}\,.
        \,$$
        \,An open sector $S_1$ with vertex $\infty$ is a set of the form
        $$\,
        S_1=S_1(\theta, \alpha, R)=\left\{\,x=r\,e^{i \delta}\,|\,
         r > R,\,\theta-\alpha/2 < \delta < \theta + \alpha/2\right\}\,.
        \,$$
        Here $\theta$ is an arbitrary real number (the bisector of the sector), $\alpha$ is a positive
        real (the opening of the sector) and $\rho$ (resp. $R$) is either a positive real number or $+\infty$
        (the radius of the sector).\\
        2.\, A closed sector $\bar{S}$ vith vertex 0 is a set of the form
             $$\,
             \bar{S}=\bar{S}(\theta, \alpha, \rho)=\left\{\,x=r\,e^{i \delta}\,|\,
             0 < r \leq \rho,\,\theta-\alpha/2 \leq \delta \leq \theta + \alpha/2\right\}\,.
             \,$$
             \,A closed  sector $\bar{S}_1$ with vertex $\infty$ is a set of the form
             $$\,
             \bar{S}_1=\bar{S}_1(\theta, \alpha, R)=\left\{\,x=r\,e^{i \delta}\,|\,
             r \geq R,\,\theta-\alpha/2 \leq \delta \leq \theta + \alpha/2\right\}\,.
             \,$$
             Here $\theta$ and $\alpha$ are as before, but $\rho$ (resp. $R$) is  a positive real number 
             (never equal to $+\infty$).
    \ede

    \bde{asy}\,1.\.
		  Let the function $f(x)$ be holomorphic on a sector $S(\theta, \alpha, \rho)$.
			The (formal) power series $\hat{f}(x)=\sum_{n=0}^{\infty} f_n\,x^n$ is said
			to represent $f(x)$ asymptotically, as $x \rightarrow 0$ in $S$, if for every
			closed sector $W \subset S$ and all $N \geq 0$ there exist a positive constant
			$C(N, W)$ such that
			 $$\,
			   \left| f(x) - \sum_{n=0}^{N-1} f_n\,x^n \right| \leq C(N, W)\,|x|^N\,,\quad
				x\in W\,.
			\,$$
			2.\,
				  Let the function $\varphi(x)$ be holomorphic on a sector $S_1(\theta, \alpha, R)$.
				  The (formal) power series $\hat{\varphi}(x)=\sum_{n\geq 0} \varphi_n\,x^{-n}$ is said
				  to represent $\varphi(x)$ asymptotically, as $x \rightarrow \infty$ in $S_1$, if for every
				  closed sector $W \subset S_1$ and all $N \geq 0$ there exist a positive constant
				  $K(N, W)$ such that
				  $$\,
				  \left| \varphi(x) - \sum_{n=0}^{N-1} \varphi_n\,x^{-n} \right| \leq K(N, W)\,|x|^{-N}\,,\quad
				  x\in W\,.
				  \,$$
		\ede
		In this case one usually write
		  \ben
		          & &
			   f(x) \sim \hat{f}(x)\,,\quad x\in S\,,\quad x \rightarrow 0\,,\\
			      & &
			      \varphi(x) \sim \hat{\varphi}(x)\,,\quad x\in S_1\,,\quad x \rightarrow\infty\,.
			\een
		A function $f(x)$ (resp. $\varphi(x))$ can have at most one asymptotic series representation as
		$x \rightarrow 0$  (resp. $x \rightarrow \infty$) in a given sector $S$ (resp. $S_1$). Moreover, 
		the set $\mathcal{A}(S)$ (resp. $\mathcal{A}(S_1)$) 
		of functions, which are holomorphic on the sector $S$ (resp. $S_1$) and admit asymptotic 
		representation on this sector forms a differential algebra. The maps
		 \ben
		   \mathcal{A}(S) &\longrightarrow& \CC[[x]]\qquad \textrm{and}\qquad
		     \mathcal{A}(S_1) \longrightarrow \CC[[x^{-1}]] \\
			     f(x)       &\longrightarrow& \hat{f}(x)\qquad\qquad\qquad\quad
			        \varphi(x)       \longrightarrow \hat{\varphi}(x)
		\een
		are  homomorphism of differential algebras. The Borel - Ritt Theorem implies that these maps are surjective maps \cite{WB, R2}.
		Unfortunately they are  not injective maps.
  
    Among formal power series $\CC[[x]]$ (resp. $\CC[[x^{-1}]]$)
  we distinguish the formal power series of Gevrey order 1.
  
   \bde{gev}
      A formal power series $\hat{f}(x)=\sum_{n=0}^{\infty} f_n\,x^n$ 
      (resp. $\hat{\varphi}(x)=\sum_{n \geq 0} \varphi_n \,x^{-1}$) is said to be of Gevrey order 1
      if there exist two positive constants $C, A > 0$ such that
        \ben
          |f_n| &<&  C\,A^n \,n! \quad \textrm{for every}\quad n\in\NN\\
          (\textrm{resp.}\quad
              |\varphi_n| &<&  C\,A^n \,n! \quad \textrm{for every}\quad n\in\NN)\,.
        \een 
   \ede
    We denote by $\CC[[x]]_1$ and $\CC[[x^{-1}]]_1$ the sets of all power series at $x=0$ 
    and at $x=\infty$ of Gevrey order 1. These sets are sub-algebras of $\CC[[x]]$ and $\CC[[x^{-1}]]$,
    respectively, as commutative differential algebras over $\CC$, \cite{L-R, S}.

   \bde{bo1}
     The formal Borel transform $\hat{\B}_1$ of order 1 of a formal power seies
		$\hat{f}(x)=\sum_{n=0}^{\infty} f_n\,x^n$ is called the formal series
		  $$\,
			  \hat{\B}_1\,\hat{f} (\xi)= \sum_{n=0}^{\infty}
				\frac{f_n}{n!}\,\xi^n\,.
			\,$$
    \ede
	 Likewise,  we introduce the formal Borel transform of a formal power series
	near $x=\infty$. 
   \bde{bo21}
	  The formal Borel transform $\hat{\B}_1$ of order 1 of a formal
		power series $\hat{\varphi}(x)=\sum_{n=0}^{\infty} \varphi_n\,x^{-n-1}$
		near $x=\infty$ is called the formal series
		 $$\,
		   \hat{\B}_1\,\hat{\varphi}(p)=\sum_{n=0}^{\infty} \frac{\varphi_n}{n!} \,p^n\,.
		\,$$
	\ede
  If $\hat{f}(x)\in\CC[[x]]_1$ (resp. $\hat{\varphi}(x)\in\CC[[x^{-1}]]_1$) then its formal Borel transform $\hat{\B}_1$
  converges in a neighborhood of the origin $\xi=0$ (resp. $p=0$) with a sum $f(\xi)$ (resp. $\varphi(p)$).

   \bde{lap-0}
    Let $f(\xi)$ be analytic and of exponential size at most 1 at $\infty$, i.e.
    $|f(\xi)| \leq A\,\exp (B |\xi|),\,\xi\in \theta$ along a direction $\theta$ from $0$ to $+\infty e^{i \theta}$.
    Then, the integral
     $$\,
       \left(\mathcal{L}_{\theta} f\right)(x)=
       \int_0^{+\infty e^{i \theta}}
       f(\xi)\,\exp\left(-\frac{\xi}{x}\right)\, d \left(\frac{\xi}{x}\right)
     \,$$ 
       is said to be Laplace complex transform $\mathcal{L}_{\theta}$ of order 1 in the direction $\theta$
       of $f$.
        \ede

     \bde{lap-in}
     Let $\varphi(p)$ be analytic and of exponential size at most 1 at $\infty$, i.e.
     $|\varphi(p)| \leq A\,\exp (B |p|),\,p\in \theta$ along a direction $\theta$ from $0$ to $+\infty e^{i \theta}$.
     Then, the integral
     $$\,
     \left(\mathcal{L}_{\theta} \varphi\right)(x)=
     \int_0^{+\infty e^{i \theta}}
     \varphi(p)\,\exp(p\,x)\, d p
     \,$$ 
     is said to be Laplace complex transform $\mathcal{L}_{\theta}$ of order 1 in the direction $\theta$
     of $\varphi$.
     \ede
              Let the function $\varphi(p)$ be analytic and satisfies the estimate
   \be\label{est}
         |\varphi(p)| \leq A_0\,e^{c_0\,|p|},\quad c_o\in\RR
   \ee
      along a direction $\theta$ from 0 to $+\infty\,e^{i \theta}$. Then the Laplace complex transform $\L_{\theta}$
      in the direction $\theta$ satisfies some useful properties.
       \ble{prop}\cite{Sa}
         Let $\varphi$ as above, $\tilde{\varphi}:=(\L_{\theta} \varphi)$ and $c\in\CC$. Then each of the functions
         $-p\,\varphi(p),\,e^{-c p}\,\varphi(p)$ or $1\ast \varphi(p)$satisfies estimates of the form \eqref{est}
         and
          \begin{itemize}
          	\item\,
          	$(\L_{\theta} (-p\,\varphi))=\frac{d \tilde{\varphi}}{d x}$\,,
          	
          	\item\,
          	$(\L_{\theta} (e^{-c p}\,\varphi))=\tilde{\varphi} (x+c)$\,,

             \item\,
             $(\L_{\theta} (1\ast \varphi)) = x^{-1}\,\tilde{\varphi}(x)$\,,
             
             \item\,
             if moreover $\varphi$ is analyticaly derivable on $\theta$ with $\frac{d \varphi}{d p}$ satisfying
             estimates of the form \eqref{est}, then
             $$\,
               \left(\L_{\theta} (\frac{d \varphi}{d p})\right)=x\,\tilde{\varphi}(x) - \varphi(0)\,.
             \,$$ 
          \end{itemize}
       \ele
       
   Now we can give the definitions of the 1-summable series in a direction.
      \bde{sum-0}
       The formal power series $\hat{f}=\sum_{n=0}^{\infty} f_n\,x^n$ is 1-summable (or Borel summable)
       in the direction $\theta$ if there exist an open sector $V$ bisected by $\theta$ whose opening is
       $> \pi$ and a holomorphic function $f(x)$ in $V$ such that for every non-negative integer $N$
        $$\,
           \left| f(x) - \sum_{n=0}^{N-1} f_n\,x^n\right| \leq C_{V_1}\,A^N_{V_1}\,N!\,|x|^N
        \,$$
        on every closed subsector $\bar{V}_1$ of $V$ with constants $C_{V_1},\,A_{V_1} > 0$ depending only
        on $V_1$. The function $f(x)$ is called the 1-sum (or Borel sum) of $\hat{f}(x)$ in the direction
        $\theta$.
      \ede
      If a series $\hat{f}(x)$ is 1-summable in all but a finite number of directions, we will say that it is 1-summable.
      
      One useful criterion for a Gevrey series of order 1 to be 1-summable is given in terms of Borel and Laplace 
      transforms:
      
       \bpr{sum-0}
        Let $\hat{f}\in \CC[[x]]_1$  $($resp. $\hat{\varphi}\in\CC[[x^{-1}]]_1)$ and let $\theta$ be a direction. The following are equivalent:
          \begin{enumerate}
          	\item\,
          	$\hat{f}$ $($resp. $\hat{\varphi})$ is 1-summable in the direction $\theta$.
          	
          	\item\,
          	The convergent power series $\hat{\B}_1\,\hat{f}(\xi)$ $($resp. $\hat{\B}_1\,\hat{\varphi}(p))$ has an analytic continuation $h$ $($resp. $g)$ in a
          	full sector $\{\, j\in\CC\,|\,0 < |j| < \infty,\,|\arg (j) - \theta| < \epsilon\, \}$ for
          	$j=\xi$ $($resp. $j=p)$.
          	In addition, this analytic continuation has exponential growth of order 1 at $\infty$ on this
          	sector, i.e. $|h(\xi)| \leq A\,\exp (B\,|\xi|)$ $($resp. $|g(p)| \leq A_0\,\exp(c_0\,|p|)$. 
          	In this case $f=(\L_{\theta} h)(x)$ $($resp. $\varphi=(\L_{\theta} g)(x))$ is its
          	1-sum in the direction $\theta$.
          \end{enumerate}
       \epr

    \bre{conv}
      If $\hat{f}(x)$ $($resp. $\hat{\varphi}(x))$ is convergent, then $\hat{f}(x)$ $($resp. $\hat{\varphi}(x))$
      is $k$-summable in the direction $\theta$ (any $k > 0$ and any $\theta$) and the classical sum and
      $k$-sum coincide in their respective domain of definition, \cite{R2, Sa}.
    \ere

  The next result, which can be found in \cite{St2} (see Lemma  4.1(I)
	in \cite{St2}) provides 1 -summability of a class of formal power series
	near $x=0$.
	 
	\ble{example}
	  Assume that $\beta_j\neq \beta_i$ and that $\alpha_j - \alpha_i \notin \ZZ_{\leq -2}$. 
	  The formal power series
		 $$\,
		  \hat{f}(x)= \sum_{n=0}^{\infty} (-1)^n\,
		  \frac{(2 + \alpha_j-\alpha_i)^{(n)}}{(\beta_j-\beta_i)^n}\,x^n
		\,$$
		is 1 - summable in any direction $\theta\neq \arg(\beta_i -\beta_j)$ from $0$ to $+\infty e^{i \theta}$.
		The function 
		 $$\,
		  f_{\theta}(x)=\int_0^{+\infty e^{i \theta}}
		  \left(1+ \frac{\xi}{\beta_j- \beta_i}\right)^{\alpha_i-\alpha_j-2}\,
		  e^{-\frac{\xi}{x}}\,d \left(\frac{\xi}{x}\right)
		 \,$$ 
		 defines its 1-sum in such a direction.
	\ele
        Here by the symbol $(a)^{(n)}$ we denote the rising factorial
        $$\,
          (a)^{(n)}=a\,(a+1)\,(a+2) \ldots (a+n-1),\qquad (a)^{(0)}=1\,.
        \,$$

    As a differential equation with an irregular point at the origin and an irregular point at $x=\infty$
    the initial equation admits unque formal fundamental matrices $\hat{\Phi}_0(x)$ and $\hat{\Phi}_{\infty}(x)$
    in the form of the theorem of Hukuhata-Turrittin \cite{W, S}
      $$\,
       \hat{\Phi}_0(x, 0)=\hat{L}(x)\,x^{\Lambda}\,\exp\left(-\frac{B}{x}\right),\quad
           \hat{\Phi}_{\infty}(x, 0)=\hat{Q}(x)\,\left(\frac{1}{x}\right)^{-\Lambda}\,\exp (G\,x)\,.
      \,$$
     Here the matrices $\Lambda, B$ and $G$ are diagonal and
      $$\,
       \Lambda=\diag (0, -2),\quad B=\diag (\beta_1, \beta_2),\quad G=\diag(\gamma_1, \gamma_2)\,.
      \,$$
      Once having formal fundamental matrices we can introduce the so called formal monodromy matrices
      $\hat{M}_0$ and $\hat{M}_{\infty}$ at the origin and at the infinity point.
      
      \bde{f-M}
       The formal monodromy matrix $\hat{M}_0$ related to the formal fundamental matrix $\hat{\Phi}_0(x)$
       is defined as
         $$\,
          \hat{\Phi}_0(x.e^{2 \pi\,i}, 0)=\hat{\Phi}_0(x, 0)\,\hat{M}_0\,.
         \,$$
         In the same manner, the formal monodromy matrix $\hat{M}_{\infty}$ related to the fromal fundamental
         matrix $\hat{\Phi}_{\infty}(x, 0)$ is defined as
           $$\,
           \hat{\Phi}_{\infty}(x^{-1}.e^{2 \pi\,i}, 0)=\hat{\Phi}_{\infty}(x, 0)\,\hat{M}_{\infty}\,.
           \,$$
           In particular,
             $$\,
              \hat{M}_0=e^{2 \pi\,i\,\Lambda}=I_2=e^{-2 \pi\,i\,\Lambda}=\hat{M}_{\infty}\,,
             \,$$
             where $I_2$ is the identity matrix of order 2.
      \ede
      
  Since the formal monodromy matrices are equal to the identity matrix, in this paper we will present the formal fundamental matrices $\hat{\Phi}_0(x, 0)$ and
  $\hat{\Phi}_{\infty}(x. 0)$ in a slight different form
   \be\label{formal}
     \hat{\Phi}_0(x, 0)  &=& \exp (G\,x)\,\hat{H}(x)\,x^{\Lambda}\,\exp\left(-\frac{B}{x}\right),\\[0.2ex]
     \hat{\Phi}_{\infty}(x, 0) &=&
     \exp\left(-\frac{B}{x}\right)\,\left(\frac{1}{x}\right)^{-\Lambda}\,\hat{P}(x)\,
     \exp( G\,x)\,.\nonumber
   \ee 
  These special forms of the formal fundamental solutions, as well as, of the actual fundamental solutions
  allow us to show in an explicit way haw these solutions are changed under the perturbation. 
   Here the matrix-function $x^{\Lambda}, \exp(G x)$ and $\exp(-B/x)$ must be regarded as formal function.
   The entries of the matrices $\hat{H}(x)$ and $\hat{P}(x)$ are formal power series in $x$ and $x^{-1}$,
   respectively.

   Utilizing the summability theory we relate to the formal fundamental matrices $\hat{\Phi}_0(x, 0)$ and
   $\hat{\Phi}_{\infty}(x, 0)$ actual fundamental matrices. More precisely,
   
    \bth{ac}(Hukuhara-Turrittin- Martinet - Ramis) The entries of the matrix $\hat{H}(x)$
    (resp. $\hat{P}(x)$) in \eqref{formal} are 1-summable in every non-singular direction $\theta$.
    If we denote by $H_{\theta}(x)$ (resp. $P_{\theta}(x)$) the 1-sum of $\hat{H}(x)$ (resp. $\hat{P}(x)$) along $\theta$
    obtained from $\hat{H}(x)$ (resp. $\hat{P}(x)$) by a Borel - Laplace transform, then
    $\Phi^{\theta}_0(x, 0)=e^{G x}\,H_{\theta}(x)\,x^{\Lambda}\,e^{-B/x}$
    (resp. $\Phi^{\theta}_{\infty}(x, 0)=e^{-B/x}\,x^{\Lambda}\,P_{\theta}(x)\,e^{G x}$) is an actual fundamental
    matrix at the origin (resp. at $x=\infty$) of the initial equation.
 \ethe

     Let $\theta$ be a singular direction of the initial equation at the origin. Let $\theta^{+}=\theta + \epsilon$ and
     $\theta^{-}=\theta-\epsilon$, where $\epsilon > 0$ is a small number, be two non-singular neighboring directions
     of the singular direction $\theta$. Denote by $\Phi^{\theta+}_0(x)$ and $\Phi^{\theta-}_0(x)$ the actual
     fundamental matrices at the origin of the initial equation corresponding to the direction $\theta^{+}$ and
     $\theta^{-}$ in the sense of \thref{ac}. Then
     
        \bde{st-0}
          With respect to the given formal fundamental matrix $\hat{\Phi}_0(x, 0)$ at the origin the Stokes matrix
          $St^{\theta}_0\in GL_2(\CC)$ corresponding to the singular direction $\theta$ is defined as
            $$\,
              St^{\theta}_0=(\Phi^{\theta+}_0(x, 0))^{-1}\,\Phi^{\theta-}_0(x, 0)\,.
            \,$$
        \ede
     
     Similarly, 
        let $\theta$ be a singular direction of the initial equation at the infinity point. Let $\theta^{+}=\theta + \epsilon$ and
        $\theta^{-}=\theta-\epsilon$, where $\epsilon > 0$ is a small number, be two non-singular neighboring directions
        of the singular direction $\theta$. Denote by $\Phi^{\theta+}_{\infty}(x)$ and $\Phi^{\theta-}_{\infty}(x)$ the actual
        fundamental matrices at the origin of the initial equation corresponding to the direction $\theta^{+}$ and
        $\theta^{-}$ in the sense of \thref{ac}. Then
        
        \bde{st-in}
        With respect to the given formal fundamental matrix $\hat{\Phi}_{\infty}(x, 0)$ at the infinity point the Stokes matrix
        $St^{\theta}_{\infty}\in GL_2(\CC)$ corresponding to the singular direction $\theta$ is defined as
        $$\,
        St^{\theta}_{\infty}=(\Phi^{\theta+}_{\infty}(x, 0))^{-1}\,\Phi^{\theta-}_{\infty}(x, 0)\,.
        \,$$
        \ede
     
    %%%%%%%%%%%%%%%%%%%%%%%%%%%%%%%%%%%
    %Fuchsian
    %%%%%%%%%%%%%%%%%%%%%%%%%%%%
    \subsection{ Fuchsian singularities }
     
       We write $\NN_0=\NN \cup \{0\}$.

    The considered in this paper Heun type equation is a Fuchsian equation of order 2 with 4 finite
    singular points  taken at $x_R=\sqrt{\varepsilon}, x_L=-\sqrt{\varepsilon}, x_{RR}=1/\sqrt{\varepsilon}$
    and $x_{LL}=-1/\sqrt{\varepsilon}$. In this paragraph we briefly introduce the needed facts and definitions 
    from the theory of the Fuchsian equations.
        
     We firstly recall  the necessary local theory of Fuchsian singularities following the book of Golubev \cite{Go}. 
     Recall that in the case of
     the scalar differential equations the regular singularity and the Fuchsian singularity are the same notion.   
   With every Fuchsian (regular) singularity of a given $n$-order scalar linear differential equation
   we associate a $n$-order algebraic equation, the so called characteristic (or inditial) equation.
   More precisely, consider an $n$-order linear differential equation
    $$\,
      y^{(n)}(x) + b_{n-1}(x)\,y^{(n-)}(x)+ \cdots + b_0(x)\,y(x)=0,\quad
      b_j(x)\in\CC(x)\,.
    \,$$
   Let $x=x_0\in\CC$ be a regular singularity for this equation. Recall that due to the theorem of
   Fuchs this means that all the functions
    $$\,
       b_{n-k}(x)\,(x-x_0)^k
    \,$$
  are holomorphic functions at $x=x_0$. Then
  
     \bde{ce}
     1.\,The $n$-order algebraic equation
        $$\,
         \rho (\rho-1) (\rho-2) \ldots (\rho-(n-1)) +
         c_{n-1}\,\rho (\rho-1)\ldots (\rho-(n-2)) + \cdots + c_1\,\rho + c_0=0\,,
        \,$$
        where 
        $$\,
           c_k=\lim_{x \rightarrow x_0} b_{n-k}(x)\,(x-x_0)^k,\quad  0 \leq k \leq n-1
        \,$$
        is called the characteristic (or the inditial) equation at the regular singularity
        $x_0\in\CC$. Its roots $\rho_k, 1 \leq k \leq n$ are called the characteristic exponents at
        the singularity $x_0$.\\
      2.\,The characteristic equation at the point $t=0$ of the equation obtained after the
      transformation $x=1/t$, is called the characteristic equation at $x=\infty$. Its roots
      are called the characteristic exponents at the regular point $x=\infty$.   
     \ede

        Denote by $\rho^j_i$ and $\rho^{jj}_i,\,i=1, 2,\,j=R, L$ the characteristic exponents at the
        singular points $x_j$ and $x_{jj},\,j=R, L$, respectively.
       In \cite{St1}, Proposition 4.6 we have proved that if we know the coefficients of the
       differential operators $L_{j, \varepsilon}$ then we can directly determine the characteristic
       coefficients at every singular points. The restriction of Proposition 4.6 to the Heun type 
       reducible equation leads to
        \bpr{ch}
          The coefficients of the operators $L_{j, \varepsilon}$ in \eqref{pe} are unique determined 
          only by the characteristic exponents
           \ben
             L_{1, \varepsilon}
              &=&
             \partial - \left(\frac{\rho^R_1}{x-x_R} + \frac{\rho^L_1}{x-x_L}
                             + \frac{\rho^{RR}_1}{x-x_{RR}} + \frac{\rho^{LL}_1}{x-x_{LL}}\right)\\[0.2ex]
              L_{2, \varepsilon}
              &=&
              \partial - \left(\frac{\rho^R_2-1}{x-x_R} + \frac{\rho^L_2-1}{x-x_L}
              + \frac{\rho^{RR}_2-1}{x-x_{RR}} + \frac{\rho^{LL}_2-1}{x-x_{LL}}\right)\,.                        
           \een
        \epr
  Thanks to \prref{ch}
    the characteristic exponents $\rho^j_i$ and $\rho^{jj}_i,\,i=1, 2,\,j=R, L$
    are
    \ben
    \rho^R_1=\frac{\beta_1}{2 \sqrt{\varepsilon}},\,
    \rho^R_2=\frac{\beta_2}{2 \sqrt{\varepsilon}};
    & &
    \rho^L_1=-\frac{\beta_1}{2 \sqrt{\varepsilon}},\,
    \rho^L_2=-\frac{\beta_2}{2 \sqrt{\varepsilon}};\\[0.4ex]
    \rho^{RR}_1=-\frac{\gamma_1}{2 \sqrt{\varepsilon}},\,
    \rho^{RR}_2=1-\frac{\gamma_2}{2 \sqrt{\varepsilon}};
    & &
    \rho^{LL}_1=\frac{\gamma_1}{2 \sqrt{\varepsilon}},\,
    \rho^{LL}_2=1+\frac{\gamma_2}{2 \sqrt{\varepsilon}}\,.
    \een
    The exponents differences $\Delta^j_{12}=\rho^j_1-\rho^j_2$ and $\Delta^{jj}_{12}=\rho^{jj}_1-\rho^{jj}_2$ corresponding to the above characteristic
    exponents are defined as follows,
    \ben
    \Delta^R_{12}=\frac{\beta_1-\beta_2}{2 \sqrt{\varepsilon}},\,
    & &
    \Delta^L_{12}=-\frac{\beta_1-\beta_2}{2 \sqrt{\varepsilon}},\\[0.5ex]
    \Delta^{RR}_{12}=-1-\frac{\gamma_1-\gamma_1}{2 \sqrt{\varepsilon}},\,
    & &	
    \Delta^{LL}_{12}=-1 + \frac{\gamma_1-\gamma_2}{2 \sqrt{\varepsilon}}\,.
    \een
    The local theory of the Fuchsian singularities ensures a necessary condition
    for the existence of the logarithmic term near the singular point $x_j$ or $x_{jj}$ in the solution of
    the Heun type equation. If the exponent difference $\Delta^j_{12}\in\NN_0$ or $\Delta^{jj}_{12}\in\NN_0$
    then the Heun type equation always admits a particular solution $\Phi_1(x, \varepsilon)$ in the
    form
     $$\,   
     \Phi_1(x, \varepsilon)=(x-x_j)^{\rho^j_1}\,h_{1, j}(x),\quad\textrm{or}\quad
           \Phi_1(x, \varepsilon)=(x-x_{jj})^{\rho^{jj}_1}\,h_{1, jj}(x)\,,
     \,$$
  where $h_{1, j}(x)$ and $h_{1, jj}(x)$ are holomorphic functions in a neighborhood of the point $x_j$ and $x_{jj}$, respectively. 
  But the second
  element $\Phi_{12}(x, \varepsilon)$ of the fundamental set of solutions near the point $x_j$ or $x_{jj}$ can
  contain logarithmic term. Classically, such a Fuchsian singularity $x_j$ (resp. $x_{jj}$) for which $\Delta^j_{12}\in\NN_0$
  (resp. $\Delta^{jj}_{12}\in\NN_0$)
  is called a resonant Fuchsian singularity \cite{I-Y}. In this paper we restrict our attention to these families of
  the perturbed equation for which there are exactly two different types of resonant Fuchsian singularities --
  one of the type $x_j, j=R, L$ and one of the type $x_{jj}, j=R, L$. Through this paper we call these
  values of the parameters $\beta_j, \gamma_j, \varepsilon$ for which the perturbed equation has
  two different types of resonant Fuchsian singularities just a double resonance.  The motivation of studying
  the Heun type equation with double resonances naturally arises from our previous work \cite{St1}
  where we have shown that during a resonance the Stokes matrices of the initial equation are connected by
  a limit with the parts of the monodromy matrices of the perturbed equation when $\sqrt{\varepsilon} \rightarrow 0$.
  In the present paper we study the double resonances with the purpose to connect again by a limit both Stokes
  matrices of the initial equation with suitable monodromy matrices of the perturbed equation.

    We finish this paragraph by introducing the notion of the monodromy matrices following the paper of
    Dubrovin and Mazzocco \cite{BD-MM}.
    Consider the Heun type equation over $X=\CC\PP^1 -\{x_R, x_L, x_{RR}, x_{LL}\}$.    
    Let $\Phi(x, \varepsilon)$ be a fundamental matrix of the Heun type equation.
    It is a multivalued analytic function on the punctured Riemann sphere $X$. This
    multivaluedness is described by the monodromy matrices. Let $x_0\in X$ be a point
    that does not lie on the same line with three of the singular points. Note that three
    of the singular points (therefore all the points) lie on the same line if and only if
    either $\sqrt{\varepsilon}\in\RR$ or $\sqrt{\varepsilon}\in i\RR$. In fact the so called from us
    double resonance holds when $\sqrt{\varepsilon}\in\RR$. Let 
    $\gamma_R, \gamma_L, \gamma_{RR}, \gamma_{LL}$  be  simple closed loops starting and ending at 
    the point $x_0$, going around the singular points
   $x_R, x_L, x_{RR}, x_{LL}$, respectively, in positive direction and not crossing each other. The so chosen
   loops fix a basis in the fundamental group $\pi_1(X, x_0)$ with base point at $x_0$ of
   the punctured Riemann sphere $X$. 
    Let $\Phi_{\gamma}(x, \varepsilon)$ be the result of the analytic continuation
    of the fundamental matrix $\Phi(x, \varepsilon)$ along the loop $\gamma\in \pi_1(x, x_0)$. Since the
    Heun type equation is a linear equation, the matrix $\Phi_{\gamma}(x, \varepsilon)$ is also a fundamental matrix 
    of the same equation. Therefore there is a unique invertable constant matrix
    $M_{\gamma}(\varepsilon)\in GL_2(\CC)$ such that
     $$\,
       \Phi_{\gamma}(x, \varepsilon)=\Phi(x, \varepsilon)\,M_{\gamma}\,.
     \,$$ 
     The matrix $M_{\gamma}$ depends only on the homotopy class $[\gamma]$ of the loop $\gamma$.
     \bde{mr}
      The antihomomorphiism mapping
        \ben
         & &
         \pi_1(X, x_0) \longrightarrow GL_2(\CC)\\[0.1ex]
         & &
         [\gamma] \longrightarrow M_{\gamma}(\varepsilon)\\[0.1ex]
         & &
         M_{\gamma_1\,\gamma_2}(\varepsilon)=M_{\gamma_2}(\varepsilon)\,M_{\gamma_1}(\varepsilon),\quad
         M_{\gamma^{-1}}(\varepsilon)=M^{-1}_{\gamma}(\varepsilon)         
        \een
        determines monodromy representation of the Heun type equation with respect to the given
        fundamental matrix.
     \ede
    
     \bde{mm}
      The images $M_j(\varepsilon)=M_{\gamma_j}(\varepsilon)$ of the generators 
      $\gamma_R, \gamma_L, \gamma_{RR}, \gamma_{LL}$ of $\pi_1(X, x_0)$ under the monodromy representation
      are called the monodromy matrices of the Heun type equation. 
      \ede
    
    In Section 4, \thref{expl-A} we will present explicitly fundamental matrices with respect to which 
    we will  compute the corresponding monodromy matrices.

%%%%%%%%%%%%%%%%%%%%%%%%%%%%%%%%%%%%%%%%
%Symmetries
%%%%%%%%%%%%%%%%%%%%%%%%%%

	\subsection{Symmetries}

  The considered in this paper  equations
	 can be  rewritten  as second-order homogeneous equation 
	 \be\label{he}
	   y'' + b_1(x, \varepsilon)\,y' + 
		b_0(x, \varepsilon)\,y=0\,.
	\ee
	 The coefficients $b_j(x, 0)$ of the equation \eqref{npe} - \eqref{op-ini} are
	 \ben
	   b_1(x, 0) &=&
		-\frac{\alpha_1 + \alpha_2}{x} - \frac{\beta_1+\beta_2}{x^2} - (\gamma_1 + \gamma_2)\,,\\[0.3ex]
			b_0(x, 0) &=&
				\frac{\alpha_1\,\gamma_2 + \alpha_2\,\gamma_1}{x} 
				+\frac{\alpha_1 + \alpha_1\,\alpha_2 + \beta_1\,\gamma_2 + \beta_2\,\gamma_1}{x^2} +
				\\[0.3ex]
					&+&
			\frac{2 \beta_1 + \alpha_1\, \beta_2 + \alpha_2\, \beta_1}{x^3}
			+ \frac{\beta_1\, \beta_2}{x^4} + \gamma_1\,\gamma_2\,.
	\een
  The equation \eqref{npe} - \eqref{op-ini}	with coefficients $b_j(x, 0)$
 is invariant under the following transformations:
 \ben
             & &
  \alpha_1 \longrightarrow - \alpha_1\,,\,\,\,
	\alpha_2 \longrightarrow - \alpha_2 - 2\,,\,\,\,
	\beta_j \longrightarrow - \gamma_j\,,\,\,\,
	\gamma_j \longrightarrow - \beta_j\,,\,\,\, 
	x \longrightarrow \frac{1}{x}\,,\,\,\,
	y \longrightarrow y,\\
           & &
  \alpha_1 \longrightarrow - \alpha_1\,,\,\,\,
	\alpha_2 \longrightarrow - \alpha_2 - 2\,,\,\,\,
	\beta_j \longrightarrow  \gamma_j\,,\,\,\,
	\gamma_j \longrightarrow  \beta_j\,,\,\,\, 
	x \longrightarrow - \frac{1}{x}\,,\,\,\,
	y \longrightarrow y,\\
           & &
  \alpha_1 \longrightarrow  \alpha_1\,,\,\,\,
	\alpha_2 \longrightarrow  \alpha_2 \,,\,\,\,
	\beta_j \longrightarrow  - \beta_j\,,\,\,\,
	\gamma_j \longrightarrow  - \gamma_j\,,\,\,\, 
	x \longrightarrow -  x\,,\,\,\,
	y \longrightarrow y.
\een

Then we have the following obvious result.
 \bpr{p1}
 Let $\phi(x)$  be  a particular solution near $x=0$ of the equation \eqref{npe} - \eqref{op-ini} with parameters 
  $\alpha_1=0,\,\,\alpha_2=-2$ and $\beta_j,\,\gamma_j,\,j=1,2$ are
	arbitrary.  Then 
      \begin{enumerate}
          \item\,
          $\phi(x^{-1})$ is a particular solution near $x=\infty$
	of the equation \eqref{npe} - \eqref{op-ini} with parameters $\alpha_1=\alpha_2=0$ and
	$\beta_j=-\gamma_j,\,j=1,2$.
            \item\,
    $\phi(-x^{-1})$ is a particular solution near $x=\infty$
	of the equation \eqref{npe} - \eqref{op-ini} with parameters $\alpha_1=\alpha_2=0$ and $\beta_j=\gamma_j$.
             \item\,
       $\phi(-x)$ is a particular solution near the origin of the equation \eqref{npe} - \eqref{op-ini} with parameters 
         $\alpha_1=0,\,\alpha_2=-2$ and $-\beta_j,\,-\gamma_j$.
        \end{enumerate}
\epr

  When $\alpha_1=0,\,\alpha_2=-2$ the transformation
   $$\,
     x \longrightarrow \frac{1}{x}
   \,$$
   takes the operators $L_{J, \varepsilon}$ from \eqref{pe-par} of the  Heun type equation  into operators of a Heun type equation 
  \be\label{the}\qquad
    L_{1, \varepsilon}
       &=&
       \partial - \frac{\gamma_1}{2 \sqrt{\varepsilon}}
       \left(-\frac{1}{x-\sqrt{\varepsilon}}+ \frac{1}{x+\sqrt{\varepsilon}}\right)
    -\frac{\beta_1}{2 \sqrt{\varepsilon}}\left(
        \frac{1}{x-\frac{1}{\sqrt{\varepsilon}}} -
        \frac{1}{x+\frac{1}{\sqrt{\varepsilon}}}\right),\\[0.35ex]
      L_{2, \varepsilon}
      &=&
      \partial + \frac{\gamma_2}{2 \sqrt{\varepsilon}}
      \left(\frac{1}{x-\sqrt{\varepsilon}} - \frac{1}{x+\sqrt{\varepsilon}}\right)
      + \left(1-\frac{\beta_2}{2 \sqrt{\varepsilon}}\right)
      \frac{1}{x-\frac{1}{\sqrt{\varepsilon}}} +\nonumber\\[0.35ex]
      & +& 
      \left(1+ \frac{\beta_2}{2 \sqrt{\varepsilon}}\right)
      \frac{1}{x+\frac{1}{\sqrt{\varepsilon}}}.\nonumber
  \ee
  Thus the Heun type equation  is invariant under the following transformation
   $$\,
    \beta_1 \longrightarrow -\gamma_1,\,\,\gamma_1 \longrightarrow  - \beta_1,\,\,   
    1 -\frac{\beta_2}{2 \sqrt{\varepsilon}} \longrightarrow \frac{\gamma_2}{2 \sqrt{\varepsilon}},\,\,
      1 +\frac{\beta_2}{2 \sqrt{\varepsilon}} \longrightarrow -\frac{\gamma_2}{2 \sqrt{\varepsilon}},\,\,
    x \longrightarrow \frac{1}{x}, y \longrightarrow y\,.
   \,$$

%%%%%%%%%%%%%%%%%%%%%%%%%%%%%%%%%%%%%%%%%%%%%%%
% DCHE
%%%%%%%%%%%%%%%%%%%%%%%%%%%%%%%%%%%%%%%%%%%%%%%

  \section {The initial equation}

  In this section we will compute the Stokes matrices of the initial equation
	at  $x=0$ and $x=\infty$. 
	In paragraph 2.2 we introduced a global fundamental matrix of the initial equation.
	Recall that its element $\Phi_{12}(x, 0)$ is represented as an integral
	 $$\,
	   \Phi_{12}(x, 0)=\Phi_1(x, 0)\,\int_{\Gamma(x, 0)}
	     \frac{\Phi_2(z, 0)}{\Phi_1(z, 0)}\, d z\,,
	 \,$$
	 where $\Phi_j(x, 0),\,j=1, 2$ are the solutions of the equations $L_{j, 0}(u)=0$.
	 In this paper we use two different fundamental matrices with respect to which we will compute 
	 the Stokes matrices: the matrix $\Phi_0(x, 0)$ at $x=0$ and the matrix $\Phi_{\infty}(x, 0)$
	 at $x=\infty$. The difference between these matrix solutions is the path of integration
	 $\Gamma(x, 0)$. Choosing 
	  $$\,
	 \Phi_1(x, 0)=e^{\gamma_1 x}\,e^{-\frac{\beta_1}{x}},\quad
	 	 \Phi_2(x, 0)=\frac{e^{\gamma_2 x}\,e^{-\frac{\beta_2}{x}}}{x^2}\,,
	  \,$$
	  the element $\Phi_{12}(x, 0)$ becomes
	  $$\,
	  	   \Phi_{12}(x, 0)=\Phi_1(x, 0)\,\int_{\Gamma(x, 0)}
	  	   \frac{e^{(\gamma_2-\gamma_1)\,z}\,e^{-\frac{\beta_2-\beta_1}{z}}}{z^2} d z\,.
	  	   \,$$
	   When $\Phi_{12}(x, 0)$ is an element of the matrix $\Phi_0(x, 0)$ the path $\Gamma(x, 0)$ is a path
	   from $0$ to $x$, approaching $0$ in the direction $\theta=\arg (\beta_2-\beta_1)$. When $\Phi_{12}(x, 0)$
	   is an element of the matrix $\Phi_{\infty}(x, 0)$ the path $\Gamma(x, 0)$ is a path from
	   $+\infty e^{i \theta}$ to $x$, approaching $+\infty e^{i \theta}$ in the direction $\theta=\arg (\gamma_1-\gamma_2)$.
        
        \bre{res-in}
	    When $\beta_1=\beta_2$ the integral
	    $$\,
	     \int_0^x \frac{e^{(\gamma_2-\gamma_1)\,z}}{z^2}\, d z
	    \,$$
	    does not exists. For this reason in this papper we study the initial and perturbed equations under assumption
	    that 
			 $$\,
			   \beta_1 \neq \beta_2 \quad \textrm{and}\quad \gamma_i,\,i=1, 2\quad \textrm{are arbitrary}\,;
			\,$$
       In fact, when $\beta_1=\beta_2$ the point $x=0$ is a resonant irregular singularity for the initial equation.
       To study such an irregular singularity we have to use some other representation at the origin of the solution
       $\Phi_{12}(x, 0)$. It seems that an integral representation but with a base point different from $0$ to be
       useful. This probem, as well as its perturbed analog are postponed for further researches.     
         \ere

%%%%%%%%%%%%%%%%%%%%%%%%%%%%%%%%%%%%%%%%%%%%%%%
%A
%%%%%%%%%%%%%%%%%%%%%%%%%%%%%%%%%%%%%%%%%%%%%%

 The  integral representation of the element $\Phi_{12}(x, 0)$ is unsuitable for the computation of the Stokes matrices.
 On the other hand, as we show bellow,  this integral is easily represented as formal series at the origin,
 as well at the infinity point. Then applying summability theory we lift these series to the corresponding
 actual solutions. With respect to these actual solutions we can easily compute the Stokes multipliers 
 $\mu^{\theta}_j, j=0, \infty$. 

 We begin by deducing  formal fundamental matrices at the origin and at the infinity point from the  integrals.

\bpr{p2}
 Assume that $\beta_1 \neq \beta_2$. Then the initial equation possesses
 an unique formal fundamental 
	matrix $\hat{\Phi}_0(x, 0)$ at the origin in the form
	 $$\,
	  \hat{\Phi}_0(x, 0)=\exp (G x)\,\hat{H}(x)\,x^{\Lambda}\,exp\left(-\frac{B}{x}\right)\,,
	 \,$$
	 where
	 $$\,
	   G=\diag(\gamma_1, \gamma_2),\quad
	     \Lambda=\diag(0, -2),\quad B=\diag(\beta_1, \beta_2)\,.
	 \,$$
	 The matrix $\hat{H}(x)$ is given by
	  \be\label{H}
	    \hat{H}(x)=\left(\begin{array}{cc}
	    	1    &\frac{x^2\,e^{(\gamma_2-\gamma_1) x} }{\beta_2-\beta_1}
	    	-\frac{(\gamma_2-\gamma_1) x^3 \hat{\psi}(x)}{\beta_2-\beta_1}\\[0.15ex]
	    	0    &1
	    	           \end{array}
	    	        \right)\,.    
	  \ee
	
	 The element $\hat{\psi}(x)$ is represented by a formal power series
	 \be\label{psi}
	  \hat{\psi}(x)=\sum_{k=1}^{\infty} \frac{(-1)^k \,k!\,S_{k-1}}{(\beta_2-\beta_1)^k}\,x^k 
		\ee
	where $S_{k-1},\,k\in\NN$ are the partial sums of the following infinite absolutely convergent series
	of numbers
	 \be\label{S}
	  S=\sum_{n=0}^{\infty} 
		\frac{(-1)^{n+1}\,(\gamma_2-\gamma_1)^n\,(\beta_2-\beta_1)^n}{n!\,(n+1)!}\,.
	\ee
	Moreover,
	\begin{enumerate}
		\item\,
		If $S \neq 0$ but $\gamma_1 \neq \gamma_2$ then the formal power series $\hat{\psi}(x)$ is divergent.
		
		\item\,
		If $S=0$ but $\gamma_1 \neq \gamma_2$ then the power series $\hat{\psi}(x)$ represents an analytic in $\CC$ function.
		
		\item\,
		If $\gamma_1=\gamma_2$ then $\hat{\psi}(x)=\frac{x}{\beta_2-\beta_1}$ but the matrix $\hat{H}(x)$ has the form
		\be\label{H-g}
		      \hat{H}(x)=\left(\begin{array}{cc}
		      	1    &\frac{x^2}{\beta_2-\beta_1}\\[0.15ex]
		      	0    &1
		      \end{array}
		      \right)\,.    
		\ee
		
	\end{enumerate}

\epr    	 

\proof 
 We have already chosen $\hat{\Phi}_1(x)$ and $\hat{\Phi}_2(x)$.
Then if $\gamma_1 \neq \gamma_2$ for $\hat{\Phi}_{12}(x)$ we have
  \ben
	  \hat{\Phi}_{12}(x)  &=&
		\hat{\Phi}_1(x)\,\int_0^x \frac{\hat{\Phi}_2(z)}{\hat{\Phi}_1(z)}\,d z=
		e^{\gamma_1 x}\,e^{-\frac{\beta_1}{x}}\,
		\int_0^x \frac{e^{(\gamma_2-\gamma_1) z}\,e^{-\frac{\beta_2-\beta_1}{z}}}{z^2}\, d z\\[0.4ex]
		          &=&
		\frac{e^{\gamma_2 x}\,e^{-\frac{\beta_2}{x}}}{\beta_2-\beta_1} -
		\frac{\gamma_2-\gamma_1}{\beta_2-\beta_1}\,e^{\gamma_1 x}\,e^{-\frac{\beta_1}{x}}\,
			\int_0^x e^{(\gamma_2-\gamma_1) z}\,e^{-\frac{\beta_2-\beta_1}{z}}\, d z				
	\een
where  the integral is taken in the direction $\arg (\beta_2-\beta_1)$. In particular case when $\gamma_1=\gamma_2$ 
the element $\hat{\Phi}_{12}(x, 0)$ becomes
  $$\,
    \hat{\Phi}_{12}(x, 0)=e^{\gamma_1 x}\,e^{-\frac{\beta_1}{x}}
    \int_0^x \frac{e^{-\frac{\beta_2-\beta_1}{z}}}{z^2}\, d z=
    \frac{e^{\gamma_1 x}\,e^{-\frac{\beta_2}{x}}}{\beta_2-\beta_1}\,,
  \,$$
  where the integral is again taken in the direction $\arg (\beta_2-\beta_1)$.
  
  Let $\gamma_1 \neq \gamma_2$.
The Taylor series at $z=0$ of the function $e^{(\gamma_2-\gamma_1) z}$ is
$$\,
   e^{(\gamma_2-\gamma_1) z}=\sum_{k=0}^{\infty} \frac{(\gamma_2-\gamma_1)^k}{k!}\,z^k\,.
\,$$
Then it is not so difficult to show that 
 \ben
    & &
  \int_0^x e^{(\gamma_2-\gamma_1) z}\,e^{-\frac{\beta_2-\beta_1}{z}}\,d z=
	            	\sum_{k=0}^{\infty} \frac{(\gamma_2-\gamma_1)^k}{k!}\,
	\int_0^x z^k\,e^{-\frac{\beta_2-\beta_1}{z}}\, d z=\\[0.4ex]
	            &=&
  (\beta_2-\beta_1)\,S\,
	\int_0^x \frac{e^{-\frac{\beta_2-\beta_1}{z}}}{z}\, d z
	- x\,e^{-\frac{\beta_2-\beta_1}{x}}\,S\,
				\sum_{k=0}^{\infty} \frac{(-1)^k k!}{(\beta_2-\beta_1)^k}\,x^k+\\[0.4ex]
							&+&
	  x\,e^{-\frac{\beta_2-\beta_1}{x}}\,
									\sum_{k=1}^{\infty} \frac{(-1)^k k!\,S_{k-1}}{(\beta_2-\beta_1)^k}\,x^k\,,
\een
where $S$ is defined by the series
$$\,
    \sum_{n=0}^{\infty} 
		\frac{(-1)^{n+1}\,(\gamma_2-\gamma_1)^n\,(\beta_2-\beta_1)^n}{n!\,(n+1)!}\,,   
\,$$
and $S_{k-1},\,k\in\NN$ are its partial sums. This series of numbers is absolutely
convergent and we denote by $S$ its sum.

It turns out that the sum of the first two addends in the above representation
of the integral $\int_0^x e^{(\gamma_2-\gamma_1) z}\,e^{-(\beta_2-\beta_1)/z}\, d z$ is equal to zero.
Indeed, consider the integral 
  $$\,
	  h(x)=e^{-\frac{\beta_1}{x}}\,\int_0^x \frac{e^{-\frac{\beta_2-\beta_1}{z}}}{z}\,d z\,.
	\,$$
	Let us present it in the following form
	$$\,
	  h(x)=e^{-\frac{\beta_2}{x}}\,
		\int_0^x \frac{e^{(\beta_2-\beta_1)\left(\frac{1}{x}-\frac{1}{z}\right)}}{z}\, d z\,.
	\,$$
	Introducing a new variable $\zeta$ via 
	$$\,
	  (\beta_2-\beta_1)\left(\frac{1}{x}-\frac{1}{z}\right)=-\frac{\zeta}{x}
	\,$$
	the integral $h(x)$ gets into
	$$\,
	  h(x)=e^{-\frac{\beta_2}{x}}\,\int_0^{\infty(\tau)}
		\frac{e^{-\frac{\zeta}{x}}}{\zeta + \beta_2 - \beta_1}\,d \zeta
	\,$$
	where $\tau=\arg (\beta_2-\beta_1)$.
	The analytic continuation of $h(x)$ on the
	$x$-plane yields an analytic in $x$ function
	$$\,
	   h_{\theta}(x)=e^{-\frac{\beta_2}{x}}\,
	   \int_0^{+\infty e^{i \theta}} \frac{e^{-\frac{\zeta}{x}}}{\zeta+\beta_2-\beta_1}\, d \zeta 
	\,$$
	in all directions $\theta$  except along the direction
	$\arg (\beta_1-\beta_2)$.
	
	Next, consider the formal power series 
	 $$\,
	\hat{f}(x)=\sum_{k=0}^{\infty} \frac{(-1)^k\,k!}{(\beta_2-\beta_1)^k}\,x^k\,,
	\,$$
    This series belongs to the class of 1-summable series from \leref{example} for $2+\alpha_2-\alpha_1=1$
    and $\beta_j-\beta_i=\beta_2-\beta_1$.
     Thus for any direction $\theta\neq \arg (\beta_2-\beta_1)$
	the function 
	 \ben
	 f_{\theta}(x) &=&
	  \int_0^{+ \infty e^{i \theta}}
	  \left(1 + \frac{\zeta}{\beta_2-\beta_1}\right)^{-1}\,e^{-\frac{\zeta}{x}}\,
	  d \left(\frac{\zeta}{x}\right)=\\[0.2ex]
	        &=&
	 \frac{\beta_2-\beta_1}{x}
	   \int_0^{+\infty\,e^{i \theta}}
		 \frac{e^{-\frac{\zeta}{x}}}{\beta_2-\beta_1+\zeta}\,d \zeta
	\een
	defines its 1 - sum in such a direction.
	Then, the explicit form of the element $\Phi_{12}(x)$ follows from
	a combination of the above observations on the integral $h(x)$ and the formal series $\hat{f}(x)$.

  Finally, the formal power series
    $$\,
		  \hat{\psi}(x)=\sum_{k=1}^{\infty} \frac{(-1)^k\,k!\,S_{k-1}}{(\beta_2-\beta_1)^k}\,x^k=
			 \sum_{k=1}^{\infty} b_k\,x^k
		\,$$
		is convergent only at $x=0$ when $S\neq 0$. Indeed, since the series \eqref{S} is convergent (and $S$ is its sum),
	 its sequence of partial sums $\{S_{k-1}\}$ satisfy the following limits
		$$\,
			\lim_{k \rightarrow \infty} S_{k-1}=S \neq 0\,.																
		\,$$
		Then, for the radius of convergence $R$ for the power series $\hat{\psi}(x)$, we have
		$$\,
		  R=\lim_{k \rightarrow \infty} \left|\frac{b_k}{b_{k+1}}\right|=
			   \lim_{k \rightarrow \infty} \frac{1}{k+1}
				\left|\frac{(\beta_2-\beta_1)\,S_{k-1}}{S_k}\right|=
				\lim_{k \rightarrow \infty} \frac{|\beta_2-\beta_1|}{k+1}=0\,.
		\,$$
		
		Let now $S=0$. Then $S_{k-1}$ and $S_k$ can be rewritten as
		 \ben
		   S_{k-1} &=&
		   -\sum_{n=k}^{\infty} \frac{(-1)^{n+1}\,(\gamma_2-\gamma_1)^n\,(\beta_2-\beta_1)^n}{n!\,(n+1)!}=\\[0.2ex]
		         &=&
		    -\frac{(-1)^{k+1}\,(\gamma_2-\gamma_1)^k\,(\beta_2-\beta_1)^k}{k!\,(k+1)!}
		    \sum_{n=0}^{\infty}
		    \frac{(-1)^n (\gamma_2-\gamma_1)^n\,(\beta_2-\beta_1)^n}{(k+1)^{(n)}\,(k+2)^{(n)}}\,,\\[0.3ex]
		       S_k &=&
		       -\sum_{n=k+1}^{\infty} \frac{(-1)^{n+1}\,(\gamma_2-\gamma_1)^n\,(\beta_2-\beta_1)^n}{n!\,(n+1)!}=\\[0.2ex]
		       &=&
		       -\frac{(-1)^{k+2}\,(\gamma_2-\gamma_1)^{k+1}\,(\beta_2-\beta_1)^{k+1}}{(k+1)!\,(k+2)!}
		       \sum_{n=0}^{\infty}
		       \frac{(-1)^n (\gamma_2-\gamma_1)^n\,(\beta_2-\beta_1)^n}{(k+2)^{(n)}\,(k+3)^{(n)}}\,.
		 \een
	 Then for the the radius of convergence $R$ we obtain
		  \ben
		    R=\lim_{k \rightarrow \infty}
		    \frac{|\beta_2-\beta_1|}{k+1}
		    \frac{|\gamma_2-\gamma_1|^k\,|\beta_2-\beta_1|^k}{k!\,(k+1)!}
		    \frac{(k+1)!\,(k+2)!}{|\gamma_2-\gamma_1|^{k+1}\,|\beta_2-\beta_1|^{k+1}}
		    \frac{\sum_{n=0}^{\infty}
		    	\frac{(-1)^n (\gamma_2-\gamma_1)^n\,(\beta_2-\beta_1)^n}{(k+1)^{(n)}\,(k+2)^{(n)}}}
		    {\sum_{n=0}^{\infty}
		    	\frac{(-1)^n (\gamma_2-\gamma_1)^n\,(\beta_2-\beta_1)^n}{(k+2)^{(n)}\,(k+3)^{(n)}}}\,.
		  \een
		  Since each of the infinite series has a limit 1 when $k \rightarrow \infty$ for $R$ we find that
		   $$\,
		     R=\lim_{k \rightarrow \infty} \frac{k+2}{|\gamma_2-\gamma_1|}=\infty\,.
		   \,$$
		   Thus the series $\hat{\psi}(x)$ is convergent in $\CC$.
    	 Unfortunately till now we can not show explicitly which is the analytic in $\CC$ function, that has
    	 $\hat{\psi}(x)$ as a series expansion around $x=0$.

		This completes the proof.
\qed

  Recall that the Bessel function of the first kind $J_{\nu}(z)$ of order $\nu$ can be
  defined
  by its series expansion around $z=0$, \cite{W}, 
  \be\label{bes}
  & &
  J_{\nu}(z)=\left(\frac{z}{2}\right)^{\nu}
  \sum_{k=0}^{\infty} \frac{(-1)^k (z/2)^{2 k}}{k!\,\Gamma(n+k+1)}\,,
  \quad \textrm{if} \quad -\nu \notin \NN\,,\\[0.2ex]
  & &
  J_{-\nu}(z)=(-1)^{\nu}\,J_{\nu}(z)\,,\quad
  \textrm{if} \quad-\nu\in\NN\,.\nonumber
  \ee
  Here $\Gamma(x)$ is the classical Gamma function (see \cite{HB-AE}).
  It is well known that for  $\nu > -1$ the Bessel function $J_{\nu}(z)$ has 
  an infinite number of zeros at that only real zeros.
  In particular, if $\nu=1$ the first five zeros of $J_1(z)$ are
  $z=0,\,z=3, 8317,\,z=7, 0156,\,z=10,1735,\,z=13,3237$.

\bre{r1}
  The condition $S \neq 0$ implies that the number
	$z=2 \sqrt{(\gamma_2-\gamma_1)\,(\beta_2-\beta_1)}$ is  not a zero of the
	Bessel function $J_1(z)$.
\ere

\bre{for}
 The formal solution $\hat{\Phi}_{12}(x)$ can be obtained as a formal particular solution
 near $x=0$ of the following non-homogenous linear ODE
    \be\label{12}
     y'(x) - \left(\frac{\beta_1}{x} + \gamma_1\right)\,y(x)=
     \frac{e^{\gamma_2 x}\,e^{-\frac{\beta_2}{x}}}{x^2}.
  \ee
 Looking for $\hat{\Phi}_{12}(x)$ in the form $e^{\gamma_2 x}\,e^{-\frac{\beta_2}{x}}\,\hat{y}(x)$,
 we find that $\hat{y}(x)$ must satisfy the equation
  $$\,
    x^2\,\hat{y}'(x) + \left[(\beta_2-\beta_1) + (\gamma_2-\gamma_1)\,x^2\right]\,\hat{y}(x)=1.
  \,$$
 The last equation admits an unique formal solution
  $$\,
     \hat{y}(x)=\sum_{k=0}^{\infty} a_k\,x^k,
  \,$$
  provided that $\beta_2 \neq \beta_1$. The coefficients $a_k$ are recursively constructed from
  the equation
   $$\,
      (k-1)\,a_{k-1} + (\beta_2-\beta_1)\,a_k + (\gamma_2-\gamma_1)\,a_{k-2}=0,\qquad
         k \geq 2,
   \,$$
   as $a_0=\frac{1}{\beta_2-\beta_1},\,a_1=0$.  Note that we are not able to give an explicit formula for the coefficients $a_k$.
    In particular,
   \ben
                     & &
        a_2=-\frac{\gamma_2-\gamma_1}{(\beta_2-\beta_1)^2},\quad
       a_3=\frac{2 (\gamma_2-\gamma_1)}{(\beta_2-\beta_1)^3},\quad
      a_4=-\frac{3!\,(\gamma_2-\gamma_1)}{(\beta_2-\beta_1)^4} + \frac{(\gamma_2-\gamma_1)^2}{(\beta_2-\beta_1)^3},\\[0.25ex]
                  & &
         a_5=\frac{4!\,(\gamma_2-\gamma_1)}{(\beta_2-\beta_1)^5} - \frac{6 (\gamma_2-\gamma_1)^2}{(\beta_2-\beta_1)^4},\quad
           a_6=-\frac{5!\,(\gamma_2-\gamma_1)}{(\beta_2-\beta_1)^6} + \frac{6.3! (\gamma_2-\gamma_1)^2}{(\beta_2-\beta_1)^5}-
              \frac{(\gamma_2-\gamma_1)^3}{(\beta_2-\beta_1)^4}, \ldots
  \een
As a result $\hat{\Phi}_{12}(x)$ becomes
   $$\,
      \hat{\Phi}_{12}(x)=\frac{e^{\gamma_2 x}\,e^{-\frac{\beta_2}{x}}}{\beta_2-\beta_1}
      + x\,e^{\gamma_2 x}\,e^{-\frac{\beta_2}{x}} \sum_{k=1}^{\infty} a_{k+1}\,x^k.
   \,$$
   Now we represent $e^{\gamma_2 x}$ as
   $$\,
       e^{\gamma_2 x}=e^{\gamma_1 x}\,e^{(\gamma_2-\gamma_1) x}=
       e^{\gamma_1 x} \sum_{k=0}^{\infty} \frac{(\gamma_2-\gamma_1)^k}{k!}\,x^k.
   \,$$
   Then we can rewrite $\hat{\Phi}_{12}(x)$ as
       $$\,
       \hat{\Phi}_{12}(x)=\frac{e^{\gamma_2 x}\,e^{-\frac{\beta_2}{x}}}{\beta_2-\beta_1}
       + x\,e^{\gamma_1 x}\,e^{-\frac{\beta_2}{x}} \sum_{k=1}^{\infty} c_k\,x^k,
       \,$$
       where the formal series $\sum c_k\,x^k$ is the product of the series
       $\sum \frac{(\gamma_2-\gamma_1)^k}{k!} x^k$ and $\sum a_{k+1} x^k$. In particular,
       $$\,
          c_k=\sum_{s=0}^{k-1} b_s\,a_{k+1-s},\qquad
          b_s=\frac{(\gamma_2-\gamma_1)^s}{s!}.
       \,$$
       So, the series $\sum_{k=1}^{\infty} c_k\,x^k$ is nothing but the series
       $-\frac{\gamma_2-\gamma_1}{\beta_2-\beta_1}\,\hat{\psi}(x)$ from \prref{p2}.

      We give here this method to illustrate the difficulty of the construction of a suitable particular formal
       solution even in the case of second order reducible equation. Fortunately, in our case we have an
      integral expression of $\Phi_{12}(x)$ whose formal representation near the orign gives us a formal solution
      $\hat{\Phi}_{12}(x)$.
\ere

     Thanks to \prref{p1} we construct a formal fundamental matrix of the initial equation at the point $x=\infty$.

       \bpr{inf}
          Let $\beta_1 \neq \beta_2$. Then the initial equation possesses
          an unique formal fundamental  matrix $\hat{\Phi}_{\infty}(x, 0)$ at $x=\infty$ in the form
          $$\,
            \hat{\Phi}_{\infty}(x, 0)=\exp \left(-\frac{B}{x}\right)\,\left(\frac{1}{x}\right)^{-\Lambda}\,\hat{P}(x)\,\exp(G x)\,,
          \,$$
          where the matrices $B, \Lambda$ and $G$ are introduced by \prref{p2}.
          The matrix $\hat{P}(x)$ is defined as follows:
            \begin{enumerate}
            	
            \item\,
            If $\gamma_1=\gamma_2$ then
                   
                   \be\label{P-g}
                   \hat{P}(x)=\left(\begin{array}{cc}
                   	1     &\frac{e^{-\frac{\beta_2-\beta_1}{x}} -1}{\beta_2-\beta_1}\\[0.1ex]
                   	0     &1
                   \end{array}
                   \right)\,.     
                   \ee
                   
              \item\,
              If $\gamma_1 \neq \gamma_2$ then     
              \be\label{P}
               \hat{P}(x)=\left(\begin{array}{cc}
                   1     &-\frac{\hat{\varphi}(x)}{x}\\[0.1ex]
                   0     &1
                                \end{array}
                           \right)\,.     
             \ee
             The element $\hat{\varphi}(x)$ is defined by the formal series 
               
               \be\label{ffs-inf}
               \hat{\varphi}(x)=\sum_{k=1}^{\infty} \frac{k!\,S_{k-1}}{(\gamma_2-\gamma_1)^k}\,x^{-k}\,.
              \ee
                        The symbol $S_{k-1}$ is defined by \prref{p2}.
                        
           Moreover,
           \begin{enumerate}
           	\item\,
           	If $\lim_{k \rightarrow \infty} S_{k-1} \neq 0$ then the formal series $\hat{\varphi}(x)$ is divergent.
           	
           	\item\,
           	If $\lim_{k \rightarrow \infty} S_{k-1} =0$ then the series $\hat{\varphi}(x)$ defines an analytic in
           	$\CC\PP^1 - \{0\}$ function.
           \end{enumerate}
           \end{enumerate}             
      \epr

          \proof
        The elements $\hat{\Phi}_j(x),\,j=1, 2$ as the solutions of the equations $L_j\,y=0$ remain the same as in \prref{p2}.
       From \prref{p1} it follows that if $\hat{\Phi}_{12}(t)$ is a particular solution at $t=0$ of the initial equation with parameters
     $\alpha_1=\alpha_2=0,\,\beta_j=-\gamma_j$ then $\hat{\Phi}_{12}(x^{-1})$ is a particular solution near $x=\infty$
        of our initial equation. For the initial equation with $\alpha_1=\alpha_2=0,\,\beta_j=-\gamma_j$ the operators $L_j$ become
                 \be\label{tr-op}
                       L_j=\partial + \left(\frac{\gamma_j}{t^2} + \beta_j\right),\qquad
                       \partial=\frac{d}{d t}.
                 \ee
               Then the corresponding entries $\hat{\Phi}_j(t)$ are given by
                              $$\,
                                      \hat{\Phi}_j(t)=e^{-\beta_j t}\,e^{\frac{\gamma_j}{t}}.
                                 \,$$
                   Following the method of construction of $\hat{\Phi}_{12}(x)$ in \prref{p2}, we find that 
                   when $\gamma_1 \neq \gamma_2$
                     \ben
                      \hat{\Phi}_{12}(t)
                                       &=&
                       -e^{-\beta_1 t}\,e^{\frac{\gamma_1}{t}}
                                \int_0^t e^{-(\beta_2-\beta_1) z}\,e^{\frac{\gamma_2-\gamma_1}{z}}\,d z\\[0.2ex]
                                      &=&
                                -  t\,e^{\frac{\gamma_2}{t}}\,e^{-\beta_1 t}
                               \sum_{k=1}^{\infty} \frac{k!\,S_{k-1}}{(\gamma_2-\gamma_1)^k}\,t^k,
                   \een
                        where the integral is taken in the direction $\arg (\gamma_1-\gamma_2)$. 
                        The  sign minus comes from the observation  that
                     the transformation $x=1/t$ changes the equation 
                  $$\,
                        y'(x) - \left(\frac{\beta_1}{x} + \gamma_1\right)\,y(x)=\frac{e^{\gamma_2 x}\,e^{-\frac{\beta_2}{x}}}{x^2},
                    \,$$
                   which has $\hat{\Phi}_{12}$ as a particular solution, into the equation 
                       $$\,
                                   \dot{y}(t) + \left(\beta_1 + \frac{\gamma_1}{t^2}\right)\,y(t)=- e^{\frac{\gamma_2}{t}}\,e^{-\beta_2\,t}.
                  \,$$ 
 The symbol $S_{k-1}$ is defined by \prref{p2}. Then putting $t=x^{-1}$ in this formula 
 for $\hat{\Phi}_{12}(t)$ we obtain a particular solution near
                 $x=\infty$ of our equation.
                 
     In the same way as in the proof of \prref{p2} we can show that the series 
     $$\,
     \hat{\varphi}(x)=  \sum_{k=1}^{\infty} \frac{k!\,S_{k-1}}{(\gamma_2-\gamma_1)^k}\,x^{-k}
       \,$$
       is divergent when $\lim_{k \rightarrow \infty} S_{k-1} \neq 0$, and has a radius of convergence $R=\infty$
       when $\lim_{k \rightarrow \infty} S_{k-1}=0$.
       
        When $\gamma_1=\gamma_2=\gamma$ the element $\hat{\Phi}_{12}(t)$ becomes
         $$\,
           \hat{\Phi}_{12}(t)=-e^{-\beta_1 t}\,e^{\frac{\gamma}{t}}\,
           \int_0^t e^{-(\beta_2-\beta_1) z}\, d z=
           \frac{e^{-\beta_1 t}\,e^{\frac{\gamma}{t}}}{\beta_2-\beta_1}\,
           \left[e^{-(\beta_2-\beta_1) t} -1\right]\,,
         \,$$
         where  the integral is taken in the rial positive axis. Putting $t=x^{-1}$ in this formula we
         obtain $\hat{\Phi}_{12}(x)$ when $\gamma_1=\gamma_2$.
         
       This ends the proof.           
  \qed

  In keeping with the theory the divergent power series $\hat{\psi}(x)$ defined by \eqref{psi}
is 1 - summable. More precisely,

\ble{sum1}
  Let $\beta_2 \neq \beta_1,\,\gamma_2 \neq \gamma_1$. Assume also that $S\neq 0$, where $S$ is given by \eqref{S}.
  Then the formal power series $\hat{\psi}(x)$
	defined by \eqref{psi} is 1 - summable in any direction 
	$\theta \neq \arg(\beta_1 - \beta_2)$.  
      The function
	 \be\label{s2}
	     \psi_{\theta}(x)=
	  	-\frac{\beta_2-\beta_1}{x}\int_0^{+\infty\,e^{i \theta}}
		 \frac{u(\zeta)\,e^{-\frac{\zeta}{x}}}{\zeta+\beta_2-\beta_1}\,d \zeta 
		 +\frac{e^{(\gamma_2-\gamma_1)\,x}-1}{(\gamma_2-\gamma_1)\,x}\,,									
	\ee
	where $u(\zeta)$ is the analytic in $\CC$ function
	 $$\,
	  u(\zeta)=\sum_{k=0}^{\infty} \frac{(\gamma_2-\gamma_1)^k}{k!\,(k+1)!}\,\zeta^k\,,
	\,$$
	defines its 1 - sum in such a direction.
\ele

	\proof 
	
	Because of convergence of the sequence $\{S_k\}$ of the partial sums
	of the series \eqref{S} there exists a positive constant $B$ such that
	  $$\,
		    |S_k| \leq B \quad
				\textrm{for all}\quad k \geq 0.
		\,$$
		Then for the coefficients $b_k=\frac{(-1)^k k!\,S_{k-1}}{(\beta_2-\beta_1)^k}$
		of the formal series $\hat{\psi}(x)$ we find the following estimate
		 $$\,
		  |b_k| \leq B\,A^k\,k!\quad
			\textrm{for all }\quad k\in\NN\,,
		\,$$
 where $A=1/|\beta_2-\beta_1|$. Therefore the series $\hat{\psi}(x)$ is of
Gevrey order 1.

As a result the formal Borel transform of the series $\hat{\psi}(x)$
 $$\,
  \psi(\zeta)=\hat{\B}_1(\hat{\psi})(\zeta)=
	\sum_{k=1}^{\infty} \frac{(-1)^k\,S_{k-1}}{(\beta_2-\beta_1)^k}\,\zeta^k
\,$$ 
is an analytic function near the origin in the Borel $\zeta$-plane. Moreover, it turns out
that
 $$\,
  \psi(\zeta)=-\left(\frac{1}{1 + \frac{\zeta}{\beta_2-\beta_1}}-1\right)\,
	\frac{J_1(2 \sqrt{(\gamma_1-\gamma_2)\,\zeta})}{\sqrt{(\gamma_1-\gamma_2)\,\zeta}}\,,
\,$$
where 
$$\,
  \frac{J_1(2 \sqrt{(\gamma_1-\gamma_2)\,\zeta})}{\sqrt{(\gamma_1-\gamma_2)\,\zeta}}=
	 \sum_{k=0}^{\infty} \frac{(\gamma_2-\gamma_1)^k}{k!\,(k+1)!}\,\zeta^k
\,$$
is the Bessel function of order 1 (see \cite{W})
 $$\,
  J_n(z)=\sum_{k=0}^{\infty} \frac{(-1)^k (z/2)^{2 k + n}}{\Gamma(k+1)\,\Gamma(n+k+1)}
\,$$
and $\Gamma(z)$ is the classical Gamma function (see \cite{HB-AE}). In particular,
the analytic in $\CC$ function $J_1(2 \sqrt{(\gamma_1-\gamma_2)\,\zeta})/\sqrt{(\gamma_1-\gamma_2)\,\zeta}$
is a particular solution near $\zeta=0$ of the following differential equation
$$\,
 \zeta\,u''(\zeta) + 2 u'(\zeta) - (\gamma_2-\gamma_1)\,u(\zeta)=0\,.
\,$$
For simplicity we denote by $u(\zeta)$ the function 
 $J_1(2 \sqrt{(\gamma_1-\gamma_2)\,\zeta})/\sqrt{(\gamma_1-\gamma_2)\,\zeta}$. Then $\psi(\zeta)$ becomes
  $$\,
   \psi(\zeta)=\frac{\zeta}{\zeta + \beta_2-\beta_1}\,u(\zeta)\,.
  \,$$

 The function $\psi(\zeta)$ has an analytic continuation along any ray $\theta \neq \arg(\beta_1-\beta_2)$ and 
   $$\,
   |\psi(\zeta)| \leq A_0\,e^{|\gamma_2-\gamma_1|\,|\zeta|}
   \,$$
   along such a ray for an appropriate constant $A_0$.
 Then the complex Laplace transform of the function $\psi(\zeta)$ gives
\ben
 \psi_{\theta}(x)   &=&
\L_{\theta} \psi(x)=
\int_0^{+\infty\,e^{i \theta}} \psi(\zeta)\,\exp\left(-\frac{\zeta}{x}\right)\,
d \left(\frac{\zeta}{x}\right)=\\[0.4ex]
                    &=&
		\frac{1}{x}\int_0^{+\infty\,e^{i \theta}}
		 \frac{\zeta\,u(\zeta)\,e^{-\frac{\zeta}{x}}}{\zeta+\beta_2-\beta_1}\,d \zeta=\\[0.2ex]
		     &=&
	\frac{1}{x} \int_0^{+\infty e^{i \theta}} u(\zeta)\,e^{-\frac{\zeta}{x}}\,d \zeta	
	-\frac{\beta_2-\beta_1}{x} \int_0^{+\infty e^{i \theta}} 
	\frac{u(\zeta)\,e^{-\frac{\zeta}{x}}}{\zeta + \beta_2-\beta_1}\,d \zeta\,. 
			\een
   The analytic function $u(\zeta)$ can be regarded as the formal Borel transform of the convergent
   series
    $$\,
     \sum_{k=0}^{\infty} \frac{(\gamma_2-\gamma_1)^k}{(k+1)!}\,x^k=
     \sum_{s=1}^{\infty} \frac{(\gamma_2-\gamma_1)^{s-1}}{s!}\,x^{s-1}=
     \frac{1}{x (\gamma_2-\gamma_1)} \left(e^{(\gamma_2-\gamma_1)\,x}-1\right)\,.
    \,$$
    Since the Borel sum of an analytic function in $\CC$ gives the same function, we find that for any ray
    $\theta \neq \arg (\beta_1-\beta_2)$ the function
      $$\,
      \psi_{\theta}(x)= 
      	-\frac{\beta_2-\beta_1}{x} \int_0^{+\infty e^{i \theta}} \frac{u(\zeta)\,e^{-\frac{\zeta}{x}}}
      	{\zeta + \beta_2-\beta_1}\,d \zeta
      	 + \frac{1}{x (\gamma_2-\gamma_1)} \left(e^{(\gamma_2-\gamma_1)\,x}-1\right)\,.
      \,$$			
defines the 1 - sum of $\hat{\psi}(x)$ in such a  direction.

This completes the proof.
\qed

   \bre{glue-0}
    Since
      $$\,
      \left|\frac{u(\zeta)}{\zeta+\beta_2-\beta_1}\right| \leq A\,e^{|\gamma_2-\gamma_1| |\zeta|}
      \,$$
      for an appropriate constant $A > 0$, the Laplace integral
       $$\,
        \int_0^{+\infty e^{i \theta}}
          \frac{u(\zeta)\,e^{-\frac{\zeta}{x}}}{\zeta +\beta_2-\gamma_1}\,d \zeta
       \,$$
       exists and defines a holomorphic function in the open disc \cite{L-R}
        $$\,
          \mathcal{D}_{\theta}(|\gamma_2-\gamma_1|)=\left\{x\in\CC\,|
      Re \,\left(\frac{e^{i \theta}}{x}\right) > |\gamma_2-\gamma_1|\right\}\,.
        \,$$
        When we move the direction $\theta\in I$ continuously the corresponding 1-sums $\psi_{\theta}(x)$
        stick each other analytically and define a holomorphic function $\tilde{\psi}(x)$ on a sector
         $$\,
          \bigcup_{\theta\in I} \tilde{\mathcal{D}}_{\theta}(|\gamma_2-\gamma_1|)
         \,$$
         whose opening is $> 2 \pi$. Here $I$ is an open interval in $\RR$ such that
         $$\,
         I= \left( -2 \pi + \arg (\beta_1-\beta_2),\,\arg (\beta_1-\beta_2)\right)\,,
         \,$$
         if  $\arg (\beta_1-\beta_2)\in [-\pi, 2 \pi)$,
         and
         $$\,
         I=\left( \arg (\beta_1-\beta_2,\,\arg(\beta_1-\beta_2) + 2 \pi\right)
         \,$$
         if $\arg (\beta_1-\beta_2)\in [0. \pi)$. 
         The set $\tilde{\mathcal{D}}_{\theta}(|\gamma_2-\gamma_1|)$ is the lifting of $\mathcal{D}_{\theta}(|\gamma_2-\gamma_1|)$
         on the Riemann surface of the natural logarithm $\tilde{\CC}$, i.e.
          $$\,
            \tilde{\mathcal{D}}_{\theta}(|\gamma_2-\gamma_1|)=\left\{
            x=r \underline{e^{i t}}\in\tilde{\CC}\,|\, r < \frac{1}{|\gamma_2-\gamma_1|},
             t\in \left(\theta - \arccos r |\gamma_2-\gamma_1|,\, \theta + \arccos r |\gamma_2-\gamma_1|\right)\right\}\,.
          \,$$
          The notation $\underline{e^{i t}}$ represents a point on $\tilde{\CC}$ which is "above" the complex number
          $e^{i t}$. Then the point $r e^{i t}$ is the canonical projection of $r \underline{e^{i t}}$ on $\CC$.
          In this way the points $r \underline{e^{i t}}$ and $r \underline{e}^{i(t +2 \pi\,m)}, m\in\ZZ$ are regarded as
          different points on $\tilde{\CC}$ but with the same projection on $\CC$ (see \cite{Sa} for details).
          
          The function $\tilde{\psi}(x)$ is a multivalued function on this sector. It is asymptotic to the series
          $\hat{\psi}(x)$ in Gevrey order 1 sense on this sector and defines the 1-sum of this series. 
          In every non-singular direction $\theta$ the multivalued function $\tilde{\psi}(x)$ has one value
          $\psi_{\theta}(x)$. Near the singular direction $\theta=\arg (\beta_1-\beta_2)$ the function
          $\tilde{\psi}(x)$ has two different values: $\psi^{+}_{\theta}(x)=\psi_{\theta+\epsilon}(x)$ and
          $\psi^{-}_{\theta}(x)=\psi_{\theta-\epsilon}(x)$ where $\epsilon > 0$ is a small number.
            \ere

   We have a similar result for the summation of the divergent series $\hat{\varphi}(x)$ from
   \eqref{ffs-inf}.

\ble{sum2}
  Let $\beta_2 \neq \beta_1,\,\gamma_2 \neq \gamma_1$. Assume also that $S\neq 0$, where $S$ is given by \eqref{S}.
  Then the formal power series $\hat{\varphi}(x)$
	defined by \eqref{ffs-inf} is 1 - summable in any direction 
	$\theta \neq \arg(\gamma_2 - \gamma_1)$.  
      The function
	 \be\label{s3}
	     \varphi_{\theta}(x) =
	  	- x\,
              \int_0^{+\infty\,e^{i \theta}} \frac{v(p)\,e^{-x\,p}}
                  {1-\frac{p}{\gamma_2-\gamma_1}}\,d  p -
              \frac{x\,(e^{-\frac{\beta_2-\beta_1}{x}}-1)}{\beta_2-\beta_1},									
	\ee
	where $v(p)$  is  the analytic in $\CC$ function
	 $$\,
	  v(p)=\sum_{k=0}^{\infty} \frac{(-1)^k\,(\beta_2-\beta_1)^k}{k!\,(k+1)!}\, p^k, 
	\,$$
	defines its 1 - sum in such a direction.
\ele 
  
     \proof
    
    Similar to the series $\hat{\psi}(x)$ the coefficients of the series $\hat{\varphi}(x)$ satisfy the 
    following estimate
     $$\,
        \left|\frac{k!\,S_{k-1}}{(\gamma_2-\gamma_1)^k}\right| \leq B\,A^k\,k!\quad
        \textrm{for all} \quad k\in\NN\,,      
     \,$$
     where $A=1/|\gamma_2-\gamma_1|$. The constant $B$ is the same as in the proof of \leref{sum1}.
     Therefore the series $\hat{\varphi}(x)$ is of Gevrey order 1.
    As a result the formal Borel transform of the series $\hat{\varphi}(x)$
     $$\,
         \varphi(p)=\hat{\B}_1\,\hat{\varphi}(p)=
               \sum_{k=0}^{\infty} \frac{(k+1)!\,S_k}{(\gamma_2-\gamma_1)^{k+1}\,k!}\,p^k=
              \sum_{k=0}^{\infty} \frac{(k+1)\,S_k}{(\gamma_2-\gamma_1)^k}\,p^k
     \,$$
              is an analytic function near the origin in the Borel $p$-plane. Observe that
        $$\,
            \varphi(p)=\frac{d}{d p} \left(\sum_{k=0}^{\infty} \frac{S_k\,p^{k+1}}{(\gamma_2-\gamma_1)^{k+1}}\right)=
                \frac{d W(p)}{d p}.
        \,$$
           The series $W(p)$ defines the function
            $$\,
                       W(p)=\left(\frac{1}{1-\frac{p}{\gamma_2-\gamma_1}} -1\right)
                      \sum_{k=0}^{\infty} \frac{(-1)^{k+1}\,(\beta_2-\beta_1)^k}{k!\,(k+1)!}\,p^k.
          \,$$
       As in the proof of \leref{sum1} the convergent in $\CC$ series defines the Bessel function of order 1
               $$\,
                           \frac{J_1(2 \sqrt{(\beta_2-\beta_1) p})}{\sqrt{(\beta_2-\beta_1) p}}=
                           \sum_{k=0}^{\infty} \frac{(-1)^k \,(\beta_2-\beta_1)^k}{k!\,(k+1)!}\,p^k.
           \,$$
         For simplicity we denote by $v(p)$ the analytic in $\CC$ function $J_1(2 \sqrt{(\beta_2-\beta_1) p})/\sqrt{(\beta_2-\beta_1) p}$. In particular, the entire function $v(p)$ satisfies the differential
           equation
         $$\,
           p\,v''(p) + 2 v'(p) + (\beta_2-\beta_1)\,v(p)=0\,.
         \,$$
          Then $W(p)$ becomes
                 $$\,
                    W(p)=\left(1 - \frac{1}{1-\frac{p}{\gamma_2-\gamma_1}} \right)\,v(p).
          \,$$
         Now applying the properties of the complex Laplace transform (see \leref{prop}) we find that
          \ben
           (\L_{\theta} \varphi)(x) &=& 
           \left(\L_{\theta} \,\frac{d W}{d p}\right) (x)=
           x\,(\L_{\theta} W)(x) - W(0)=\\[0.15ex]
                     &=&
                 x \int_0^{+\infty e^{i \theta}}
                 \left(1- \frac{1}{1- \frac{p}{\gamma_2-\gamma_1}}\right)\,v(p)\,e^{-x p}\, d p=\\[0.2ex]
                 &=&
             x \int_0^{+\infty e^{i \theta}} v(p)\,e^{-x p}\, d p
                      	- x\,
                      	\int_0^{+\infty\,e^{i \theta}} \frac{v(p)\,e^{-x\,p}}
                      	{1-\frac{p}{\gamma_2-\gamma_1}}\,d  p \,.
          \een
         Since $v(p)/(1-\frac{p}{\gamma_2-\gamma_1})$ has exponential growth $\leq 1$ at $\infty$, i.e.
          $$\,
         \left|\frac{v(p)}{1- \frac{p}{\gamma_2-\gamma_1}}\right| \leq A_0\,e^{|\beta_2-\beta_1| |p|}\,,
          \,$$
          along any direction $\theta \neq \arg (\gamma_2-\gamma_1)$ form 0 to $+\infty e^{i \theta}$ the integral
           $$\,
             \int_0^{+\infty e^{i \theta}}
             \frac{v(p)\,e^{-x p}}{1-\frac{p}{\gamma_2-\gamma_1}}\, d p
           \,$$
           defines a holomorphic function on (see \cite{Sa})
            $$\,
             \pi_{\theta}=\left\{x\in\CC\,|\, Re (x.e^{i \theta}) > |\beta_2-\beta_1|\right\}\,.
            \,$$
         The analytic in $\CC$ function 
          $$\,
           v(p)= \sum_{k=0}^{\infty} \frac{(-1)^k \,(\beta_2-\beta_1)^k}{k!\,(k+1)!}\,p^k
          \,$$
         can be regarded as the formal Borel transform of the convergent series
          \ben
             \sum_{k=0}^{\infty} \frac{(-1)^k \,(\beta_2-\beta_1)^k}{(k+1)!}\,x^{-k-1}
                 &=&
             -\frac{1}{\beta_2-\beta_1} \sum_{s=1}^{\infty}
             \frac{(-1)^s\,(\beta_2-\beta_1)^s}{s!}\,x^{-s}=\\[0.2ex]
                &=&
                -\frac{1}{\beta_2-\beta_1}\left[e^{-\frac{\beta_2-\beta_1}{x}} -1\right]\,.
          \een
         Since the Borel sum of an analytic function in $\CC\PP^1 - \{0\}$ gives the same function,  
        for any ray $\theta\neq \arg(\gamma_2-\gamma_1)$ the function
           $$\,
                \varphi_{\theta}(x) =
	  	- x\,
              \int_0^{+\infty\,e^{i \theta}} \frac{v(p)\,e^{-x\,p}}
                  {1-\frac{p}{\gamma_2-\gamma_1}}\,d  p -
              \frac{x\,(e^{-\frac{\beta_2-\beta_1}{x}}-1)}{\beta_2-\beta_1}								
           \,$$
     defines the 1-sum of the formal series $\hat{\varphi}(x)$.

    This ends the proof.
\qed
  
       \bre{glu-inft}
           Following \cite{Sa} we denote by $\tilde{\pi}_{\theta}$ the lifting of $\pi_{\theta}$ of
           the Riemann surface of the natural logarithm  $\tilde{\CC}$, i.e.
           $$\,
             \tilde{\pi}_{\theta} : =\left\{\,x=r \underline{e^{i t}}\in \tilde{\CC}\,|\, r > |\beta_2-\beta_1|,
             t\in \left(-\theta - \arccos \frac{|\beta_2-\beta_1|}{r}, 
             -\theta + \arccos \frac{|\beta_2-\beta_1|}{r}\right)\right\}\,,
           \,$$
           Denote also by 
           $\tilde{D}(I, |\beta_2-\beta_1|)$ the sector on $\tilde{\CC}$
             $$\,
               \tilde{D}(I, |\beta_2-\beta_1|) := \bigcup_{\theta\in I} \tilde{\pi}_{\theta}\,.
             \,$$
             Here $I$ is an open interval in $\RR$ such that
             $$\,
               I= \left( -2 \pi + \arg (\gamma_2-\gamma_1),\,\arg (\gamma_2-\gamma_1)\right)\,,
             \,$$
             if  $\arg (\gamma_2-\gamma_1)\in [-\pi, 2 \pi)$,
             and
             $$\,
               I=\left( \arg (\gamma_2-\gamma_1,\,\arg(\gamma_2-\gamma_1) + 2 \pi\right)
             \,$$
             if $\arg (\gamma_2-\gamma_1)\in [0. \pi)$. When we move the direction $\theta\in I$
             continuously the corresponding 1-sums $\varphi_{\theta}(x)$ stick each other analytically and define
             a holomorphic function $\tilde{\varphi}(x)$ on the sector $\tilde{D}(I, |\beta_2-\beta_1|)$
             of opening $> 2 \pi$. On this sector the multivalued function $\tilde{\varphi}(x)$ defines the
             Borel sum of the series $\hat{\varphi}(x)$ and it is asymptotic to this series in Gevrey order 1. 
             In every non-singular direction $\theta$ it has one
             value $\varphi_{\theta}(x)$. Near the singular direction $\theta=\arg (\gamma_2-\gamma_1)$
             the multivalued function $\tilde{\varphi}(x)$ has two differnt values : 
             $\varphi^{+}_{\theta}=\varphi_{\theta+\epsilon}(x)$ and 
             $\varphi^{-}_{\theta}(x)=\varphi_{\theta-\epsilon}(x)$ for a small number $\epsilon > 0$.
       \ere

   Denote by $F(x)$ the actual matrix-function $x^{\Lambda}\,\exp(-B/x)$, where the matrices $\Lambda$ and $B$
   are defined by \prref{p2}.
  The first existence result says 
 
\bpr{p-A}
 Let $\beta_1 \neq \beta_2$. 
 \begin{enumerate}
 	\item\,
 	Assume that $S=0$ where $S$ is defined by \eqref{S} but $\gamma_1 \neq \gamma_2$. 
 	Then the initial equation possesses an unique actual
 	fundamental matrix $\Phi(x, 0)$ at the origin in the form
 	\be\label{ac}
 	   \Phi(x, 0)=\exp (G x)\,H(x)\,F(x)\,,
 	 \ee
 	 where $H(x)=\hat{H}(x)$, defined by \eqref{H}, is an analytic in $\CC$ matrix-function. The matrix $F(x)$
 	 is the branch of $x^{\Lambda}\,\exp(-B/x)$ for $\arg (x)$ and the matrix $G$ is introduced by \prref{p2}.
 	 
 	 \item\,
 	 Assume that $\gamma_1=\gamma_2$. Then the initial equation possesses an unique actual fundamental matrix
 	 at the origin in the form \eqref{ac}, where the matrix $H(x)=\hat{H}(x)$ is defined by \eqref{H-g}.
 	 The matrix $F(x)$
 	 is the branch of $x^{\Lambda}\,\exp(-B/x)$ for $\arg (x)$ and the matrix $G$ is introduced by \prref{p2}.
 	 
 	 \item\,
 	 Assume that $S \neq 0$ where $S$ is defined by \eqref{S} but $\gamma_1 \neq \gamma_2$.
 Then for every non-singular direction $\theta$ the  initial equation
possesses an unique actual fundamental matrix $\Phi^{\theta}_0(x, 0)$ at the origin of the form
 \be\label{f}
   \Phi^{\theta}_0(x, 0)=\exp(G x)\,H_{\theta}(x)\,F_{\theta}(x)\,,
   \quad
   \Phi^{\theta+2 \pi}_0(x, 0)=\Phi^{\theta}_0(x, 0)\,\hat{M}=\Phi^{\theta}_0(x, 0)\,.
\ee
	The matrix $F_{\theta}(x)$ is  the branch 
	of  $x^{\Lambda}\,\exp(-B/x)$ for $\arg(x)=\theta$ and the matrix $G$ is introduced by \prref{p2}. 
	The matrix $H_{\theta}(x)$ is defined by 
	$$\,
	   H_{\theta}(x)=\left(\begin{array}{cc}
		     1    &\frac{x^2}{\beta_2-\beta_1} 
		      +x^2\,\phi_{\theta}(x)\\[0.2ex]
				  0   &1
					      \end{array}
					\right)\,,			
	\,$$
	where
	$$\,
	  \phi_{\theta}(x)=\int_0^{+\infty\,e^{i \theta}}
		 \frac{u(\zeta)\,e^{-\frac{\zeta}{x}}}{\zeta+\beta_2-\beta_1}\,d \zeta
	\,$$
	with
	$$\,
	  u(\zeta)=\sum_{k=0}^{\infty} \frac{(\gamma_2-\gamma_1)^k}{k!\,(k+1)!}\,\zeta^k\,.
	\,$$
	
	For the singular direction $\theta=\arg (\beta_1-\beta_2)$ the initial equation admits two actual fundamental
	matrices at the origin
	 $$\,
	   (\Phi^{\theta}_0(x, 0))^{\pm}=\Phi^{\theta \pm \epsilon}_0(x, 0)\,,
	 \,$$
	 where the matrices $\Phi^{\theta \pm \epsilon}_0(x, 0)$ are given by \eqref{f}
	  and  $\epsilon > 0$ is a small number.
	  \end{enumerate}
\epr

\bth{A-S}
  With respect to the actual fundamental matrix at the origin given by \prref{p-A} and extended by \reref{glue-0},
	the initial equation  has one singular direction $\theta=\arg(\beta_1-\beta_2)$.
	The corresponding Stokes matrix $St^{\theta}_0$ is
	  $$\,
		  St^{\theta}_0=\left(\begin{array}{cc}
			     1       &-2 \pi\,i(\gamma_2-\gamma_1)\,S\\
					 0       &1
					              \end{array}
									\right)\,,			
		\,$$
		where $S$ is introduced by \eqref{S}.
	
\ethe

\proof

Let $\theta^+=\theta+\epsilon$ and $\theta^-=\theta-\epsilon$ for a small positive $\epsilon$
be two directions  which are slightly to the left
and to the right of the singular direction $\theta=\arg(\beta_1-\beta_2)$. Let
$\Phi^{\theta+}_0(x, 0)$ and $\Phi^{\theta-}_0(x, 0)$ be the associated actual
fundamental matrices at $x=0$ on these directions (in sense of \prref{p-A}).

Then comparing $\Phi^{\theta+}_0(x, 0)$ and $\Phi^{\theta -}_0(x, 0)$ we find
 $$\,
   \Phi^{\theta -}_0(x, 0) - \Phi^{\theta +}_0(x, 0)=
	   \left(\begin{array}{cc}
		       1     &\mu^{\theta}_0\\
					 0     &1
					 \end{array}
		\right)\,,			
\,$$
where
\ben
  \mu^{\theta}_0    &=&
     (\gamma_2-\gamma_1)\,e^{\gamma_1\,x}\,e^{-\frac{\beta_2}{x}}
	  [\phi_{\theta}^-(x) - \phi_{\theta}^+(x)]=\\[0.2ex]
		       & =&
		2 \pi\,i\,(\gamma_2-\gamma_1)\,e^{\gamma_1\,x}\,e^{-\frac{\beta_2}{x}}\,
		\res\left(\frac{u(\zeta)\,e^{-\frac{\zeta}{x}}}{\zeta+\beta_2-\beta_1}\,;\,\zeta=\beta_1-\beta_2\right)=
		    \\[0.35ex]
				  &=&
	  2 \pi\,i\,(\gamma_2-\gamma_1)\,u(\beta_1-\beta_2)\,\Phi_1(x, 0)\,.				
\een
Observe that $u(\beta_1-\beta_2)=-S$ where $S$ is the sum of the absolutely convergent number series
  $$\,
    \sum_{k=0}^{\infty} \frac{(-1)^{k+1}\,(\gamma_2-\gamma_1)^k\,(\beta_2-\beta_1)^k}{k!\,(k+1)!}\,,
  \,$$
  introduced by \eqref{S}.
Now, a direct application of the definition of the Stokes matrix gives us the result. 

Note that when either $S=0$ or $\gamma_1=\gamma_2$ the Stokes matrix $St^{\theta}_0=I_2$. This phenomenon witness the fact that the
series $\hat{\psi}(x)$ from \eqref{psi} is convergent when $S=0$ or $\gamma_1=\gamma_2$. Moreover it is in concordance with a theorem
on the analytic dependence of the Stokes matrix $St^{\theta}_0$ on the parameter  (see \cite{HLR}).

This ends the proof. 
\qed

The second existence result states
\bpr{p-A-inft}
 Let that $\beta_1 \neq \beta_2$. 
   \begin{enumerate}
   	
   	\item\,
   	Assume that $\gamma_1=\gamma_2$. Then the initial equation possesses an unique
   	actual fundamental matrix $\Phi_{\infty}(x, 0)$ at $x=\infty$ in the form
   	 \be\label{act-inft}
   	\Phi_{\infty}(x, 0)=\exp \left(-\frac{B}{x}\right)\,\left(\frac{1}{x}\right)^{-\Lambda}\,
   	P(x)\,\exp(G x)\,,
   	 \ee
   	where $P(x)=\hat{P}(x)$ is defined by \eqref{P-g}. The matrices $\Lambda, B$ and $G$ are defined by \prref{p2} and $\exp(G x)$ is the branch of $\exp(G x)$
   	for $\arg (x)$.
   
   \item\,
   Assume that $\gamma_1 \neq \gamma_2$ but $S=0$ where $S$ is introduced by \eqref{S}. Then the initial equation possesses an unique
   actual fundamental matrix $\Phi_{\infty}(x, 0)$ at $x=\infty$ in the form \eqref{act-inft},
  where $P(x)=\hat{P}(x)$, defined by \eqref{P}, is an analytic in $\CC\PP^1 - \{0\}$ matrix-function.
  The matrices $\Lambda, B$ and $G$ are defined by \prref{p2} and $\exp(G x)$ is the branch of $\exp(G x)$
  for $\arg (x)$.
  
  \item\,
    Assume that $\gamma_1 \neq \gamma_2$ but $S \neq 0$ where $S$ is introduced by \eqref{S}.
 Then for every non-singular direction
 $\theta$ the initial equation
possesses an unique actual fundamental matrix $\Phi^{\theta}_{\infty}(x, 0)$ near $x=\infty$ of the form
 \be\label{d-inf}
   \Phi^{\theta}_{\infty}(x, 0)=\exp\left(-\frac{B}{x}\right)\,\left(\frac{1}{x}\right)^{-\Lambda}\,P_{\theta}(x)\,(\exp(G x))_{\theta}\,.
\ee
  The matrices $B, \Lambda$ and  $G$ are given in \prref{p2} and $(\exp(G x))_{\theta}$ is the branch
  of $\exp (G x)$ for $\arg (x)=\theta$.
	The matrix $P_{\theta}(x)$ is defined by 
	$$\,
	   P_{\theta}(x)=\left(\begin{array}{cc}
		     1    &\frac{e^{-\frac{\beta_2-\beta_1}{x}} -1}{\beta_2-\beta_1} + \omega_{\theta}(x)\\[0.2ex]
		     0    &1
					      \end{array}
					\right)\,,			
	\,$$
	where
	$$\,
	  \omega_{\theta}(x)=\int_0^{+\infty\,e^{i \theta}}
		 \frac{v(p)\,e^{-x p}}{1- \frac{p}{\gamma_2-\gamma_1}}\,d  p
	\,$$
	with
	$$\,
	  v(p)=\sum_{k=0}^{\infty} \frac{(-1)^k\,(\beta_2-\beta_1)^k}{k!\,(k+1)!}\,p^k\,.
	\,$$
	
	For the singular direction $\theta=\arg (\gamma_2-\gamma_1)$  the initial equation admits two actual fundamental
	matrices at $x=\infty$
	 $$\,
	   (\Phi^{\theta}_{\infty}(x, 0))^{\pm}= \Phi^{\theta \pm \epsilon}_{\infty}(x, 0)\,,
	 \,$$
	 where the matrices $\Phi^{\theta \pm \epsilon}_{\infty}(x, 0)$ are given by \eqref{d-inf} and $\epsilon > 0$ is a 
	 small number.
	 \end{enumerate}
\epr

\bth{A-S-inft}
  With respect to the actual fundamental matrix near $x=\infty$ given by \prref{p-A-inft},
 and extended by \reref{glu-inft} the initial equation has one singular direction $\theta=\arg(\gamma_2-\gamma_1)$.
	The corresponding Stokes matrix $St^{\theta}_{\infty}$ is
	  $$\,
		  St^{\theta}_{\infty}=\left(\begin{array}{cc}
			     1       &2 \pi\,i(\gamma_2-\gamma_1)\,S\\
					 0       &1
					              \end{array}
									\right)\,,			
		\,$$
		where $S$ is introduced by \eqref{S}.
	
\ethe

   \proof
   The proof is similar to the proof of \thref{A-S}.
   
   We again observe that, according to expectation, $St^{\theta}_{\infty}=I_2$ when either $S=0$ or $\gamma_1=\gamma_2$. 
   This fact one more time confirms
   the convergence of the series $\hat{\varphi}(x)$ from \eqref{ffs-inf} and is in keeping with a theorem on analytic dependence
   of the Stokes matrices on the parameter (see \cite{HLR}). 
   \qed

  %%%%%%%%%%%%%%%%%%%%%%%%%%%%%%%%%%%%%%%%%%%%%%%%%%%%
	% HE
	%%%%%%%%%%%%%%%%%%%%%%%%%%%%%%%%%%%%%%%

    \section{ Heun type equation }
		
     In this section we will compute the monodromy matrices of the Heun type equation.

	The Heun type equation is invariant under the transformation
               $$\,
                    \sqrt{\varepsilon} \longrightarrow -\sqrt{\varepsilon}.
                \,$$
          So through this section we assume that
              \be\label{rest}
                  \arg (\sqrt{\varepsilon}) \in \left(- \frac{\pi}{2}, \frac{\pi}{2}\right].
              \ee
             
     Observe that in this case 
                    $$\,
                  \arg (\frac{1}{\sqrt{\varepsilon}}) \in \left[- \frac{\pi}{2}, \frac{\pi}{2}\right).
              \,$$
             Observe also that 
               $\arg (\sqrt{\varepsilon}) = \arg (\frac{1}{\sqrt{\varepsilon}})$ if and only if
               $\sqrt{\varepsilon}\in\RR^{+}$, and  $\arg (\sqrt{\varepsilon}) = \arg (-\frac{1}{\sqrt{\varepsilon}})$ if 
               and only if $\sqrt{\varepsilon}\in i\,\RR^{+}$.

   %%%%%%%%%%%%%%%%%%%%%%%%%%%%%%%%%%%%%%%%%%%%%%%%%%%%%%%
	 % A
   %%%%%%%%%%%%%%%%%%%%%%%%%%%%%%%%%%%%%%%%%%%%%%%%%%%%%%%%%	

    In the introduction we have  denoted by $x_j, j=R, L$ the singular points $\sqrt{\varepsilon}$ and $-\sqrt{\varepsilon}$, and
     by $x_{jj}, j=R, L$ the singular points $1/\sqrt{\varepsilon}$ and $-1/\sqrt{\varepsilon}$.
     Recall that when $\beta_1 \neq \beta_2,\,\gamma_1 \neq \gamma_2$ the initial equation admits two Stokes matrices -- 
     one of them correspond to the origin and the other corresponds to the infinity point. In order
     to realize both Stokes matrices as a limit of two different matrices $e^{2 \pi\,i T_j}$ we
     concern with this Heun type equation which has exactly two resonant singular points. One of them will be
     of the type $x_j$, and the other will be of the type $x_{jj}$.
    Throughout this paper we call such a resonance for which two singular points of different type
    are resonant, {\bf a double resonance}. We  distinguish 4 different types of double resonances.
      \begin{itemize}
      	\item\,
      	The double resonance of type ${\bf (A.1)}$. It is defined by the condition
      	      \ben {\bf A.1}\qquad
      	      \Delta^L_{12}=\frac{\beta_2-\beta_1}{2 \sqrt{\varepsilon}},\quad
      	      \Delta^{LL}_{12}+1=\frac{\gamma_1-\gamma_2}{2 \sqrt{\varepsilon}} \in\NN_0.
      	      \een
      	 \item\,
      	 	The double resonance of type ${\bf (A.2)}$. It is defined by the condition
      	 	   \ben {\bf A.2}\qquad
      	 	   \Delta^L_{12}=\frac{\beta_2-\beta_1}{2 \sqrt{\varepsilon}},\quad
      	 	   \Delta^{RR}_{12}+1=\frac{\gamma_2-\gamma_1}{2 \sqrt{\varepsilon}} \in\NN_0.
      	 	   \een     
      	 \item\,
      	     	The double resonance of type ${\bf (A.3)}$. It is defined by the condition	   
                       \ben {\bf A.3}\qquad
                       \Delta^R_{12}=\frac{\beta_1-\beta_2}{2 \sqrt{\varepsilon}},\quad
                       \Delta^{LL}_{12}+1=\frac{\gamma_1-\gamma_2}{2 \sqrt{\varepsilon}} \in\NN_0.
                       \een
         \item\,
            	The double resonance of type ${\bf (A.4)}$. It is defined by the condition	                 
                        \ben {\bf A.4}\qquad
                        \Delta^R_{12}=\frac{\beta_1-\beta_2}{2 \sqrt{\varepsilon}},\quad
                        \Delta^{RR}_{12}+1=\frac{\gamma_2-\gamma_1}{2 \sqrt{\varepsilon}} \in\NN_0.
                        \een
      \end{itemize} 
      
         In accordance with the initial equation we consider two different fundamental matrix 
         $\Phi_j(x, \varepsilon),\,j=0, \infty$ of the perturbed equation. The fundamental matrix
         $\Phi_0(x, \varepsilon)$ corresponds to the fundamental matrix $\Phi_0(x, 0)$ of the initial equation at the origin.
         Therefore the integral that defines the function $\Phi_{12}(x, \varepsilon)$ must be
         taken in the direction $\arg (\beta_2-\beta_1)$ with the same base point $x$ as in the definition
         of the matrix $\Phi_0(x, 0)$.
          The fundamental matrix
          $\Phi_{\infty}(x, \varepsilon)$ corresponds to the fundamental matrix $\Phi_{\infty}(x, 0)$ of the initial equation at $x=\infty$.
          Therefore its element $\Phi_{12}(x, \varepsilon)$ contains an integral that must be taken
          in the direction $\arg (\gamma_1-\gamma_2)$ with the same base point $x$ as in the definition of the matrix
          $\Phi_{\infty}(x, 0)$. The next theorem provides the explicit form of
          the fundamental matrices $\Phi_0(x, \varepsilon)$ and $\Phi_{\infty}(x, \varepsilon)$ introduced by
          \thref{global}.
          
          \bth{expl-A}
           Assume that $\beta_1 \neq \beta_2$. Then both fundamental matrices
           $\Phi_0(x, \varepsilon)$ and $\Phi_{\infty}(x, \varepsilon)$   
           have the same elements $\Phi_j(x, \varepsilon), j=1, 2$
    
		 \ben 
		  \Phi_1(x, \varepsilon)=
			\left(\frac{x-\sqrt{\varepsilon}}{x+\sqrt{\varepsilon}}\right)^{\frac{\beta_1}{2 \sqrt{\varepsilon}}}\,
			\left(\frac{\frac{1}{\sqrt{\varepsilon}}+x}{\frac{1}{\sqrt{\varepsilon}}-x}\right)
			  ^{\frac{\gamma_1}{2 \sqrt{\varepsilon}}}\,,
				 \Phi_2(x, \varepsilon)=
			\frac{(x-\sqrt{\varepsilon})^{\frac{\beta_2}{2 \sqrt{\varepsilon}}-1}}
			 {(x+\sqrt{\varepsilon})^{\frac{\beta_2}{2 \sqrt{\varepsilon}}+1}}\,
			\left(\frac{\frac{1}{\sqrt{\varepsilon}}+x}{\frac{1}{\sqrt{\varepsilon}}-x}\right)
			  ^{\frac{\gamma_2}{2 \sqrt{\varepsilon}}}.
		\een
           The element $\Phi_{12}(x, \varepsilon)$
            of the fundamental matrix $\Phi_0(x, \varepsilon)$ is defined as	
		  \ben	  
				\Phi_{12}(x, \varepsilon)=\Phi_1(x, \varepsilon)\,
					\int_{\epsilon}^x 
					\frac{(z-\sqrt{\varepsilon})^{\frac{\beta_2-\beta_1}{2 \sqrt{\varepsilon}}-1}}
					{(z+\sqrt{\varepsilon})^{\frac{\beta_2-\beta_1}{2 \sqrt{\varepsilon}}+1}}\,
					\left(\frac{\frac{1}{\sqrt{\varepsilon}}+z}{\frac{1}{\sqrt{\varepsilon}}-z}\right)
					^{\frac{\gamma_2-\gamma_1}{2 \sqrt{\varepsilon}}}\,d z\,,
		 \een 
		 where $\epsilon=\sqrt{\varepsilon}$ during a double resonance of type {\bf A.1} and {\bf A.2}.
		 During a resonance of type {\bf A.3} or {\bf A.4}  $\epsilon=-\sqrt{\varepsilon}$. The integral 
         is taken in the direction $\arg (\beta_2-\beta_1)$ as the base point $x$ is a point near the origin.

         The element $\Phi_{12}(x, \varepsilon)$ of the fundamental matrix $\Phi_{\infty}(x, \varepsilon)$
           has the form
                  \ben	  
                  \Phi_{12}(x, \varepsilon)=\Phi_1(x, \varepsilon)\,
                  \int_{\epsilon}^x 
                  \frac{(z-\sqrt{\varepsilon})^{\frac{\beta_2-\beta_1}{2 \sqrt{\varepsilon}}-1}}
                  {(z+\sqrt{\varepsilon})^{\frac{\beta_2-\beta_1}{2 \sqrt{\varepsilon}}+1}}\,
                  \left(\frac{\frac{1}{\sqrt{\varepsilon}}+z}{\frac{1}{\sqrt{\varepsilon}}-z}\right)
                  ^{\frac{\gamma_2-\gamma_1}{2 \sqrt{\varepsilon}}}\,d z\,,
                  \een 
                  where $\epsilon=1/\sqrt{\varepsilon}$ during a resonance of type {\bf A.1} and {\bf A.3},
                  and $\varepsilon=-1/\sqrt{\varepsilon}$ during a double resonance of type {\bf A.2} 
                  and {\bf A.4}. The integral is taken in the direction $\arg (\gamma_1-\gamma_2)$
                  as the base point $x$ is a point near $x=\infty$.                
         \ethe 
           
           It turns out that the perturbed equation with the so determined fundamental matrices $\Phi_0(x, \varepsilon)$
           and $\Phi_{\infty}(x, \varepsilon)$ has a double resonance  for very special values of the 
           parameter of perturbation $\varepsilon$.

           \bpr{two}
           Assume that $\beta_1 \neq \beta_2$. Then 
            the Heun type equation with the fundamental matrices $\Phi_0(x, \varepsilon)$ and $\Phi_{\infty}(x, \varepsilon)$,
            introduced by \thref{expl-A}, has a double resonance if and only if $\varepsilon\in\RR^{+}$.
           \epr

           \proof

      The statement follows immediately either from the conditions
        $$\,
          \Delta^{RR}_{21}-1=\frac{\gamma_1-\gamma_2}{2 \sqrt{\varepsilon}}\in\NN,\quad
          \arg (\frac{1}{\sqrt{\varepsilon}})=\arg (\gamma_1-\gamma_2)
        \,$$           
        that define the double resonance of type {\bf A.1} or {\bf A.3},
        or from the conditions
           $$\,
           \Delta^{LL}_{21}-1=\frac{\gamma_2-\gamma_1}{2 \sqrt{\varepsilon}}\in\NN,\quad
           \arg (-\frac{1}{\sqrt{\varepsilon}})=\arg (\gamma_1-\gamma_2)
           \,$$           
           that define the double resonance of type {\bf A.2} or {\bf A.4}.
           Indeed these conditions imply that
             $$\,
                \arg (\frac{1}{\sqrt{\varepsilon}})=\arg (\sqrt{\varepsilon}).
             \,$$
             But the latter is true if and only if $\sqrt{\varepsilon}\in\RR^{+}$.
             And in particular, $\varepsilon\in\RR^{+}$.
           
            In very particular case when $\gamma_1=\gamma_2$ the elements $\Phi_{12}(x, \varepsilon)$ of both
            matrices $\Phi_0(x, \varepsilon)$ and $\Phi_{\infty}(x, \varepsilon)$ do not contain logarithmic terms.
            Moreover, one can assume that the path of integration of the element $\Phi_{12}(x, \varepsilon)$ of
            $\Phi_{\infty}(x, \varepsilon)$ can be taken in any direction, that does not cross the points 
            $x_j, j=R, L$. In this paper on purpose completeness of the text we fix always this path in
            the direction $\arg (\gamma_1-\gamma_2)$. So when $\gamma_1=\gamma_2$ the direction will be the real
            positive axis.
            
           This ends the proof.   
           \qed

         \bre{res-per}
           When $\beta_1=\beta_2$ we meet with the same problem as with the initial equation (see \reref{res-in}). 
           The function 
           $\Phi_{12}(x, \varepsilon)$ as an element of the matrix $\Phi_0(x, \varepsilon)$ becomes
            $$\,
              \Phi_{12}(x, \varepsilon)=
              	\Phi_1(x, \varepsilon)\,
              	\int_{\epsilon}^x 
              	\frac{1}
              	{(z- \sqrt{\varepsilon})\,(z+\sqrt{\varepsilon})}\,
              	\left(\frac{\frac{1}{\sqrt{\varepsilon}}+z}{\frac{1}{\sqrt{\varepsilon}}-z}\right)
              	^{\frac{\gamma_2-\gamma_1}{2 \sqrt{\varepsilon}}}\,d z
            \,$$
            when $\beta_1=\beta_2$. Here $\epsilon$ is either $\sqrt{\varepsilon}$ or $-\sqrt{\varepsilon}$. 
            But this integral does not exist. Recall that when $\beta_1=\beta_2$ the point $x=0$ is a resonant
            irregular singular point for the initial equation (see \reref{res-in}). To overcome the above problem we have to use some other
            kind of solution or some other kind of perturbation. For this reason the so called resonant irregular
            singularities are out of scope of this paper.
         \ere                 

               Thanks to \prref{two} we fix the paths of integration including in the solution
               $\Phi_{12}(x, \varepsilon)$ as follows:
                  \begin{itemize}
                  	\item\,
                  	During a double resonance of type {\bf A.1} for both fundamental matrices
                  	$\Phi_0(x, \varepsilon)$ and $\Phi_{\infty}(x, \varepsilon)$ the path is taken
                  	in the real positive axis. In this case $d_R=d_{RR}=0$.
                    In particular $\arg (\beta_2-\beta_1)=\arg (\gamma_1-\gamma_2)=0$.

                      	\item\,
                      	During a double resonance of type {\bf A.2} for the fundamental matrix
                      	$\Phi_0(x, \varepsilon)$ the path of integration is taken in the real positive
                      	axis. For the matrix $\Phi_{\infty}(x, \varepsilon)$ the path is taken
                      	in the real negative axis. Here $d_R=d_{LL}=0$.
                      	  In particular $\arg (\beta_2-\beta_1)=\arg (\gamma_2-\gamma_1)=0$.

                        	\item\,
                        	During a double resonance of type {\bf A.3} for the  fundamental matrix
                        	$\Phi_0(x, \varepsilon)$ the path of integration is taken in the real
                        	negative axis. For the matrix $\Phi_{\infty}(x, \varepsilon)$ the path is taken
                        	in the real positive axis. Here $d_L=d_{RR}=0$.
                        	  In particular $\arg (\beta_1-\beta_2)=\arg (\gamma_1-\gamma_2)=0$.

                           	\item\,
                           	During a double resonance of type {\bf A.4} for both fundamental matrices
                           	$\Phi_0(x, \varepsilon)$ and $\Phi_{\infty}(x, \varepsilon)$ the path is taken
                           	in the real negative axis. Here $d_L=d_{LL}=0$.
                           	  In particular $\arg (\beta_1-\beta_2)=\arg (\gamma_2-\gamma_1)=0$.        
                          \end{itemize}

                The following  statement describes the behavior of the fundamental matrices $\Phi_0(x, \varepsilon)$
                and $\Phi_{\infty}(x, \varepsilon)$ near the singular points.
                \bth{F}
                Assume that $\beta_1 \neq \beta_2$.
                Then during a double resonance 
                the fundamental matrix $\Phi_0(x, \varepsilon)$ of the Heun type  equation
                is represented near the singular points $x_j,\,j=R, L$ as follows,
                \ben
                \Phi_0(x, \varepsilon)=\left( I_L(\varepsilon) + \mathcal{O}(x-x_L)\right)\,
                (x-x_L)^{\frac{1}{2} \Lambda + \frac{1}{2 x_L} B}\,(x-x_L)^{T_L}
                \een
                in a neighborhood of $x_L$ which does not contain the other singular points, and
                    \ben
                    \Phi_0(x, \varepsilon)=\left( I_R(\varepsilon) + \mathcal{O}(x-x_R)\right)\,
                    (x-x_R)^{\frac{1}{2} \Lambda + \frac{1}{2 x_R} B}\,(x-x_R)^{T_R}
                    \een
                    in a neighborhood of $x_R$ which does not contain the other singular points.
                    
                    Similarly,  during a double resonance 
                    the fundamental matrix $\Phi_{\infty}(x, \varepsilon)$ of the Heun type  equation
                    is represented near the singular points $x_{jj},\,j=R, L$ as follows,
                    \ben
                    \Phi_{\infty}(x, \varepsilon)=\left( I_{LL}(\varepsilon) + \mathcal{O}(x-x_{LL})\right)\,
                    (x-x_{LL})^{-\frac{x_{LL}}{2} G}\,(x-x_{LL})^{T_{LL}}
                    \een
                    in a neighborhood of $x_{LL}$ which does not contain the other singular points, and
                    \ben
                    \Phi_{\infty}(x, \varepsilon)=\left( I_{RR}(\varepsilon) + \mathcal{O}(x_{RR}-x)\right)\,
                    (x_{RR}-x)^{-\frac{x_{RR}}{2} G}\,(x_{RR}-x)^{T_{RR}}
                    \een
                    in a neighborhood of $x_{RR}$ which does not contain the other singular points.
                    The matrices $I_j(\varepsilon) + \mathcal{O}(x-x_j)$ are holomorphic matrix-functions
                    near the points $x_j, j=L, R, LL, RR$, respectively. The matrices $\Lambda, B, G$ are the matrices,
                    associated with the initial equation and defined in the previous section. The matrices $T_j$ and $T_{jj}$
                     are given by
                \be\label{T}
                T_j=\left(\begin{array}{cc}
                	0   &d_j\\
                	0   &0
                \end{array}
                \right),\quad
                      T_{jj}=\left(\begin{array}{cc}
                      	0   &d_{jj}\\
                      	0   &0
                      \end{array}
                      \right)\,.			
                \ee
                The elements $d_j$ and $d_{jj}$ are defined as follows
                \be\label{d}
                d_j=\res \left(\frac{\Phi_2(x, \varepsilon)}{\Phi_1(x, \varepsilon)},\,x=x_j\right),\quad
                     d_{jj}=\res \left(\frac{\Phi_2(x, \varepsilon)}{\Phi_1(x, \varepsilon)},\,x=x_{jj}\right)\,.
                \ee
                \ethe

                \proof
                The elements $\Phi_i(x, \varepsilon),\,i=1, 2$ are the same for both fundamental matrices
                $\Phi_0(x, \varepsilon)$ and $\Phi_{\infty}(x, \varepsilon)$. In a neighborhood of the
                singular point $x_j,\,j=L, R, LL$, which does not contain the other singular point, they have the form
                 $$\,
                  \Phi_i(x, \varepsilon)=(x-x_j)^{m_{i, j}}\,h_{i, j}(x), 
                 \,$$
              and
                     $$\,
                     \Phi_i(x, \varepsilon)=(x_{RR}-x)^{m_{i, RR}}\,h_{i, RR}(x) 
                     \,$$
                in a neighborhood of the point $x_{RR}$, which does not contain the rest singular points. 
                Here $m_{1, L}=-\beta_1/2 \sqrt{\varepsilon},
                m_{2, L}=-\beta_2/2 \sqrt{\varepsilon}-1, m_{1, R}=\beta_1/2 \sqrt{\varepsilon},
                m_{2, R}=\beta_2/2 \sqrt{\varepsilon}-1, m_{1, LL}=\gamma_1/2 \sqrt{\varepsilon},
               m_{2, LL}=\gamma_2/2 \sqrt{\varepsilon}, m_{1, RR}=-\gamma_1/2 \sqrt{\varepsilon},
                m_{2, RR}=-\gamma_2/2 \sqrt{\varepsilon}$.
                The functions $h_{i, j}(x)$ are holomorphic in the same neighborhood of the singular points
                $x_j$.
                
                Consider the element $\Phi_{12}(x, \varepsilon)$ of the matrix $\Phi_0(x, \varepsilon)$.
                During a double resonance of type {\bf A.1} and {\bf A.2} it has the simple form
                       $$\,
                       \Phi_{12}(x, \varepsilon)=(x-x_R)^{m_{2, R}+1}\,u_{12, R}(x),
                       \,$$                 
                       in a neighborhood of the point $x_R$, which does not contain the point $x_L$. 
                       The function $u_{12, R}(x)$ is a holomorphic function
                        in the same neighborhood of the point $x_R$.
                  At the same time  in a neighborhood of the  point $x_L$, which does not contain the
                  point $x_R$, the solution $\Phi_{12}(x, 0)$ has the form   
                 \ben
                   \Phi_{12}(x, \varepsilon)  &=&
                   (x-x_L)^{m_{1, L}}\,h_{1, L}(x)
                   \left[d_L\,\log (x-x_L) + (x-x_L)^{-\frac{\beta_2-\beta_1}{2 \sqrt{\varepsilon}}}\,
                   g_{1, L}(x)\right]\\[0.1ex]
                                              &=&
                    d_L\,(x-x_L)^{m_{1, l}}\,\log (x-x_L)\,h_{1, L}(x) +
                    (x-x_L)^{-\frac{\beta_2}{2\sqrt{\varepsilon}}}\,f_{1, L}(x)\,.                          
                                    \een
                   Here $g_{1, L}(x)$ and $f_{1, L}(x)$ are holomorphic functions
                   in the same neighborhood of the point $x_L$.

                In the same manner during a double resonance of type {\bf A.3} and {\bf A.4} the
                element $\Phi_{12}(x, \varepsilon)$ of $\Phi_0(x, \varepsilon)$ is represented in a neighborhood 
                of the point $x_L$, which does not contain the point $x_R$,
                as
                      $$\,
                      \Phi_{12}(x, \varepsilon)=(x-x_L)^{m_{2, L}+1}\,v_{12, L}(x),
                      \,$$
                      where $v_{12, L}(x)$ is a holomorphic function in the same neighborhood of the point $x_L$.
                At the same time in a neighborhood of the point $x_R$, which does not contain the point $x_L$,  this element is represented
                      as
                 \ben
                   \Phi_{12}(x, \varepsilon) &=&
                   (x-x_R)^{m_{1, R}}\,h_{1, R}(x)
                   \left[d_R\,\log(x-x_R) + (x-x_R)^{\frac{\beta_2-\beta_1}{2 \sqrt{\varepsilon}}}\,
                   g_{1, R}(x)\right]=\\[0.1ex]
                                              &=&
                   d_R\,(x-x_R)^{m_{1, R}}\,\log(x-x_R)\,h_{1, R}(x)+
                   (x-x_R)^{\frac{\beta_2}{2 \sqrt{\varepsilon}}}\,f_{1, R}(x)\,.
                 \een
                 Here $g_{1, R}(x)$ and $f_{1, R}(x)$ are holomorphic functions in the same neighborhood of $x_R$.

                 Consider now the element $\Phi_{12}(x, \varepsilon)$ of the fundamental matrix 
                 $\Phi_{\infty}(x, \varepsilon)$. During a double resonance of type {\bf A.1} and {\bf A.3}
                 it is represented simply in a neighborhood of  the point $x_{RR}$, which does not contain the point
                 $x_{LL}$, as
                   $$\,
                      \Phi_{12}(x, \varepsilon)=(x_{RR}-x)^{m_{2, RR}+1}\,w_{12, RR}(x)\,.
                   \,$$
                   Here $w_{12, RR}$ is a holomorphic function in the same neighborhood.
                   At the same time in a neighborhood of the point $x_{LL}$ which does not contain the point
                   $x_{RR}$ this element is represented as 
                  \ben
                    \Phi_{12}(x, \varepsilon)
                         &=&
                      (x-x_{LL})^{m_{1, LL}}\,h_{1, LL}(x)
                      \left[d_{LL}\,\log(x-x_{LL}) + (x-x_{LL})^{\frac{\gamma_2-\gamma_1}{2 \sqrt{\varepsilon}}+1}\,
                      g_{1, LL}(x)\right]   \\[0.15ex]
                       &=&
                       d_{LL}(x-x_{LL})^{m_{1, LL}}\log(x-x_{LL}) h_{1, LL}(x)+
                       (x-x_{LL})^{\frac{\gamma_2}{2 \sqrt{\varepsilon}}+1} f_{1, LL}(x)\,.
                  \een
               Here $g_{1, LL}(x)$ and $f_{1, LL}(x)$ are holomorphic functions in the neighborhood of the point $x_{LL}$.

                     In the same manner during a double resonance of type {\bf A.2} and {\bf A.4} the
                     element $\Phi_{12}(x, \varepsilon)$ of the matrix $\Phi_{\infty}(x, \varepsilon)$ is represented 
                     in a neighborhood of the point $x_{RR}$, which does not contain the point $x_{LL}$, as
                     \ben
                     \Phi_{12}(x, \varepsilon) &=&
                     (x_{RR} - x)^{m_{1, RR}}\,h_{1, RR}(x)
                     \left[d_{RR}\,\log(x_{RR}-x) + (x_{RR}-x)^{-\frac{\gamma_2-\gamma_1}{2 \sqrt{\varepsilon}}+1}\,
                     g_{1, RR}(x)\right]\\[0.1ex]
                     &=&
                     d_{RR}\,(x_{RR}-x)^{m_{1, RR}}\,\log(x_{RR}-x)\,h_{1, RR}(x)+
                     (x_{RR}-x)^{-\frac{\gamma_2}{2 \sqrt{\varepsilon}}+1}\,f_{1, RR}(x)\,.
                     \een
                     Here $g_{1, RR}(x)$ and $f_{1, RR}(x)$ are holomorphic functions in the same neighborhood of $x_{RR}$.
                     At the same time in a neighborhood of the point $x_{LL}$, which does not contain the point $x_{RR}$  this element has the simple form
                     $$\,
                     \Phi_{12}(x, \varepsilon)=(x-x_{LL})^{m_{2, LL}+1}\,q_{12, LL}(x)\,,
                     \,$$
                     where $q_{12, LL}(x)$ is a holomorphic function in the same neighborhood of the point $x_{LL}$.
                     
                     Any solution $\Phi_{12}(x, \varepsilon)$ that contains logarithmic terms is obtained from a solution
                     $\Phi_{12}(x, \varepsilon)$ defined by the corresponding integral after its analytic continuation
                     along a curve that does not cross a singular point except the base point $\epsilon$.
                     
                     This ends the proof. \qed
                
                   \bre{ex}
                     The matrices $\Lambda/2 + B/2 x_j$ and $-x_{jj}\,G/2$ are  expressed in terms of the
                     characteristic exponents $\rho^j_i$ and $\rho^{jj}_i, i=1, 2, j=R, L$ as follows,
                      \ben
                       \frac{1}{2}\,\Lambda + \frac{1}{2 x_j}\,B=
                         \left(\begin{array}{cc}
                         	 \rho^j_1     &0\\
                         	  0           &\rho^j_2-1
                         \end{array}
                         \right),\quad
                         -\frac{x_{jj}}{2}\,G=
                                   \left(\begin{array}{cc}
                                   	\rho^{jj}_1     &0\\
                                   	0           &\rho^{jj}_2-1
                                   \end{array}
                                   \right)\,.
                      \een
                      
                   \ere

                This local structure of the space of solutions of the Heun type equation defines
                a corresponding monodromy representation.
                \bth{M}
                 During a double resonance 
                the monodromy matrices $M_j(\varepsilon),\,j=R, L$ of the Heun type equation
                with respect to the fundamental matrix $\Phi_0(x, \varepsilon)$, introduced by \thref{expl-A} and
                \thref{F} are
                given by,
                \be\label{M-f}
                M_j(\varepsilon)=e^{\pi\,i\,(\Lambda + \frac{1}{x_j}\,B)}\,e^{2 \pi\,i\,T_j}=
                e^{2 \pi\,i\,T_j}\,e^{\pi\,i\,(\Lambda + \frac{1}{x_j}\,B)}\,.
                \ee
                Similarly,  during a double resonance 
                the monodromy matrices $M_{jj}(\varepsilon),\,j=R, L$ of the Heun type equation
                with respect to the fundamental matrix $\Phi_{\infty}(x, \varepsilon)$, introduced by \thref{expl-A} and
                \thref{F} are
                given by,
                \be\label{M-in}
                M_{jj}(\varepsilon)=e^{-\pi\,i\,x_{jj}\,G}\,e^{2 \pi\,i\,T_{jj}}=
                e^{2 \pi\,i\,T_{jj}}\,e^{-\pi\,i\,x_{jj}\,G}\,.
                \ee
                \ethe
                
                \proof
                The statement for the matrices $M_j(\varepsilon), j=R, L$ follows from \thref{F}, \deref{mm} and the observation
                that during a double resonance the matrices $e^{\pi\,i (\Lambda + \frac{1}{x_j} B)}$ and 
                $(x-x_j)^{T_j}$ commute. 
                
                The statement for the matrices $M_{jj}, j=R, L$ again follows from \thref{F}, \deref{mm} and the observation that
                during a double resonance the matrices $e^{-\pi\,i\,x_{jj} G}$ and $(x-x_{jj})^{T_{jj}}$
                (resp. $(x_{RR}-x)^{T_{RR}}$) commute.
                \qed

                \bre{mixed}
                  During a double resonance the first column of the matrices $\Phi_0(x, \varepsilon)$ and $\Phi_{\infty}(x, \varepsilon)$
                  is a eigenvector of the monodromy operators that are represented by the monodromy matrices $M_j(\varepsilon)$
                  and $M_{jj}(\varepsilon)$, respectively. The corresponding eigenvalues are $e^{2 \pi\,i\,\rho^j_1}$ and
                  $e^{2 \pi\,i\,\rho^{jj}_1}$, respectively. The numbers $d_j$ and $d_{jj}$ (when they are different from zero)
                  block the second column of being eigenvector of the monodromy operators. Following Duval \cite{Du} and
                  Hurtubise et al. \cite{HLR} we call such a basis of solutions a "mixed" basis of solutions.
                \ere

       We finish this indentation explaining the computation of the number $d_{RR}$.         
       In this paper we will applay the same technique used in \cite{St1} to compute the residue of the complex
       function $f(x)$. Namely, let $x_0$ be a pole for $f(x)$ of order $s$, i.e.
        $$\,
          f(x)=\frac{\varphi(x)}{(x-x_0)^s}\,,
        \,$$
  where $\varphi(x)$ is a holomorphic function in a neighborhood of $x_0$ and $\varphi(x_0) \neq 0$. Then
    \be\label{R}
      \res \left(f(x),\,x=x_0\right)=\frac{\varphi^{(s-1)}(x_0)}{(s-1)!}\,.
    \ee
    Applying formula \eqref{R} the number $d_{RR}$ becomes the coefficient before $\frac{1}{x-x_{RR}}$ in
    $\frac{\Phi_2(x, \varepsilon)}{\Phi_1(x, \varepsilon)}$. But after integration we get
     $$\,
      \int_{-1/\sqrt{\varepsilon}}^{x}\,\, \frac{d_{RR}}{t-\frac{1}{\sqrt{\varepsilon}}} \, d t=
      -d_{RR} \int_{-1/\sqrt{\varepsilon}}^{x} \,\,\frac{ d t}{\frac{1}{\sqrt{\varepsilon}}-t} =
      d_{RR}\,\log\left(\frac{1}{\sqrt{\varepsilon}}-x\right) -d_{RR}\,\log \frac{2}{\sqrt{\varepsilon}}\,.
     \,$$
     So, the number
      $$\,
       d_{RR}=\res \left(\frac{\Phi_2(x, \varepsilon)}{\Phi_1(x, \varepsilon)},\,x=x_{RR}\right)
      \,$$            
      is exactly the coefficient before $\log (x_{RR}-x)$.
               
  To the end of this section we will compute the non-zero numbers $d_j$ during every type of double
  resonance. The numbers $d_j,\,j=R, L$ will be computed with respect to the fundamental matrix
  $\Phi_0(x, \varepsilon)$. The numbers $d_{jj},\,j=R, L$ will be computed with respect to the
  fundamental matrix $\Phi_{\infty}(x, \varepsilon)$.

  %%%%%%%%%%%%%%%%%%%%%%%%%%%%%%%%%%%%%%%%%%%%%%%%%%%%%%%%%%%
  %d
  %%%%%%%%%%%%%%%%%%%%%%%%%%%%%%%%%%%%%%%%%%%%%%%%%%%%%%%%%%%%%

    We start with the computation of the numbers $d_L$ and $d_{LL}$ during a double resonance of
    type {\bf A.1}.

    \bth{A.1-M}
		  Let $\frac{\beta_2-\beta_1}{2 \sqrt{\varepsilon}}\in\NN$.		  
			Then with respect to the fundamental matrix $\Phi_0(x, \varepsilon)$ for the number $d_L$ 
			 we have
			 \begin{enumerate}
			 
			 \item\,
			 If $\gamma_1=\gamma_2$ then $d_L=0$.
			 
			 \item\,
			 If $\frac{\gamma_1-\gamma_2}{2 \sqrt{\varepsilon}}\in\NN$ then
				\ben
				  d_L=-\frac{\gamma_1-\gamma_2}{2 \sqrt{\varepsilon}}\,
					\left(\frac{1+\varepsilon}{1-\varepsilon}\right)^{\frac{\gamma_1-\gamma_2}{2 \sqrt{\varepsilon}}}
					\sum_{k=1}^{\frac{\beta_2-\beta_1}{2 \sqrt{\varepsilon}}}
					\frac{(2 \sqrt{\varepsilon})^{k-1}}{\Gamma(k) \Gamma(k+1)}
					\frac{\Gamma(\frac{\beta_2-\beta_1}{2 \sqrt{\varepsilon}})}
					     {\Gamma(\frac{\beta_2-\beta_1}{2 \sqrt{\varepsilon}}-k+1)}
					\left(\frac{\sqrt{\varepsilon}}{1+\varepsilon}\right)^k\,A\,,		
				\een
				where
				$$\,
				  A=\sum_{s=0}^k
					  \left(\begin{array}{l}
						   k\\
							 s
							    \end{array}
						\right)
				\frac{\Gamma(\frac{\gamma_1-\gamma_2}{2 \sqrt{\varepsilon}}+s)}
						 {\Gamma(\frac{\gamma_1-\gamma_2}{2 \sqrt{\varepsilon}}+1-k+s)}
						 \left(\frac{1+\varepsilon}{1-\varepsilon}\right)^s\,.			
				\,$$
				\end{enumerate}
		\ethe

   \proof
	
  	 According to \eqref{d} the number $d_L$ is defined by
		$$\,
		   d_L=\res\left(
			 \frac{(x-\sqrt{\varepsilon})^{\frac{\beta_2-\beta_1}{2 \sqrt{\varepsilon}}-1}}
			      {(x+\sqrt{\varepsilon})^{\frac{\beta_2-\beta_1}{2 \sqrt{\varepsilon}}+1}}
			\left(\frac{\frac{1}{\sqrt{\varepsilon}}-x}{\frac{1}{\sqrt{\varepsilon}}+x}
			\right)^{\frac{\gamma_1-\gamma_2}{2 \sqrt{\varepsilon}}}\,;\,
						x=x_L\right)\,.
		\,$$
	
	Then for the number $d_L$ we find that
  	\ben
	  d_L=
		\frac{1}{\left(\frac{\beta_2-\beta_1}{2 \sqrt{\varepsilon}}\right)!}
		D^{\frac{\beta_2-\beta_1}{2 \sqrt{\varepsilon}}}
				\left((x-\sqrt{\varepsilon})^{\frac{\beta_2-\beta_1}{2 \sqrt{\varepsilon}}-1}
				\left(\frac{\frac{1}{\sqrt{\varepsilon}}+x}{\frac{1}{\sqrt{\varepsilon}}-x}
				\right)^{\frac{\gamma_2-\gamma_1}{2 \sqrt{\varepsilon}}}\right)_
				{x=x_L}\,,
				\een
				where we have denoted by $D$ the differential operator $\frac{d}{d x}$.
				Obviously when $\gamma_1=\gamma_2$ the number $d_L$ becomes zero.
				Let $\gamma_1 \neq \gamma_2$.
				Now applying the  Leibnitz's rule
				\ben
		      D^n\left(f(x)\,g(x)\right)=
				\sum_{k=0}^n
				   \left(\begin{array}{c}
					      n\\
								k
								 \end{array}
						\right)
						D^{n-k}(f(x))\,D^k(g(x))\,,		
				\een
		we obtain the wanted expression for $d_L$.
	\qed

    \bth{A.1-M-inf}
		  Let  $\frac{\beta_2-\beta_1}{2 \sqrt{\varepsilon}} \in \NN$.
		  Then with respect to the fundamental matrix $\Phi_{\infty}(x, \varepsilon)$
		the number $d_{LL}$ is given as follows,
		  \begin{enumerate}
		  	
		  \item\,
		  If $\gamma_1=\gamma_2$ then $d_{LL}=0$.

		  \item\,
		  If $\frac{\gamma_1-\gamma_2}{2 \sqrt{\varepsilon}}\in\NN$ then
                \ben
                 d_{LL}=
                 \frac{2 \sqrt{\varepsilon}}{\beta_2-\beta_1}
                 \frac{1}{1-\varepsilon^2}
                           \left(\frac{1+\varepsilon}{1-\varepsilon}\right)^{\frac{\beta_2-\beta_1}{2 \sqrt{\varepsilon}}}
                    \sum_{k=0}^{\frac{\gamma_1-\gamma_2}{2 \sqrt{\varepsilon}} - 1}
                          \frac{1}{k!\,(k+1)!} \frac{(2 \sqrt{\varepsilon})^{1+k} (\sqrt{\varepsilon})^k}{(1+\varepsilon)^k}
                   \frac{\Gamma(\frac{\gamma_1-\gamma_2}{2 \sqrt{\varepsilon}} +1)}
                         {\Gamma(\frac{\gamma_1-\gamma_2}{2 \sqrt{\varepsilon}}-k)}\,A,
               \een
	  where
             $$\,
                         A=\sum_{s=0}^k
                               \left(\begin{array}{c}
                                            k\\
                                            s
                                       \end{array}\right)
                 \frac{\Gamma(\frac{\beta_2-\beta_1}{2 \sqrt{\varepsilon}}+1+s)}
                 {\Gamma(\frac{\beta_2-\beta_1}{2 \sqrt{\varepsilon}}-k+s)}
                  \left(\frac{1+\varepsilon}{1-\varepsilon}\right)^s.
              \,$$
              \end{enumerate}
 \ethe	
	
                  \proof
              Let $\gamma_1 \neq \gamma_2$.
              To compute the number $d_{LL}$ we apply the same procedure as in the construction of the formal solution 
              of the initial equation  near $x=\infty$. 
             Instead of the  solution $\Phi_{12}(x, \varepsilon)$ given by \thref{expl-A} we will use the solution 
             $\Phi_{12}(x, \varepsilon)$ built by the operators $L_{j, \varepsilon}$ from
                 by \eqref{the}. 
                 Observe that the limit $\varepsilon \rightarrow 0$ takes these operators $L_{j, \varepsilon}$  into the operators 
               $L_j$ defined by \eqref{tr-op}.
            Observe also  that after the transformation $x=1/t$ the function $\Phi_{12}(x, \varepsilon)$ becomes a particular solution
                  of the following non-homogeneous equation
                      \be\label{inf}
                             \dot{y}(t)     &+&
                                     \left[\frac{\beta_1}{2 \sqrt{\varepsilon}} 
                                \left(-\frac{1}{t-\frac{1}{\sqrt{\varepsilon}}} + \frac{1}{t+\frac{1}{\sqrt{\varepsilon}}}\right)
                                            +\frac{\gamma_1}{2 \sqrt{\varepsilon}}
                               \left(\frac{1}{t-\sqrt{\varepsilon}} - \frac{1}{t+\sqrt{\varepsilon}}\right)\right]\,y(t)
                               \nonumber\\[0.25ex]
                                            &=&
                             -\frac{1}{\varepsilon} \frac{\left(\frac{1}{\sqrt{\varepsilon}}-t\right)^{\frac{\beta_2}{2 \sqrt{\varepsilon}} -1}}
                                           {\left(\frac{1}{\sqrt{\varepsilon}} + t\right)^{\frac{\beta_2}{2 \sqrt{\varepsilon}} + 1}}
                              \left(\frac{t+\sqrt{\varepsilon}}{t-\sqrt{\varepsilon}}\right)
                              ^{\frac{\gamma_2}{2 \sqrt{\varepsilon}}}.
                    \ee
                  Then the number $d_{LL}$ of the perturbed equation becomes the number $d_L$ of the transformed equation. As a result, we redefine
                      the number $d_{LL}$ as
                    \ben
                    d_{LL}=-\frac{1}{\varepsilon}\,
                                  \res\left(\left(\frac{t-\sqrt{\varepsilon}}{t+\sqrt{\varepsilon}}\right)^{\frac{\gamma_1-\gamma_2}{2 \sqrt{\varepsilon}}}
                                 \frac{(\frac{1}{\sqrt{\varepsilon}} -t)^{\frac{\beta_2-\beta_1}{2 \sqrt{\varepsilon}} -1}}
                                       {(\frac{1}{\sqrt{\varepsilon}} + t)^{\frac{\beta_2-\beta_1}{2 \sqrt{\varepsilon}} +1}}\,;\,t=-\sqrt{\varepsilon}\right).
                    \een
                    Then
                     \ben
                     d_{LL}     &=& 
                   -\frac{1}{\varepsilon} \frac{1}{\left(\frac{\gamma_1-\gamma_2}{2 \sqrt{\varepsilon}} -1\right)!}\\[0.15ex]
                                    &\times&
                        \sum_{k=0}^{\frac{\gamma_1-\gamma_2}{2 \sqrt{\varepsilon}} -1}
                           \left(\begin{array}{c}
                                     \frac{\gamma_1-\gamma_2}{2 \sqrt{\varepsilon}} -1 \\
                                             k
                                  \end{array}\right)
                       \left(D^{\frac{\gamma_1-\gamma_2}{2 \sqrt{\varepsilon}} -1 -k} \left((t-\sqrt{\varepsilon})^{\frac{\gamma_1-\gamma_2}{2 \sqrt{\varepsilon}}}\right)
                            D^k \left(
                                  \frac{(\frac{1}{\sqrt{\varepsilon}} -t)^{\frac{\beta_2-\beta_1}{2 \sqrt{\varepsilon}} -1}}
                                       {(\frac{1}{\sqrt{\varepsilon}} + t)^{\frac{\beta_2-\beta_1}{2 \sqrt{\varepsilon}} +1}}\right)\right)_{t=-\sqrt{\varepsilon}},
                    \een
                  where we have denoted by $D$ the differential operator $\frac{d}{d t}$.
                Now it is not difficult to show that the number $d_{LL}$ has the pointed form.
                
                When $\gamma_1=\gamma_2$ the solution $\Phi_{12}(x, \varepsilon)$ given by \thref{expl-A} does not contain
                $\log(x-x_{LL})$ since it does not have a singularity at the point $x_{LL}$. 

                This ends the proof.
     \qed

                \bre{inf-sol}
                   We will use the same equation \eqref{inf} to compute the numbers $d_{jj}$ for the rest
                   double resonances.
                 \ere

          The next two theorems give the numbers $d_L$ and $d_{RR}$ during a double resonance of type
          {\bf A.2}.
          
         \bth{A.2-M}
         Let  $\frac{\beta_2-\beta_1}{2 \sqrt{\varepsilon}}\in\NN$. 
         Then  with respect to the fundamental matrix $\Phi_0(x, \varepsilon)$
         for the number $d_L$ we have
             
             \begin{enumerate}
             	
             \item\,
             If $\gamma_1=\gamma_2$ then $d_L=0$.
             
             \item\,
             If   $\frac{\gamma_2-\gamma_1}{2 \sqrt{\varepsilon}}\in\NN$ then         
         \ben
         d_L= \frac{\gamma_2-\gamma_1}{2 \sqrt{\varepsilon}}
         \left(\frac{1-\varepsilon}{1+\varepsilon}\right)^{\frac{\gamma_2-\gamma_1}{2 \sqrt{\varepsilon}}}
         \sum_{k=1}^{\frac{\beta_2-\beta_1}{2 \sqrt{\varepsilon}}}
         \frac{(- 2 \sqrt{\varepsilon})^{k-1}}{(k-1)!\,k!} 
         \left(\frac{\sqrt{\varepsilon}}{1-\varepsilon}\right)^k 
         \frac{\Gamma(\frac{\beta_2-\beta_1}{2 \sqrt{\varepsilon}})}
         {\Gamma(\frac{\beta_2-\beta_1}{2 \sqrt{\varepsilon}}-k+1)}\,A,
         \een
         where
         $$\,
         A=\sum_{s=0}^k
         \left(\begin{array}{c}
         k\\
         s
         \end{array}\right)
         \frac{\Gamma(\frac{\gamma_2-\gamma_1}{2 \sqrt{\varepsilon}}+s)}
         {\Gamma(\frac{\gamma_2-\gamma_1}{2 \sqrt{\varepsilon}}-k+s+1)}
         \left(\frac{1-\varepsilon}{1+\varepsilon}\right)^s.
         \,$$
          \end{enumerate}
         \ethe

              \bth{A.2-M-inf}
              Let  $\frac{\beta_2-\beta_1}{2 \sqrt{\varepsilon}}\in\NN$.
              Then  with respect to the fundamental matrix $\Phi_{\infty}(x, \varepsilon)$
              the number $d_{RR}$ is defined as follows:
              
                 \begin{enumerate}
                 	
                 \item\,
                 If $\gamma_1=\gamma_2$ then $d_{RR}=0$.
                 
                 \item\,
                 If $\frac{\gamma_2-\gamma_1}{2 \sqrt{\varepsilon}}\in\NN$ then
              \ben
              d_{RR}= -\frac{2 \sqrt{\varepsilon}}{\beta_2-\beta_1}
              \frac{1}{1-\varepsilon^2}
              \left(\frac{1-\varepsilon}{1+\varepsilon}\right)^{\frac{\beta_2-\beta_1}{2 \sqrt{\varepsilon}}}
              \sum_{k=0}^{\frac{\gamma_2-\gamma_1}{2 \sqrt{\varepsilon}}-1}
              \frac{(-1)^k ( 2 \sqrt{\varepsilon})^{k+1}}{k!\,(k+1)!} 
              \left(\frac{\sqrt{\varepsilon}}{1-\varepsilon}\right)^k 
              \frac{\Gamma(\frac{\gamma_2-\gamma_1}{2 \sqrt{\varepsilon}}+1)}
              {\Gamma(\frac{\gamma_2-\gamma_1}{2 \sqrt{\varepsilon}}-k)}\,A,
              \een
              where
              $$\,
              A=\sum_{s=0}^k
              \left(\begin{array}{c}
              k\\
              s
              \end{array}\right)
              \frac{\Gamma(\frac{\beta_2- \beta_1}{2 \sqrt{\varepsilon}}+s+1)}
              {\Gamma(\frac{\beta_2-\beta_1}{2 \sqrt{\varepsilon}}-k+s)}
              \left(\frac{1-\varepsilon}{1+\varepsilon}\right)^s.
              \,$$
                 \end{enumerate}
              \ethe	
              
   In the next two theorems we present explicitly the numbers $d_R$ and $d_{LL}$ during
   a double resonance of type {\bf A.3}.

        \bth{A.3-M}
        Let  $\frac{\beta_1-\beta_2}{2 \sqrt{\varepsilon}}\in\NN$. 
        Then  with respect to the fundamental matrix $\Phi_0(x, \varepsilon)$
        for the number $d_R$ we have
          \begin{enumerate}
          	
          \item\,
          If $\gamma_1=\gamma_2$ then $d_R=0$.
          
          \item\,
          If $\frac{\gamma_1-\gamma_2}{2 \sqrt{\varepsilon}}\in\NN$ then
          
        \ben
        d_R= \frac{\gamma_1-\gamma_2}{2 \sqrt{\varepsilon}}
        \left(\frac{1-\varepsilon}{1+\varepsilon}\right)^{\frac{\gamma_1-\gamma_2}{2 \sqrt{\varepsilon}}}
        \sum_{k=1}^{\frac{\beta_1-\beta_2}{2 \sqrt{\varepsilon}}}
        \frac{(-1)^k ( 2 \sqrt{\varepsilon})^{k-1}}{(k-1)!\,k!} 
        \left(\frac{\sqrt{\varepsilon}}{1-\varepsilon}\right)^k 
        \frac{\Gamma(\frac{\beta_1-\beta_2}{2 \sqrt{\varepsilon}})}
        {\Gamma(\frac{\beta_1-\beta_2}{2 \sqrt{\varepsilon}}-k+1)}\,A,
        \een
        where
        $$\,
        A=\sum_{s=0}^k
        \left(\begin{array}{c}
        k\\
        s
        \end{array}\right)
        \frac{\Gamma(\frac{\gamma_1-\gamma_2}{2 \sqrt{\varepsilon}}+s)}
        {\Gamma(\frac{\gamma_1-\gamma_2}{2 \sqrt{\varepsilon}}-k+s+1)}
        \left(\frac{1-\varepsilon}{1+\varepsilon}\right)^s.
        \,$$
           \end{enumerate}
        \ethe

        \bth{A.3-M-inf}
        Let $\frac{\beta_1-\beta_2}{2 \sqrt{\varepsilon}}\in\NN$. 
        Then  with respect to the fundamental matrix $\Phi_{\infty}(x, \varepsilon)$
        for the number $d_{LL}$ we have
        
          \begin{enumerate}
          	
           \item\,
           If $\gamma_1=\gamma_2$ then $d_{LL}=0$.
           
           \item\,
           If $\frac{\gamma_1-\gamma_2}{2 \sqrt{\varepsilon}}\in\NN$ then
        \ben
        d_{LL}= 
        -\frac{2 \sqrt{\varepsilon}}{\beta_1-\beta_2}
        \frac{1}{1-\varepsilon^2}
        \left(\frac{1-\varepsilon}{1+\varepsilon}\right)^{\frac{\beta_1-\beta_2}{2 \sqrt{\varepsilon}}}
        \sum_{k=0}^{\frac{\gamma_1-\gamma_2}{2 \sqrt{\varepsilon}}-1}
        \frac{ ( -2 \sqrt{\varepsilon})^{k+1}}{k!\,(k+1)!} 
        \left(\frac{\sqrt{\varepsilon}}{1-\varepsilon}\right)^k 
        \frac{\Gamma(\frac{\gamma_1-\gamma_2}{2 \sqrt{\varepsilon}}+1)}
        {\Gamma(\frac{\gamma_1-\gamma_2}{2 \sqrt{\varepsilon}}-k)}\,A,
        \een
        where
        $$\,
        A=\sum_{s=0}^k
        \left(\begin{array}{c}
        k\\
        s
        \end{array}\right)
        \frac{\Gamma(\frac{\beta_1- \beta_2}{2 \sqrt{\varepsilon}}+s+1)}
        {\Gamma(\frac{\beta_1-\beta_2}{2 \sqrt{\varepsilon}}-k+s)}
        \left(\frac{1-\varepsilon}{1+\varepsilon}\right)^s.
        \,$$
           \end{enumerate}
        \ethe

  The last two theorems of this subsection give the numbers $d_R$ and $d_{RR}$ during a double resonance of type
  {\bf A.4}.

               \bth{A.4-M}
               Let $\frac{\beta_1-\beta_2}{2 \sqrt{\varepsilon}}\in\NN$. 
              Then  with respect to the fundamental matrix $\Phi_0(x, \varepsilon)$
              for the number $d_R$ we have
              
               \begin{enumerate}
               	
               \item\,
                  If $\gamma_1=\gamma_2$ then $d_R=0$.
                  
               \item\,
               If    $\frac{\gamma_2-\gamma_1}{2 \sqrt{\varepsilon}}\in\NN$ then 
               \ben
               d_R= \frac{\gamma_2-\gamma_1}{2 \sqrt{\varepsilon}}
               \left(\frac{1+\varepsilon}{1-\varepsilon}\right)^{\frac{\gamma_2-\gamma_1}{2 \sqrt{\varepsilon}}}
               \sum_{k=1}^{\frac{\beta_1-\beta_2}{2 \sqrt{\varepsilon}}}
               \frac{(2 \sqrt{\varepsilon})^{k-1}}{(k-1)!\,k!} 
               \left(\frac{\sqrt{\varepsilon}}{1+\varepsilon}\right)^k 
               \frac{\Gamma(\frac{\beta_1-\beta_2}{2 \sqrt{\varepsilon}})}
               {\Gamma(\frac{\beta_1-\beta_2}{2 \sqrt{\varepsilon}}-k+1)}\,A,
               \een
               where
               $$\,
               A=\sum_{s=0}^k
               \left(\begin{array}{c}
               k\\
               s
               \end{array}\right)
               \frac{\Gamma(\frac{\gamma_2-\gamma_1}{2 \sqrt{\varepsilon}}+s)}
               {\Gamma(\frac{\gamma_2-\gamma_1}{2 \sqrt{\varepsilon}}-k+s+1)}
               \left(\frac{1+\varepsilon}{1-\varepsilon}\right)^s.
               \,$$
                \end{enumerate}
               \ethe

               \bth{A.4-M-inf}
               Let  $\frac{\beta_1-\beta_2}{2 \sqrt{\varepsilon}}\in\NN$.
               Then  with respect to the fundamental matrix $\Phi_{\infty}(x, \varepsilon)$
              for the number $d_{RR}$ we have
              
                \begin{enumerate}
                	
                 \item\,
                  If $\gamma_1=\gamma_2$ then $d_{RR}=0$.
                  
                  \item\,
                  If $\frac{\gamma_2-\gamma_1}{2 \sqrt{\varepsilon}}\in\NN$ then
                  
               \ben
               d_{RR}= 
                 -\frac{2 \sqrt{\varepsilon}}{\beta_1-\beta_2}
               \frac{1}{1-\varepsilon^2}
               \left(\frac{1+\varepsilon}{1-\varepsilon}\right)^{\frac{\beta_1-\beta_2}{2 \sqrt{\varepsilon}}}
               \sum_{k=0}^{\frac{\gamma_2-\gamma_1}{2 \sqrt{\varepsilon}}-1}
               \frac{(2 \sqrt{\varepsilon})^{k+1}}{k!\,(k+1)!} 
               \left(\frac{\sqrt{\varepsilon}}{1+\varepsilon}\right)^k 
               \frac{\Gamma(\frac{\gamma_2-\gamma_1}{2 \sqrt{\varepsilon}}+1)}
               {\Gamma(\frac{\gamma_2-\gamma_1}{2 \sqrt{\varepsilon}}-k)}\,A,
               \een
               where
               $$\,
               A=\sum_{s=0}^k
               \left(\begin{array}{c}
               k\\
               s
               \end{array}\right)
               \frac{\Gamma(\frac{\beta_1-\beta_2}{2 \sqrt{\varepsilon}}+s+1)}
               {\Gamma(\frac{\beta_1-\beta_2}{2 \sqrt{\varepsilon}}-k+s)}
               \left(\frac{1+\varepsilon}{1-\varepsilon}\right)^s.
               \,$$
                 \end{enumerate}
               \ethe

	%%%%%%%%%%%%%%%%%%%%%%%%%%%%%%%%%%%%%%%%%%%%%%%%%%%%%%
	% Main results
	%%%%%%%%%%%%%%%%%%%%%%%%%%%%%%%%%%%%%%%%%%%%%%%%
	\section{ Main results }

   In this section we establish the main result of this paper. Our goal is to connect by a radial limit $\sqrt{\varepsilon} \rightarrow 0$ the monodromy 
   matrices $M_j(\varepsilon), j=R, L$ and $M_{jj}(\varepsilon), j=R, L$ of the Heun type equation with the Stokes matrices
  of the initial equation. As we mentioned in the introduction in this paper we consider the initial equation with two
  different fundamental matrix solutions without studying the connection between them. These are 
  the fundamental solution at the origin and the fundamental solution at the infinity point. In the same
  manner we consider the Heun type equation with two different fundamental matrix solutions
  $\Phi_0(x, \varepsilon)$ and $\Phi_{\infty}(x, \varepsilon)$ again without studying their connection.
  So we can image that we deal with two pairs of initial and perturbed equation. The initial equation
  of the first pair has an irregular singularity at the origin of Poincar\`{e} rank 1 and we
  split it into two Fuchsian singularities $x_L$ and $x_R$. The initial equation of the second pair has an
  irregular singularity at $x=\infty$ again of Poincar\`{e} rank 1 and we split it into two new
  Fuchsian singular points $x_{LL}$ and $x_{RR}$.

   In next two paragraphs we express the so called unfolded Stokes matrices of each imaginary pair as a
   part of  monodromy matrices of the same pair during a double resonance.
  
    \subsection{The unfolded Stokes matrices $St_j(\varepsilon)$ as the matrices $e^{2 \pi\,i\,T_j}$}
    
    As we saw in Section 3 the initial equation has only one singular direction $\theta=\arg (\beta_1-\beta_2)$
    at the origin. Following \cite{CL-CR, CL-CR1}
    we consider the initial equation together with its Stokes matrix at the origin in the ramified domain
    of the origin
    $\{ x\in\CC\PP^1\,:\,  \arg (\beta_1-\beta_2) -\kappa < \arg (x) < \arg(\beta_1-\beta_2) + k\}$
    where $0 < \kappa < \pi/2$. We cover this domain by two open sectors
      \ben
        \Omega_{1, 0}
           &=&
        \left\{ x=r e^{i \delta}\,|\, 0 < r < \rho,\,
         \arg (\beta_1-\beta_2) - \kappa < \delta < \arg (\beta_1-\beta_2) + \pi + \kappa \right\},\\
         \Omega_{2, 0}
           &=&
          \left\{ x=r e^{i \delta}\,|\, 0 < r < \rho,\,
           - (\arg (\beta_1-\beta_2) + \pi + \kappa) < \delta <  \arg (\beta_1-\beta_2)+ \kappa \right\}\,.\nonumber
      \een
     Denote by $\Omega_R$ and $\Omega_L$ the connection components of the intersection
     $\Omega_{1, 0} \cap \Omega_{2, 0}$. 
     The radius $\rho$ of the sectors $\Omega_{j, 0}$ is so chosen that the only singular points which belong
     to $\Omega_R$ and  $\Omega_L$ to be the points $x_R$ and $x_L$, respectively. Recall that during a double resonance
     $x_R, x_L, \beta_2-\beta_1\in\RR$. So, $\arg (\beta_1-\beta_2)\in\{0, \pi\}$. 
     
     Let us extend the actual fundamental matrix solution $\Phi_0(x, 0)$ introduced by \prref{p-A} to the whole
     $\Omega_{1, 0}\cup \Omega_{2, 0}$.      
     From the sectorial normalization theorem of Sibuya \cite{S} if follows that over the sector 
     $\Omega_{1, 0}$ the matrix $H_1(x)=H^{+}_{\theta}(x)$ 
     is an analytic matrix function asymptotic at the origin to the matrix $\hat{H}(x)$.
      Similarly, over the sector
     $\Omega_{2, 0}$ the matrix $H_2(x)=H^{-}_{\theta}(x)$ 
     is an analytic matrix function asymptotic at the origin to the matrix $\hat{H}(x)$.
     Then the matrix
       $$\,
         \Psi^0_j(x)=\exp (G x)\,H_j(x)\,F(x),\quad j=1, 2
       \,$$
     with the corresponding branch of $F(x)$ is an actual fundamental matrix at 
     the origin over the sector $\Omega_{j, 0}$, respectively. Note that if 
     $\arg (\beta_1-\beta_2)=0$
     we can observe the Stokes phenomenon on $\Omega_R$. If 
     $\arg (\beta_1-\beta_2)=\pi$
     we can observe the Stokes phenomenon over $\Omega_L$. 

    Let us turn around the origin in the positive sense, starting from the sector $\Omega_{1, 0}$.
    On the first sector $\Omega_j, j=R, L$ that we cross we can not observe the Stokes phenomenon.
    On the next sector $\Omega_j, j=R, L$ that we cross we define the Stokes matrix  
    $St_{R}$ (resp. $St_L$) as
   $$\,
     ((\Phi^0_{\theta}(x, 0))^{+})^{-1}\,(\Phi^0_{\theta + 2 \pi}(x, 0))^{-}=
      (\Psi^0_1(x))^{-1}\,\Psi^0_2(x)\,\hat{M}=St_j\,\hat{M}=St^{\theta}_0\,\hat{M}=St^{\theta}_0\,,
   \,$$
  where  $St^{\theta}_0$ is the Stokes matrix defined by \thref{A-S} and $\hat{M}=I_2$. 
  
    At the same time we consider the Heun type equation on the  sectorial domains 
    $\Omega_{1, 0}(\varepsilon)$ and $\Omega_{2, 0}(\varepsilon)$.
    They are obtained from the open sectors $\Omega_{1, 0}$ and 
    $\Omega_{2, 0}$
    by making a cut between the singular points $x_L$ and $x_R$ through the real axis.
    The point $x_0=0$ belongs to this cut.
    When $\varepsilon \rightarrow 0$ the sectorial
    domains $\Omega_{j, 0}(\varepsilon)$ tend to the sectors $\Omega_{j, 0}, j=1, 2$, respectively.
    The sectorial domains $\Omega_{1, 0}(\varepsilon)$ and $\Omega_{2, 0}(\varepsilon)$ intersect in the
    sectors $\Omega_L(\varepsilon)$ and $\Omega_R(\varepsilon)$ and along the cut. The singular points
    $x_j, j=R, L$ belong to this cut. 

     Let us represent the fundamental matrix $\Phi_0(x, \varepsilon)$ of the Heun type equation
     in a slightly different form
      $$\,
        \Phi_0(x, \varepsilon)=G(x, \varepsilon)\,H(x, \varepsilon)\,F(x, \varepsilon)\,,
      \,$$
     where
     $$\,
       G(x, \varepsilon)=
       (x-x_{LL})^{-\frac{x_{LL}}{2} G}
       (x_{RR}-x)^{-\frac{x_{RR}}{2 } G}
     \,$$
     and
     $$\,
       F(x, \varepsilon)=(x-x_L)^{\frac{1}{2} \Lambda + \frac{1}{2 x_L} B}
       (x-x_R)^{\frac{1}{2} \Lambda + \frac{1}{2 x_R} B}\,.
     \,$$
     The matrix $H(x, \varepsilon)=\{h_{ij}(x, \varepsilon)\}_{i, j=1}^2$ is given by
     \ben
       & &
       h_{ii}(x, \varepsilon)=1,\quad h_{ij}(x, \varepsilon)=0\quad \textrm{for}\quad i > j,\\[0.15ex]
       & &
       h_{12}(x, \varepsilon)=(x-x_L)^{\frac{\beta_2-\beta_1}{2 \sqrt{\varepsilon}}+1}
       (x-x_R)^{-\frac{\beta_2-\beta_1}{2 \sqrt{\varepsilon}}+1}
       \int_{\Gamma_0(x, \varepsilon)}
       \frac{\Phi_2(z, \varepsilon)}{\Phi_1(z, \varepsilon)}\, d z\,.
     \een
     The path of integration $\Gamma_0(x, \varepsilon)$ is a path either from $x_R$ to $x$ or
     from $x_L$ to $x$, taken in the direction $\arg(\beta_2-\beta_1)$. Analytic continuation of the
     path $\Gamma_0(x, \varepsilon)$ around the origin in the positive sense yields two branches
     $h_{12}^{-}(x, \varepsilon)$ and $h_{12}^{+}(x, \varepsilon)$ of the element $h_{12}(x, \varepsilon)$
     near the singular direction $\arg (\beta_1-\beta_2)$. The branch $h_{12}^{-}(x, \varepsilon)$
     corresponds to a path taken in the direction $\arg (\beta_1-\beta_2) - \epsilon$, and the branch
     $h_{12}^{+}(x, \varepsilon)$ corresponds to a path taken in the direction $\arg (\beta_1-\beta_2) + \epsilon$.
     Here $\epsilon > 0$ is a small number. When $\Gamma_0(x, \varepsilon)$ crosses the singular direction
     $\arg (\beta_1-\beta_2)$ we rather observe the Stokes phenomenon than the linear monodromy.
     This phenomenon is measured geometrically by the so called unfolded Stokes matrix \cite{CL-CR1}. 
     In concordance with the initial equation we
     fix on the sector $\Omega_{1, 0}(\varepsilon)$ the fundamental matrix of the Heun type equation as
      $$\,
        \Psi^0_1(x, \varepsilon)=G(x, \varepsilon)\,H_1(x, \varepsilon)\,F_1(x, \varepsilon)\,,
      \,$$
      where $H_1(x, \varepsilon)=\{h_{ij}^{+}(x, \varepsilon)\}_{i, j=1}^2$. The matrix
       $F_1(x, \varepsilon)$ is  the branch on $\Omega_{1, 0}(\varepsilon)$ of the matrix
       $F(x, \varepsilon)$. Similarly, on the sector $\Omega_{2, 0}(\varepsilon)$
      we fix the fundamental matrix $\Phi_0(x, \varepsilon)$ as
          $$\,
          \Psi^0_2(x, \varepsilon)=G(x, \varepsilon)\,H_2(x, \varepsilon)\,F_2(x, \varepsilon)\,,
          \,$$
          where $H_2(x, \varepsilon)=\{h_{ij}^{-}(x, \varepsilon)\}_{i, j=1}^2$. 
      The matrix $F_2(x, \varepsilon)$ coincides with $F_1(x, \varepsilon)$ on the sector $\Omega_j(\varepsilon), j=R, L$,
      which contains the non-singular direction $\arg (\beta_2-\beta_1)$. On the sector $\Omega_j(\varepsilon), j=R, L$
      which contains the singular direction $\arg (\beta_1-\beta_2)$ these two matrices again concide since 
      $\hat{M}=I_2$ and therefore $F_2(x, \varepsilon)=F_1(x, \varepsilon)\,\hat{M}=F_1(x, \varepsilon)$. 
      Let is turn around the origin in the positive
      sense starting from the sector $\Omega_{1, 0}(\varepsilon)$ and the solution $\Psi^0_1(x, \varepsilon)$
      on it. When $\Gamma_0(x, \varepsilon)$ crosses the direction $\arg (\beta_2-\beta_1)$ we can not observe the
      Stokes phenomenon on the corresponding $\Omega_j(\varepsilon), j=R, L$. There the solutions $\Psi^0_1(x, \varepsilon)$
      and $\Psi^0_2(x, \varepsilon)$ coincide. Now we continue analytically the solution $\Psi^0_2(x, \varepsilon)$.
      When $\Gamma_0(x, \varepsilon)$ crosses the singular direction $\arg(\beta_1-\beta_2)$ we already observe the
      Stokes phenomenon. The jump of the solution $\Psi^0_2(x, \varepsilon)$ to the solution $\Psi^0_1(x, \varepsilon)$
      is measured by the unfolding Stokes matrix $St_j(\varepsilon)$
        $$\,
          (\Psi^0_1(x, \varepsilon))^{-1}\,\Psi^0_2(x, \varepsilon)=St_j(\varepsilon)\,\hat{M}=
          St_j(\varepsilon)\,,
        \,$$         
        since $\hat{M}=I_2$.
          In this formula $St_j(\varepsilon)=St_R(\varepsilon)$ if $\arg(\beta_1-\beta_2)=0$, and
          $St_j(\varepsilon)=St_L(\varepsilon)$ if $\arg (\beta_1-\beta_2)=\pi$.
       
       In Proposition 4.31 of \cite{CL-CR} Lambert and Rousseau represent the monodromy matrices $M_j(\varepsilon)$
       of the perturbed equation as a product of the unfolded Stokes matrices $St_j(\varepsilon)$ and
       monodromy matrices of the branch of $F(x, \varepsilon)$.          
       Recently, in his remarkable paper \cite{Kl-2} Klime\v{s} provides an explicit expression of the unfolded
       Stokes matrices $St_j(\varepsilon), j=L, R$ in terms of the monodromy matrices $M_j(\varepsilon), j=L, R$
       and the matrices $e^{\pi\,i(\Lambda +B/x_j)}, j=L, R$. Here we present and use the result of 
       Klime\v{s} in our formulation from \cite{St1}. As in \cite{St1} we consider both equations provided
       that the matrices $F(x)$ and $F(x, \varepsilon)$ are not changed when we cross the first intersection
       staring from the sectors $\Omega_1$ and $\Omega_1(\varepsilon)$. In fact since the formal monodromy
       matrix $\hat{M}$ is equal to the identity matrix $I_2$, the matrices $F(x)$ and $F(x, \varepsilon)$
       are not changed even after one tour around the origin.          

      \bpr{un}          
        Let $M_j(\varepsilon)$ and $St_j(\varepsilon), j=L, R$ be the monodromy matrices and the unfolded
        Stokes matrices of the Heun type equation with respect to the fundamental solution $\Psi_1^0(x, \varepsilon)$
        on the sector $\Omega_{1, 0}(\varepsilon)$. Then depending on the $\arg (\beta_1-\beta_2)$ they 
        satisfy the following
         relations:
         \begin{enumerate}
          \item\,
          If $\arg (\beta_1-\beta_2)=0$
        on the  sector $\Omega_{1, 0}(\varepsilon)$ (upper sector) 
          $$\,
            M_L(\varepsilon)=e^{\pi\,i(\Lambda + \frac{1}{x_L}\,B)}\,St_L(\varepsilon),\quad
           M_R(\varepsilon)=St_R(\varepsilon)\,e^{\pi\,i(\Lambda + \frac{1}{x_R}\,B)}\,.
          \,$$
       On  sector $\Omega_{2, 0}(\varepsilon)$ (lower sector) 
                   $$\,
            M_L(\varepsilon)=St_L(\varepsilon)\,e^{\pi\,i(\Lambda + \frac{1}{x_L}\,B)},\quad
            M_R(\varepsilon)=e^{\pi\,i(\Lambda + \frac{1}{x_R}\,B)}\,St_R(\varepsilon)\,.
          \,$$
          \item\,
          If $\arg (\beta_1-\beta_2)=\pi$ on the sector $\Omega_{1, 0}(\varepsilon)$ (lower sector)
               $$\,
               M_L(\varepsilon)=St_L(\varepsilon)\,e^{\pi\,i(\Lambda + \frac{1}{x_L}\,B)},\quad
               M_R(\varepsilon)=e^{\pi\,i(\Lambda + \frac{1}{x_R}\,B)}\,St_R(\varepsilon)\,.
               \,$$
               On the sector $\Omega_{2, 0}(\varepsilon)$ (upper sector)
                 $$\,
                 M_L(\varepsilon)=e^{\pi\,i(\Lambda + \frac{1}{x_L}\,B)}\,St_L(\varepsilon),\quad
                 M_R(\varepsilon)=St_R(\varepsilon)\,e^{\pi\,i(\Lambda + \frac{1}{x_R}\,B)}\,.
                 \,$$ 
       \end{enumerate}
     \epr
     
     As an immediate consequence we have
       \bco{T}
         During a double resonance
         the unfolded Stokes matrices $St_j(\varepsilon)$ and the matrices 
     $e^{2 \pi\,i\,T_j}, j=L, R$ satisfy the following relation
    $$\,
    St_L(\varepsilon)=e^{2 \pi\,i\,T_L},\quad
    St_R(\varepsilon)=e^{2 \pi\,i\,T_R}\,.
    \,$$
    \eco

    \proof
    
     From \prref{un} we have that
      $$\,
        M_L(\varepsilon)=St_L(\varepsilon)\,e^{\pi\,i (\Lambda +\frac{1}{x_L}\,B)}
      \,$$
     on the sectro $\Omega_{1, 0}(\varepsilon)$ (resp. $\Omega_{2, 0}(\varepsilon)$) when $\arg (\beta_1-\beta_2)=\pi$
     (resp. $\arg (\beta_1- \beta_2)=0$).
     On the other hand during a double resonance for the monodromy matrix $M_L(\varepsilon)$ given by \eqref{M-f} we have
       $$\,
        M_L(\varepsilon)=e^{2 \pi\,i\,T_L}\,e^{\pi\,i (\Lambda +\frac{1}{x_L}\,B)}=
        e^{\pi\,i (\Lambda +\frac{1}{x_L}\,B)}\,e^{2 \pi\,i\,T_L}\,.
       \,$$ 
     Comparing both expressions for the monodromy matrix $M_L(\varepsilon)$ we find that
       $$\,
         St_L(\varepsilon)=e^{2 \pi\,i\,T_L}\,.
       \,$$
     Moreover, this relation remains valid even if
     $$\,
         M_L(\varepsilon)=e^{\pi\,i (\Lambda +\frac{1}{x_L}\,B)}\,St_L(\varepsilon)\,,
         \,$$ 
     since the matrices $e^{2 \pi\,i\,T_L}$ and $e^{\pi\,i (\Lambda + \frac{1}{x_L}\, B)}$ commite during a double
     resonance.
     
     In the same manner one can obtain the relation between $St_R(\varepsilon)$ and 
     $e^{\pi\,i (\Lambda + \frac{1}{x_R}\,B)}$.
     \qed

   %%%%%%%%%%%%%%%%%%%%%%%%%%%%%%%%%%%%%%%%%%%%%%%%%%%%%%%%%%%%%%%
   %infty
   %%%%%%%%%%%%%%%%%%%%%%%%%%%%%%%%%%%%%%%%%%%%%%%%%%%%%%%%%%%%%%%%  
      \subsection {The unfolded Stokes matrices $St_{jj}(\varepsilon)$ as the matrices $e^{2 \pi\,i\,T_{jj}}$}

      Recall that during a double resonance $x_{RR}, x_{LL}, \gamma_2-\gamma_1\in\RR$.
      So we consider the initial and the perturbed equations in a neighborhood of the real infinity point.
    Due to the symmetries of the initial and the perturbed equations we apply in a neighborhood of the real 
    infinity point the same construction as in the previus paragraph.
    Recall that the initial equation admits only one singular direction $\theta=\arg (\gamma_2-\gamma_1)$ at the
    infinity point. So,
    now we consider the initial equation together with its Stokes matrix at the infinity point
    in the ramified domain of the real infinity point 
      $\{ x\in\CC\PP^1\,: |x| > \varrho,\, \arg (\gamma_2-\gamma_1) -\kappa < \arg (x) < \arg (\gamma_2-\gamma_1) + k\}$
      where $0 < \kappa < \pi/2$. Just above we cover this domain by two open sectors
    \ben
           \Omega_{1, \infty}
           &=&
           \left\{ x=r e^{i \delta}\,|\,  r > \varrho,\,
           \arg (\gamma_2-\gamma_1) - \kappa < \delta < \arg(\gamma_2-\gamma_1) + \pi + \kappa \right\},\\
           \Omega_{2, \infty}
           &=&
           \left\{ x=r e^{i \delta}\,|\,  r > \varrho,\,
          -(\arg (\gamma_2-\gamma_1) + \pi + \kappa) < \delta < \arg(\gamma_2-\gamma_1) + \kappa \right\}\,.
           \een   
           In this paper we consider this covering as a mirror image of the corresponding covering of the origin.
           The idea is that we consider both equations on the Riemann sphere. Then if we face  this covering
           its right side will be an continuation of the real negative side.
         Then under $\arg (\gamma_2-\gamma_1)$ we mean the angle between the point $\gamma_2-\gamma_1$ and the left
         real axis in the positive sense with the base point at $x=\infty$. In particular if $\gamma_2-\gamma_1 < 0$ then  the sector $\Omega_{1, \infty}$
         is the upper sector. Conversly, if $\gamma_2-\gamma_1 > 0$ then $\Omega_{1, \infty}$ is the lower sector.
       Denote by $\Omega_{RR}$ and $\Omega_{LL}$ the connection components of the intersection
       $\Omega_{1, \infty} \cap \Omega_{2, \infty}$. 
      The radius $\varrho$ of the sectors $\Omega_j$ is so chosen that the only singular points  which
      belong to $\Omega_{RR}$ and $\Omega_{LL}$ to be $x_{RR}$ and $x_{LL}$, respectively. 
      Note that due to the mirror image, this time $\Omega_{RR}$ is the left intersection. We denote it again by
      $\Omega_{RR}$ to underline that it contains the point $x_{RR}$.

         Applying the normalization theorem of Sibuya \cite{S} to the solution $\Phi^{\infty}_{\theta}(x, 0)$ 
         from \prref{p-A-inft} we obtain over the sector $\Omega_{1, \infty}$
         the analytic matrix function $P_1(x)=P^{+}_{\theta}(x)$ which is asymptotic 
         at $x=\infty$ to the formal matrix $\hat{P}(x)$. Similarly, over the sector $\Omega_{2, \infty}$
         the matrix $P_2(x)=P^{-}_{\theta}(x)$ is an analytic matrix function asymptotic at
         $x=\infty$ to the formal matrix $\hat{P}(x)$.  Then the matrix
         $$\,
            \Psi^{\infty}_j(x)=\exp\left(-\frac{B}{x}\right)\,\left(\frac{1}{x}\right)^{-\Lambda}\,
            P_j(x)\,\exp(G x),\quad
            j=1, 2        
         \,$$
         with the corresponding branches of $(1/x)^{-\Lambda}$ and $\exp(G x)$ is an actual fundamental matrix at
         $x=\infty$ over the sector $\Omega_{j, \infty}$ respectively. If $\arg (\gamma_2-\gamma_1)=0$
         we can observe the Stokes phenomenon on $\Omega_{RR}$. If $\arg (\gamma_2-\gamma_1)=\pi$ we can
         observe the Stokes phenomenon over $\Omega_{LL}$.
      
        Let us turn around $x=\infty$ in the positive sense by analytic continuation. We start from the sector
        $\Omega_{1, \infty}$ and the solution $\Psi^{\infty}_1(x)$.On the first sector $\Omega_{jj}, j=R, L$
         that we cross we can not observe the Stokes phenomenon. On the next sector $\Omega_{jj}, j=R, L$
         that we cross we define the Stokes matrix $St_{jj}, j=R, L$ as
          $$\,
           (\Psi^{\infty}_1(x))^{-1}\,\Psi^{\infty}_2(x)=St_{jj}=St^{\theta}_{\infty}\,,
          \,$$  
   where $St^{\theta}_{\infty}$ is the Stoke matrix computed by \thref{A-S-inft}.

        At the same time we consider the perturbed equation on the whole 
          $\Omega_{1, \infty}(\varepsilon) \cup \Omega_{2, \infty}(\varepsilon)$.
        The sectorial  domains      $\Omega_{1, \infty}(\varepsilon)$ and $\Omega_{2, \infty}(\varepsilon)$
        are obtained from the open sectors $\Omega_{1, \infty}$ and $\Omega_{2, \infty}$ making a cut between
        the singular points $x_{RR}$ and $x_{LL}$ through the real axis. The point $x_0=\infty$
         belongs to this cut.
         When $\varepsilon \rightarrow 0$ the sectorial
        domains $\Omega_{j, \infty}(\varepsilon)$ tend to the sectors $\Omega_{j, \infty}, j=1, 2$, respectively.
        The sectorial domains $\Omega_{1, \infty}(\varepsilon)$ and $\Omega_{2, \infty}(\varepsilon)$ intersect in the
        sectors $\Omega_{LL}(\varepsilon)$ and $\Omega_{RR}(\varepsilon)$ and along the cut. The singular points
        $x_{jj}, j=R, L$ belong to this cut.

      The fundamental matrix $\Phi_{\infty}(x, \varepsilon)$ of the Heun type equation can be rewritten in the form
       $$\,
         \Phi_{\infty}(x, \varepsilon)=F(x, \varepsilon)\,P(x, \varepsilon)\,G(x, \varepsilon)\,,
         \,$$     
  where $G(x, \varepsilon)$ and $F(x, \varepsilon)$ are defined as above.
    The matrix 
$P(x, \varepsilon)=\{p_{ij}(x, \varepsilon)\}_{ij=1}^2$ is given by
   \ben
      & &
      p_{ii}=1, \quad p_{ij}(x, \varepsilon)=0\quad \textrm{for}\quad i > j,\\[0.15ex]
         & &
         p_{12}(x, \varepsilon)=
         (x-x_{LL})^{\frac{\gamma_1-\gamma_2}{2 \sqrt{\varepsilon}}}
         (x_{RR}-x)^{-\frac{\gamma_1-\gamma_2}{2 \sqrt{\varepsilon}}}
         \int_{\Gamma_{\infty}(x, \varepsilon)}
         \frac{\Phi_2(z, \varepsilon)}{\Phi_1(z, \varepsilon)}\, d z\,.
   \een
   The path $\Gamma_{\infty}(x, \varepsilon)$ is a path either from $x_{LL}$ or $x_{RR}$ to $x$ taken
   in the direction $\arg (\gamma_1-\gamma_2)$. When we continue analytically the path $\Gamma_{\infty}(x, \varepsilon)$ 
   by turning around $x=\infty$ in the positive sense, we obtain two branches $p_{12}^{-}(x, \varepsilon)$
   and $p_{12}^{+}(x, \varepsilon)$ of the element $p_{12}(x, \varepsilon)$. The branch $p_{12}^{-}(x, \varepsilon)$
   corresponds to the path taken in the direction $\arg (\gamma_2-\gamma_1)-\epsilon$, and the branch
   $p_{12}^{+}(x, \varepsilon)$ corresponds to the path taken in the direction $\arg (\gamma_2-\gamma_1)+\epsilon$.
   Here $\epsilon > 0$ is a small number. 
   We redefine over the sector $\Omega_{1, \infty}(\varepsilon)$ the fundamental matrix $\Phi_{\infty}(x, \varepsilon)$
   as 
   $$\,
     \Psi^{\infty}_1(x, \varepsilon)=F(x, \varepsilon)\,P_1(x, \varepsilon)\,G(x, \varepsilon)\,,
   \,$$
   where $P_1(x, \varepsilon)=\{p^{+}_{ij}(x, \varepsilon)\}_{i, j=1}^2$. Similarly, over the sector
   $\Omega_{2, \infty}(\varepsilon)$ we redefine the fundamental matrix $\Phi_{\infty}(x, \varepsilon)$ as
   $$\,
   \Psi^{\infty}_2(x, \varepsilon)=F(x, \varepsilon)\,P_2(x, \varepsilon)\,G(x, \varepsilon)\,,
   \,$$
   where $P_2(x, \varepsilon)=\{h^{-}_{ij}(x, \varepsilon)\}_{i, j=1}^2$. 
  Let us turn around $x=\infty$ in the positive sense starting from the sector $\Omega_{1, \infty}(\varepsilon)$
  and the solution $\Psi^{\infty}_1(x, \varepsilon)$ on it. When the path $\Gamma_{\infty}(x, \varepsilon)$ crosses
  the direction $\arg (\gamma_1-\gamma_2)$ we can not observe the Stokes phenomenon. On the sector
  $\Omega_{jj}(\varepsilon)$ that is crossed by the direction $\arg (\gamma_1-\gamma_2)$ the solutions
  $\Psi^{\infty}_1(x, \varepsilon)$ and $\Psi^{\infty}_2(x, \varepsilon)$ coincide. When $\Gamma_{\infty}(x, \varepsilon)$
  crosses the singular direction $\theta=\arg (\gamma_2-\gamma_1)$ the solution $\Psi^{\infty}_2(x, \varepsilon)$
  jumps to the solution $\Psi^{\infty}_1(x, \varepsilon)$
    $$\,
      (\Psi^{\infty}_1(x, \varepsilon))^{-1}\,\Psi^{\infty}_2(x, \varepsilon)=St_{jj}(\varepsilon)\,.
    \,$$ 
This jump is measured geometrically by the unfolded Stokes matrix $St_{jj}(\varepsilon)$.
   If $\arg (\gamma_2-\gamma_1)=0$ we have that $j=R$, and if $\arg (\gamma_2-\gamma_1)=\pi$ then
   $j=L$. Due to the symmetries of the perturbed equation we can
  replace the origin by the infinity point in \prref{un}. Then the analog of \prref{un} gives us an explicit connection between
  the monodromy matrices $M_{jj}(\varepsilon)$, the unfolded Stokes matrices $St_{jj}(\varepsilon)$ and
  the matrices $e^{-\pi\,i\,x_{jj}\,G}$

       \bpr{un-inf}          
       Let $M_{jj}(\varepsilon)$ and $St_{jj}(\varepsilon), j=L, R$ be the monodromy matrices and the unfolded
       Stokes matrices of the Heun type equation with respect to the fundamental solution on the 
       sector $\Omega_{1, \infty}(\varepsilon)$. Then depending on the $\arg (\gamma_2-\gamma_1)$ they satisfy the following
       relations:
       \begin{enumerate}
       	\item\,
       	If $\arg (\gamma_2-\gamma_1)=0$
       	on the  sector $\Omega_{1, \infty}(\varepsilon)$ (lower sector) 
       	$$\,
       	M_{LL}(\varepsilon)=e^{-\pi\,i\,x_{LL}\,G}\,St_{LL}(\varepsilon),\quad
       	M_{RR}(\varepsilon)=St_{RR}(\varepsilon)\,e^{-\pi\,i\,x_{RR}\,G}\,.
       	\,$$
       	On  sector $\Omega_{2, \infty}(\varepsilon)$ (upper sector) 
       	$$\,
       	M_{LL}(\varepsilon)=St_{LL}(\varepsilon)\,e^{-\pi\,i\,x_{LL}\,G},\quad
       	M_{RR}(\varepsilon)=e^{-\pi\,i\,x_{RR}\,G}\,St_{RR}(\varepsilon)\,.
       	\,$$
       	\item\,
       	If $\arg (\gamma_2-\gamma_1)=\pi$ on the sector $\Omega_{1, \infty}(\varepsilon)$ (upper sector)
       	$$\,
       	M_{LL}(\varepsilon)=St_{LL}(\varepsilon)\,e^{-\pi\,i\,x_{LL}\,G},\quad
       	M_{RR}(\varepsilon)=e^{-\pi\,i\,x_{RR}\,G}\,St_{RR}(\varepsilon)\,.
       	\,$$
       	On the sector $\Omega_{2, \infty}(\varepsilon)$ (lower sector)
       	$$\,
       	M_{LL}(\varepsilon)=e^{-\pi\,i\,x_{LL}\,G}\,St_{LL}(\varepsilon),\quad
       	M_{RR}(\varepsilon)=St_{RR}(\varepsilon)\,e^{-\pi\,i\,x_{RR}\,G}\,.
       	\,$$ 
       \end{enumerate}
       \epr
       
          As an immediate consequence we have
          \bco{T-inf}
          During a double resonance
          the unfolded Stokes matrices $St_{jj}(\varepsilon)$ and the matrices 
          $e^{2 \pi\,i\,T_{jj}}, j=L, R$ satisfy the following relation
          $$\,
          St_{LL}(\varepsilon)=e^{2 \pi\,i\,T_{LL}},\quad
          St_{RR}(\varepsilon)=e^{2 \pi\,i\,T_{RR}}\,.
          \,$$
          \eco
       
        \proof
        The proof is similar to the proof of \coref{T}.\qed
%%%%%%%%%%%%%%%%%%%%%%%%%%%%%%%%%%%%%%%%
%main result
%%%%%%%%%%%%%%%%%%%%%%%%%%%%%%%%%%%%

   \subsection{The main results}
   It turns out that the matrices $e^{2 \pi\,i\,T_j}$ and $e^{2 \pi\,i\,T_{jj}}$ have limits when
   $\sqrt{\varepsilon} \rightarrow 0$ radially. In the following lemma 
   we  compute the limits of the numbers $d_j$ and $d_{jj}$ given in the previous section.

	\ble{lim-A}
	 Assume that $\beta_1 \neq \beta_2$. Then during a double resonance
	  the number $d_j$ either is equal to zero or has the following limit when
	 $\sqrt{\varepsilon} \rightarrow 0$
	  \ben
		  \lim_{\sqrt{\varepsilon} \rightarrow 0} d_j
                               =
				(\gamma_2-\gamma_1)
			\sum_{k=0}^{\infty}
			\frac{(-1)^k\,(\gamma_2-\gamma_1)^k\,(\beta_2-\beta_1)^k}{k!\,(k+1)!}\,.
     \een
     Similarly,  during a double resonance
     the number $d_{jj}$ either is equal to zero or has the following limit when
     $\sqrt{\varepsilon} \rightarrow 0$
        \ben 			
               \lim_{\sqrt{\varepsilon} \rightarrow 0} d_{jj}
                               =
				-(\gamma_2-\gamma_1)
			\sum_{k=0}^{\infty}
			\frac{(-1)^k\,(\gamma_2-\gamma_1)^k\,(\beta_2-\beta_1)^k}{k!\,(k+1)!}.
		\een
	\ele

  \proof
	
	We will study in detail only the limit of the numbers $d_L$ and $d_{LL}$ computed 
	in \thref{A.1-M} and \thref{A.1-M-inf}. The limit of the rest numbers $d_j$ and $d_{jj}$ obtained in the
	previous section is computed in the same manner.
	 
	Let us firstly compute the limit of the finite sum $A$ when $\sqrt{\varepsilon} \rightarrow 0$ introduced
	in \thref{A.1-M}.
	We can rewrite $A$ as
	 $$\,
	  A=\sum_{s=0}^k 
		   \left(\begin{array}{c}
			     k\\
					 s
					   \end{array}
				\right)
		\frac{\Gamma(\frac{\gamma_1-\gamma_2}{2 \sqrt{\varepsilon}}+s)}
				 {\Gamma(\frac{\gamma_1-\gamma_2}{2 \sqrt{\varepsilon}})\,
				 \left(\frac{\gamma_1-\gamma_2}{2 \sqrt{\varepsilon}}\right)^s}
		\frac{\Gamma(\frac{\gamma_1-\gamma_2}{2 \sqrt{\varepsilon}})
		      \left(\frac{\gamma_1-\gamma_2}{2 \sqrt{\varepsilon}}\right)^{1-k+s}}
				 {\Gamma(\frac{\gamma_1-\gamma_2}{2 \sqrt{\varepsilon}}+1-k+s)}
		\frac{\left(\frac{\gamma_1-\gamma_2}{2 \sqrt{\varepsilon}}\right)^s}
			   {\left(\frac{\gamma_1-\gamma_2}{2 \sqrt{\varepsilon}}\right)^{1-k+s}}
			   \left(\frac{1+\varepsilon}{1-\varepsilon}\right)^s
			   \,.
	\,$$
  Then using the limit
	 \be\label{G}
	   \lim_{|z| \rightarrow \infty}
		 \frac{\Gamma(z+\alpha)}{\Gamma(z)\,z^{\alpha}}=1\,,
	\ee
	we find that
	 $$\,
	  \lim_{\sqrt{\varepsilon} \rightarrow 0} A=
	  \left(\frac{2}{1-\varepsilon}\right)^k
			\left(\frac{\gamma_1-\gamma_2}{2 \sqrt{\varepsilon}}\right)^{k-1}\,.
		\,$$
 As a result the number $d_L$ is reduced to the form
  \ben
	  d_L=(\gamma_2-\gamma_1)\left(\frac{1+\varepsilon}{1-\varepsilon}\right)^
		   {\frac{\gamma_1-\gamma_2}{2 \sqrt{\varepsilon}}}
			\sum_{k=1}^{\frac{\beta_2-\beta_1}{2 \sqrt{\varepsilon}}}
			\frac{(2 \sqrt{\varepsilon})^{k-1}}{(k-1)! k!}
			\left(\frac{\sqrt{\varepsilon}}{1+\varepsilon}\right)^k
			\left(\frac{2}{1-\varepsilon}\right)^k
			\left(\frac{\gamma_2-\gamma_1}{2 \sqrt{\varepsilon}}\right)^k
			\frac{\Gamma(\frac{\beta_2-\beta_1}{2 \sqrt{\varepsilon}})}
			     {\Gamma(\frac{\beta_2-\beta_1}{2 \sqrt{\varepsilon}}-k+1)}\,.
	\een
 Now, choosing $\log 1=0$, it is easy to see that
  \ben
	  \lim_{\sqrt{\varepsilon} \rightarrow 0}
		  \left(\frac{1+\varepsilon}{1-\varepsilon}\right)^
			{\frac{\gamma_1-\gamma_2}{2 \sqrt{\varepsilon}}}=1\quad
			\textrm{and}\quad
	 \lim_{\sqrt{\varepsilon} \rightarrow 0} (1-\varepsilon^2)^{-k}=1\,.		
	\een
Dropping out these two expressions, we rewriting $d_L$ as
\ben
  d_L=(\gamma_2-\gamma_1)
  \sum_{k=1}^{\frac{\beta_2-\beta_1}{2 \sqrt{\varepsilon}}}
	\frac{(2 \sqrt{\varepsilon})^{k-1} (2 \sqrt{\varepsilon})^k}
	     {(k_1)! k!}
	     \left(\frac{\gamma_1-\gamma_2}{2 \sqrt{\varepsilon}}\right)^k 
	     \left(\frac{\beta_2-\beta_1}{2 \sqrt{\varepsilon}}\right)^{k-1}
				\frac{\Gamma(\frac{\beta_2-\beta_1}{2 \sqrt{\varepsilon}}) 
				\left(\frac{\beta_2-\beta_1}{2 \sqrt{\varepsilon}}\right)^{-k+1}}
			     {\Gamma(\frac{\beta_2-\beta_1}{2 \sqrt{\varepsilon}}-k+1)}\,.
\een
  Again applying the limit \eqref{G}, we find that
	  $$\,
		  \lim_{\sqrt{\varepsilon} \rightarrow 0} d_L=
			(\gamma_2-\gamma_1)
			\sum_{k=0}^{\infty}
			\frac{(-1)^k\,(\gamma_2-\gamma_1)^k\,(\beta_2-\beta_1)^k}{k!\,(k+1)!}\,.
		\,$$
		
		In the same manner we can rewrite the finite sum $A$  from \thref{A.1-M-inf} as 
		\ben
		A=\sum_{s=0}^k \left(\begin{array}{c}
			             k\\
			             s\end{array}\right)
		\frac{\Gamma(\frac{\beta_2-\beta_1}{2 \sqrt{\varepsilon}}+1+s)}
		     {\Gamma(\frac{\beta_2-\beta_1}{2 \sqrt{\varepsilon}})\,
		     	\left(\frac{\beta_2-\beta_1}{2 \sqrt{\varepsilon}}\right)^{1+s}}
		\frac{\Gamma(\frac{\beta_2-\beta_1}{2 \sqrt{\varepsilon}})\,
			\left(\frac{\beta_2-\beta_1}{2 \sqrt{\varepsilon}}\right)^{1+s}}
		     {\Gamma(\frac{\beta_2-\beta_1}{2 \sqrt{\varepsilon}}-k+s)}
		\frac{\left(\frac{\beta_2-\beta_1}{2 \sqrt{\varepsilon}}\right)^{-k+s}}
		     {\left(\frac{\beta_2-\beta_1}{2 \sqrt{\varepsilon}}\right)^{-k+s}}
		     \left(\frac{1+\varepsilon}{1-\varepsilon}\right)^s\,.          
		\een
		Again using the limit \eqref{G} we find that
		$$\,
		  \lim_{\sqrt{\varepsilon} \rightarrow 0}\,A=
		  \left(\frac{\beta_2-\beta_1}{2 \sqrt{\varepsilon}}\right)^{k+1}\,
		  \left(\frac{2}{1-\varepsilon}\right)^k\,.
		\,$$
		As a result the number $d_{LL}$ is reduced to
		\ben
		  d_{LL}=\frac{1}{1-\varepsilon^2} 
		  \left(\frac{1+\varepsilon}{1-\varepsilon}\right)^{\frac{\beta_2-\beta_1}{2 \sqrt{\varepsilon}}}
		  \sum_{k=0}^{\frac{\gamma_1-\gamma_2}{2 \sqrt{\varepsilon}}}
		  \frac{1}{k!\,(k+1)!}
		  \frac{(2 \sqrt{\varepsilon})^{1+k}\,(\sqrt{\varepsilon})^k}
		      {(1+\varepsilon)^k}
		  \left(\frac{\beta_2-\beta_1}{2 \sqrt{\varepsilon}}\right)^k\,
		  \left(\frac{2}{1-\varepsilon}\right)^k
		  \frac{\Gamma(\frac{\gamma_1-\gamma_2}{2 \sqrt{\varepsilon}}+1)}
		       {\Gamma(\frac{\gamma_1-\gamma_2}{2 \sqrt{\varepsilon}}-k)}    
		\een
		Now it is not difficult to show that
		$$\,
		  \lim_{\sqrt{\varepsilon} \rightarrow 0}\,d_{LL}=
		  (\gamma_1-\gamma_2)
		  \sum_{k=0}^{\infty}
		    \frac{(-1)^k\,(\gamma_2-\gamma_1)^k\,(\beta_2-\beta_1)^k}{k!\,(k+1)!}\,.
		\,$$
		
 This completes the proof.\qed

   In Theorem 4.25 in \cite{CL-CR1} Lambert and Rousseau prove that the unfolded Stokes matrices 
   $St_j(\varepsilon), j=R, L$ depend analytically on the parameter of perturbation $\varepsilon$
   and converge when $\varepsilon\rightarrow 0$ to the Stokes matrices $St_j, j=R, L$ of the initial equation.
   Due to the symmetries we have the same result for the unfolded Stokes matrices $St_{jj}, j=R, L$
   and the Stokes matrices $St_{jj}, j=R, L$ of the initial equation. 
  Then thanks to \leref{lim-A}, \coref{T}, \coref{T-inf} and Theorem 4.25 in \cite{CL-CR1}
     we state the main result of this paper
  
   \bth{main}
    Assume  that $\varepsilon\in\RR^{+}$.
    Assume also that $\beta_2-\beta_1, \gamma_2-\gamma_1\in\RR$ are fixed such that $\beta_1 \neq \beta_2$
    and $\frac{\beta_2-\beta_1}{2 \sqrt{\varepsilon}},\,\frac{\gamma_2-\gamma_1}{2 \sqrt{\varepsilon}}\in\ZZ$.
     Then depending on the position of $\beta_2-\beta_1$ and $\gamma_2-\gamma_1$ toward 0 
     the matrices $T_j, T_{jj},\,j=R, L$ of the Heun type equation and the Stokes matrices 
     $St^{\theta}_k,\,k=0, \infty$ of the initial equation are connected as follows,
     \begin{enumerate}
     	\item\,
     	If $\beta_2-\beta_1, \gamma_2-\gamma_1\in\RR^{+}$ then
     	 $$\,
     	    e^{2 \pi\,i\,T_L} \longrightarrow St^{\pi}_0,\quad
     	      e^{2 \pi\,i\,T_{RR}} \longrightarrow St^{0}_{\infty}\,,
     	 \,$$
     	 when $\sqrt{\varepsilon} \rightarrow 0$.
     	 	\item\,
     	 	If $\beta_1-\beta_2, \gamma_1-\gamma_2\in\RR^{+}$ then
     	 	$$\,
     	 	e^{2 \pi\,i\,T_R} \longrightarrow St^0_0,\quad
     	 	e^{2 \pi\,i\,T_{LL}} \longrightarrow St^{\pi}_{\infty}\,,
     	 	\,$$
     	 	when $\sqrt{\varepsilon} \rightarrow 0$.
     	 		\item\,
     	 		If $\beta_2-\beta_1, \gamma_1-\gamma_2\in\RR^{+}$ then
     	 		$$\,
     	 		e^{2 \pi\,i\,T_L} \longrightarrow St^{\pi}_0,\quad
     	 		e^{2 \pi\,i\,T_{LL}} \longrightarrow St^{\pi}_{\infty}\,,
     	 		\,$$
      		when $\sqrt{\varepsilon} \rightarrow 0$.
         	\item\,
         	If $\beta_1-\beta_2, \gamma_2-\gamma_1\in\RR^{+}$ then
         	$$\,
         	e^{2 \pi\,i\,T_R} \longrightarrow St^0_0,\quad
         	e^{2 \pi\,i\,T_{RR}} \longrightarrow St^{0}_{\infty}\,,
         	\,$$
         	when $\sqrt{\varepsilon} \rightarrow 0$.
     \end{enumerate}
   \ethe
  
    \vspace{3ex}

    {\bf Funding Information}\, The author was partially supported by Grant DN 02-5/2016 of the Bulgarian Fond
    "Scientific Research".

	%%%%%%%%%%%%% References %%%%%%%%%%%%%%%%%%%%%%%%%%%%%%%%%%%%%%%%%%%%%%%
\begin{small}
    
\end{small}
%%%%%%%%%%%%%%%%%%%%%%%%%%%%%%%%%%%%%%%%%%%%%%%%%%%%%%%%%%%%%%%

\end{document}